\documentclass{article}
\usepackage{amsmath,amssymb,amsthm, amsfonts}
\usepackage{etaremune, enumerate, float, verbatim}
\usepackage{algorithm}
\usepackage{algorithmic}
\usepackage{url}
\usepackage{hyperref}
\usepackage{graphicx}
\usepackage{url}
\usepackage[mathlines]{lineno}
\usepackage{dsfont} 
\usepackage{mathrsfs} 
\usepackage{soul} 
\usepackage{color}

\def\blu{\color{blue}}
\definecolor{green}{rgb}{0,0.7,0}

\usepackage{tikz}

\numberwithin{figure}{section}   
\numberwithin{table}{section}
\numberwithin{equation}{section} 

\newtheorem{thm}{Theorem}[section]
\newtheorem{cor}[thm]{Corollary}
\newtheorem{lem}[thm]{Lemma}
\newtheorem{prop}[thm]{Proposition}

\newtheorem{obs}[thm]{Observation}
 \newtheorem{alg}[thm]{Algorithm}

\theoremstyle{definition}
\newtheorem{rem}[thm]{Remark}

\theoremstyle{definition}
\newtheorem{defn}[thm]{Definition}

\theoremstyle{definition}
\newtheorem{ex}[thm]{Example}

\newcommand {\du} {\,\sqcup\,}

\newcommand{\bit}{\begin{itemize}}
\newcommand{\eit}{\end{itemize}}
\newcommand{\ben}{\begin{enumerate}}
\newcommand{\een}{\end{enumerate}}
\newcommand{\beq}{\begin{equation}}
\newcommand{\eeq}{\end{equation}}
\newcommand{\bea}{\begin{eqnarray*}}
\newcommand{\eea}{\end{eqnarray*}}
\newcommand{\bean}{\begin{eqnarray}}
\newcommand{\eean}{\end{eqnarray}}
\newcommand{\bpf}{\begin{proof}}
\newcommand{\epf}{\end{proof}\ms}
\newcommand{\bmt}{\begin{bmatrix}}
\newcommand{\emt}{\end{bmatrix}}
\newcommand{\ms}{\medskip}
\newcommand{\noi}{\noindent}
\newcommand{\beqa}{\begin{array}}
\newcommand{\eeqa}{\end{array}}
\newcommand{\ol}{\overline}

\newcommand{\lf}{\left\lfloor}
\newcommand{\rf}{\right\rfloor}
\newcommand{\lp}{\left(}
\newcommand{\rp}{\right)}

\newcommand{\Xb}{X}
\newcommand{\Xc}{X}
\newcommand{\Xo}{X}
\newcommand{\Xr}{X}

\newcommand{\Z}{\operatorname{Z}}
\newcommand{\Zs}{\operatorname{Z}_-}
\newcommand{\Zp}{\operatorname{Z}_+}

\newcommand{\Dc}{\gamma_c}
\newcommand{\pd}{\gamma_p} 
\newcommand{\vc}{\tau}
\newcommand{\zbar}{\ol\Z}
\newcommand{\zsbar}{\ol\Z_-}
\newcommand{\zpbar}{\ol\Z_+}

\newcommand{\dbar}{\ol\gamma}
\newcommand{\pdbar}{\ol\pd}
\newcommand{\vcbar}{\ol\vc}
\newcommand{\ibar}{\underline{\alpha}}

\newcommand{\xbar}{\ol X}

\newcommand{\xbbar}{\ol X}

\newcommand{\xtar}{\mathfrak{X}}
\newcommand{\xbtar}{\mathfrak{X}}
\newcommand{\xrtar}{\mathfrak{X}}
\newcommand{\X}{\mathfrak{X}}
\newcommand{\ztar}{\mathfrak{Z}}
\newcommand{\zstar}{\mathfrak{Z}^-}
\newcommand{\zptar}{\mathfrak{Z}^+}
\newcommand{\dtar}{\mathfrak{D}}
\newcommand{\dctar}{\mathfrak{D}^c}
\newcommand{\pdtar}{\mathfrak{P}}
\newcommand{\vctar}{\mathfrak{C}}
\newcommand{\itar}{\mathfrak{I}}
\newcommand{\irtar}{\mathfrak{Ir}}
\newcommand{\zirtar}{\mathfrak{Zir}}

\newcommand{\xytar}{{\mathfrak X}_Y}
\newcommand{\yxtar}{{\mathfrak Y}_X}

\newcommand{\zso}{z^-_0}
\newcommand{\zpo}{z^+_0}

\newcommand{\xxo}{x_0}
\newcommand{\xbo}{x_0}

\newcommand{\ddo}{d_0}
\newcommand{\zzo}{z_0}
\newcommand{\pdo}{pd_0}
\newcommand{\vco}{\tau_0}
\newcommand{\io}{\alpha_0}

\newcommand{\ulzo}{\underline{z_0}}
\newcommand{\ulxo}{\underline{x_0}}
\newcommand{\ulxbo}{\underline{x_0}}

\newcommand{\ulzso}{\underline{z^-_0}}
\newcommand{\ulzpo}{\underline{z^+_0}}
\newcommand{\ulvco}{\underline{\tau_0}}
\newcommand{\ulio}{\underline{\alpha_0}}

\newcommand{\ulpdo}{\underline{pd_0}}
\newcommand{\uldo}{\underline{d_0}}

\newcommand{\dist}{\operatorname{dist}}
\newcommand{\diam}{\operatorname{diam}}

\DeclareMathOperator{\ir}{ir}
\DeclareMathOperator{\IR}{IR}

\DeclareMathOperator{\zir}{zir}
\DeclareMathOperator{\ZIR}{ZIR}

\newcommand{\ytar}{\mathfrak{Y}}
\newcommand{\yyo}{y_0}
\newcommand{\ulyo}{\overline{y_0}}
\newcommand{\ybar}{\underline{Y}}

\title{TAR reconfiguration for vertex set parameters}
\author{Bryan Curtis\thanks{Department of Mathematics, Iowa State University,
Ames, IA 50011, USA (bcurtis1@iastate.edu)}\and
Mary Flagg\thanks{Department of Mathematics, Statistics and Computer Science, University of St. Thomas, Houston, TX 77006, USA (flaggm@stthom.edu).} \and
   Leslie
Hogben\thanks{American Institute of Mathematics, Pasadena, CA 91125, USA
(hogben@aimath.org) and Department of Mathematics, Iowa State University,
Ames, IA 50011, USA.}}

\begin{document}

\maketitle 

\begin{abstract}  This paper  surveys results about token addition and removal (TAR) reconfiguration for several  well-known vertex set parameters including domination, power domination, standard zero forcing, and PSD zero forcing. We also expand the range of parameters to which universal $X$-set TAR graph results apply,  for  $X$-sets and their TAR graphs introduced in  [B.~Bjorkman, C.~Bozeman, D.~Ferrero, M.~Flagg, C.~Grood, L.~Hogben, B.~Jacob, C.~Reinhart,  Power domination reconfiguration,  \url{https://arxiv.org/abs/2201.01798}]  and 
 [N.H.~Bong, J.~Carlson, B.~Curtis, R.~Haas, L.~Hogben, Isomorphisms  and  properties of TAR reconfiguration graphs for zero forcing and other $X$-set parameters,  \emph{Graphs Combin.} {39} (2023), Paper No. 86].   Here we examine which of the $X$-set axioms are needed for which results.  With this new framework, the main  results apply to skew zero forcing and vertex covering, and results about TAR reconfiguration graphs of these parameters are presented.  While $X$-sets are defined for parameters that take the minimum cardinality over the $X$-sets  of a graph, and $X$-set results are restricted to such minimizing parameters, our expansion of the universal perspective allows these results to be applied to parameters that take the maximum value among relevant sets, called $Y$-sets.  Maximizing parameters to which the main results apply include independence number, (upper) irredundance number, and (upper) zero forcing irredundance number; 
 TAR reconfiguration results are presented for these parameters.  We also show that the equivalence of connectedness in certain token jumping reconfiguration graphs and certain TAR reconfiguration graphs for independent sets established in [M.~Kami\'nski, P.~Medvedev, M.~Milani\v c. Complexity of independent set reconfiguration problems. {\em J. Theoretical Computer Science} 439 (2012), 9--15.] extends to $X$-set and $Y$-set parameters.
\end{abstract}

 \noi {\bf Keywords} reconfiguration;  token addition and removal; TAR; vertex set; $X$-set;  $Y$-set; domination; power domination; zero forcing; PSD forcing; skew forcing; vertex cover; independence number; irredundance  

\noi{\bf AMS subject classification} 05C69, 05C45, 05C50, 05C57,  05C60, 05C70, 68R10

\section{Introduction}

The study of reconfiguration examines relationships among solutions to a problem. The \emph{reconfiguration graph} for the problem has as its vertices solutions to a problem and edges are determined by a \emph{reconfiguration rule} that describes  relationships between the solutions. The reconfiguration rule can be viewed as describing a single step in the process of transforming  one solution  
 to another in which each intermediate state is also a solution.  The ability to transform one solution to another solution is equivalent to having a path between the two solutions in the reconfiguration graph, i.e., the two solutions are in the same connected component of the reconfiguration graph.  
In \cite{Nish18}, Nishimura  surveys recent work on structural and algorithmic (complexity) questions   across a broad variety of parameters. She describes three types of reconfiguration rules, including  token addition and removal.  

For many reconfiguration problems arising from  graphs (including all those we study here), a solution can be represented as a subset of the vertices of a graph $G$; we call such a problem a \emph{vertex set problem}. The \emph{token addition and removal (TAR) reconfiguration graph} for a vertex set problem has an edge between two sets  if and only if one can be obtained from the other by the addition or removal of a single vertex.   Both the TAR graph (which includes all feasible subsets of vertices of $G$) and the $k$-TAR graph, which allows only feasible subsets of at most  (or at least) $k$ vertices, are studied. 

A set $S\subseteq V(G)$ is a \emph{dominating  set} of a graph $G$ if every vertex of $G$ is in $S$ or a neighbor of a vertex in $S$.  The \emph{domination number}  of $G$, denoted by $\gamma(G)$, is the minimum cardinality of a dominating set of $G$.  The domination number has been widely studied; see for example, the book by Haynes, Hedetniemi, and Slater \cite{HHS98} and more recent books by Haynes, Hedetniemi, and Henning, such as \cite{HHH23}.  Haas and Seyffarth initiated the study of reconfiguration graphs for domination in \cite{HS14}. 
 Many connectedness results for $X$-set parameters in \cite{PDrecon} use the same ideas as  the analogous results for domination.

A universal definition of an $X$-set for a parameter  determined by subsets of vertices was given in \cite{PDrecon} and \cite{XTARiso}; parameters to which this definition applies are here called original $X$-set parameters and include the domination number, the power domination number,  the standard zero forcing number,   and the PSD zero forcing number, but not to the skew forcing number (see Section \ref{ss:param} for definitions of these parameters). 
Bjorkman et al.\  were the first to use the universal approach with $X$-sets. In addition to establishing  connectedness results,  they showed that  for an on original $X$-set parameter and  a graph $G$ with no isolated vertices, the order of $G$ can be determined from its $X$-TAR graph.   Power domination TAR graphs were  also studied  in detail in that paper (as were token jumping reconfiguration graphs, also known as token exchange).  In \cite{XTARiso} numerous additional results about $X$-TAR graphs were established, including that  for an on original $X$-set parameter and  graphs $G$ and $G'$ with no isolated vertices, an isomorphism of the $\Xo$-TAR graphs can be used to find a relabeling of the vertices of $G'$ so that $G$ and $G'$ have the same $X$-sets; standard zero forcing  TAR graphs were also studied  in that paper. 

Skew zero  forcing is an obvious omission from the list of parameters to which the results in \cite{PDrecon} and \cite{XTARiso} apply, because skew zero forcing allows the empty set to be a skew forcing set, and this was prohibited by one of the rules for $X$-sets in \cite{PDrecon,XTARiso} (see Definition \ref{def:original+assumed}).
In Section \ref{sec:proof-requirements} we determine which of the  axioms in the original definition of $X$-set are needed for each result of \cite{PDrecon} and \cite{XTARiso}. 
We present both a   very general definition of a  \emph{super $X$-set parameter}  (Definition \ref{d:superX}), which uses only one of the original $X$-set parameter axioms, and a stronger  \emph{robust $X$-set parameter} definition   (Definition \ref{Robust-parameters-2}),  which uses versions of  three of the  five axioms for original $X$-set parameters.  Many connectedness results  stated in \cite{PDrecon} (and derived from earlier work such as \cite{HS14}) and a result from \cite{KMM12} about token jumping graphs and TAR graphs  are established for every super $X$-set parameter, but  the isomorphism results need not be true for such a parameter.
Almost all results of \cite{PDrecon,XTARiso} are true for every robust $X$-set parameter, allowing the extensions of these results to  additional parameters including skew forcing number and vertex cover number. 

In Section \ref{s:Dom-PD-Z}, we survey known results about TAR graphs for various parameters  related to domination and zero forcing,  primarily about connectedness of TAR subgraphs, isomorphism properties, uniqueness, and existence of Hamilton cycles or paths.  Parameters that have been studied include domination number, power domination number, and standard zero forcing number; we also present a small additional amount of material for these parameters. 
TAR graphs for PSD forcing and  skew forcing   are examined here in more detail in Sections  \ref{s:PSD} and  \ref{s:skew},    where examples illustrating connectedness, uniqueness, and  Hamiltonicity are presented, in addition to the main order and isomorphism results and connectedness bounds  that follow from the revised $X$-set parameter  definitions.    Reconfiguration of vertex covers was studied in \cite{IDHPSUU11,INZ15} and other works.  Complexity is the focus of much of this prior work on vertex cover reconfiguration, but some structural results have been obtained and we  provide a partial survey of these results in Section \ref{s:vtxcover}.   We do not survey complexity results for any of the parameters discussed.  

We  initially take a very general approach in Section \ref{sec:proof-requirements}  by defining a vertex set parameter, which allows the main results to be extended to many vertex set parameters that take the maximum size of a set rather than the minimum; such parameters are here called $Y$-set parameters.  General results for TAR graphs of $Y$-set parameters are discussed in Section \ref{s:maxYparam}, together with a more detailed discussion on TAR graphs for independence, irredunundance,  zero forcing irredundance. A complementation technique is introduced to  translate $X$-set parameter results to $Y$-set parameter results.

In this arXiv version we have included the proofs of some results from other papers when the result is  not stated as generally as it is here but essentially the same proof works; these are colored {\blu blue}. In the next section we list notation we will use throughout. Section \ref{ss:param} provides definitions of the parameters that are discussed.

 \subsection{Notation and terminology}\label{ss:term}

  Every graph discussed is  simple, undirected, finite, and has a nonempty vertex set. A graph $G=(V(G),E(G))$ consists of a set of vertices $V(G)$ and a set of edges $E(G)$ which are two element subsets of $V(G)$. The edge $\{v,w\} \in E(G)$ is denoted $vw$. The order of the graph is $|V(G)|$. A graph is \emph{odd} or \emph{even} according as its order is odd or even. The subgraph of $G$ induced by a set $W$ of vertices is denoted by $G[W]$.
 For $S\subseteq V(G)$, define $G - S=G[V(G)\setminus S]$. When $S$ is a single vertex $v\in V(G)$, the notation $G - \{v\}$ is simplified to $G - v$.

If $uw \in E(G)$, then $v$ and $w$ are said to be \emph{adjacent} or \emph{neighbors}.  The \emph{neighborhood} of $v$ is the set of neighbors of $v$ and is denoted by $N_G(v)$ and the \emph{closed neighborhood} of $v$ is $N_G[v]=N_G(v)\cup\{v\}$. The closed neighborhood of a set $S\subseteq V(G)$ is  $N_G[S]=\cup_{x\in S}N_G[x]$.  The subscript $G$ may be omitted when it is clear from context. The degree of the vertex $v$ is  $\deg_G(v)=|N(v)|$. For a graph $G$, $\Delta(G)=\max\{\deg_G(v):v \in V(G)\}$ and $\delta(G)=\min\{\deg_G(v):v \in V(G)\}$.

Given distinct vertices $v_0$ and $v_{\ell}$ in $V(G)$,  a \emph{path} of length $\ell$  from $v_0$ to $v_{\ell}$ is a sequence of distinct vertices $(v_0,\ldots, v_{\ell})$ such that 
$v_i$ is a neighbor of $v_{i+1}$ for every integer $i$, $0\leq i\leq \ell -1$. A \emph{cycle} on $k$ vertices, denoted $(c_1, \dots, c_k)$, is a sequence of distinct vertices $c_1, \dots, c_k$ with an edge between $c_i$ and $c_{i+1}$ for all $1 \leq i \leq k-1$ and the edge $c_1c_k$.   A \emph{Hamilton cycle} (respectively,  \emph{Hamilton path})  of a graph $G$ is a cycle (respectively, path) that includes all vertices of $G$.

A graph $G$ is \emph{connected} if there exists a path between any two distinct vertices of $G$. If $u$ and $v$ are distinct vertices in a connected graph $G$, the \emph{distance} between $u$ and $v$, denoted by $\dist_G(u,v)$, is defined to be the minimum length 
 over all paths between $u$ and $v$. The \emph{diameter} of a connected graph $G$ is  the maximum value of $\dist_G(u,v)$ over all pairs of distinct vertices $u$ and $v$ of $G$   and is denoted by $\diam(G)$.  The {\emph{connected components} of a graph $G$ are {the  maximal} connected induced subgraphs of $G$. Two graphs $G$ and $H$ are \emph{disjoint} if $V(G)\cap V(H)=\emptyset$ and the \emph{disjoint union} of $G$ and $H$ is denoted by $G\sqcup H$.
 The  \textit{symmetric difference} of sets $A$ and $B$ is denoted by  $A
\ominus B=(A\cup B)\setminus (A\cap B)$.

Suppose $G_1$ and $G_2$ are disjoint graphs. The \emph{Cartesian product} of $G_1$ and $G_2$, denoted by $G_1\square G_2$, is the graph with $V(G\square  G_2)=V(G_1)\times V(G_2)$ such that $(v_1,v_2)$ and $(u_1,u_2)$ are adjacent if and only if  
  $v_1=u_1$ and $v_2u_2\in E(G_2)$, or
  $v_2=u_2$ and $v_1u_1\in E(G_1)$.
The \emph{join} of $G_1$ and $G_2$ is the graph $G_1 \vee G_2$ with vertex set $V(G_1 \vee G_2)=V(G_1) \cup V(G_2)$ and edge set $E(G_1) \cup E(G_2) \cup \{ v_1v_2:v_1 \in V(G_1), v_2 \in V(G_2)\}$.

Given an integer $n \geq 1$, the path, cycle and complete graph on $n$ vertices are denoted by $P_n$, $C_n$ and $K_n$, respectively. In examples, the vertices of $P_n$ and $C_n$ will be labeled  with the integers $1,\dots,n$, and  \emph{path (cycle) order} refers to the convention that the sequence of vertices in the path (cycle) is $(1,2,\dots,n)$.  The empty graph $\ol{K_n}$ is the graph with $n$ vertices and no edges.  
We denote the the complete bipartite graph by  $K_{p,q}$ and assume $1\le p\le q$ and the partite sets are $A=\{a_1,\dots,a_p\}$ and $B=\{b_1,\dots,b_q\}$.



\subsection{Parameters}\label{ss:param}
 In this section we define most of the parameters we discuss.  
Standard zero forcing was introduced in multiple applications, including control of quantum systems and as an upper bound to the maximum nullity of a real symmetric matrix $A$ whose off-diagonal entries $a_{ij}$ are nonzero 
 if and only if the graph has edge $ij$.  Other types of zero forcing were defined to serve as similar bounds for maximum nullity among positive semidefinite or skew-symmetric matrices (the distinguishing factor for skew  zero forcing is that the diagonal entries must all be zero).   

Starting with an initial set of  blue vertices $S$, a zero forcing process colors vertices blue by repeated applications of  a \emph{color change rule}; the color change rule determines the type of zero forcing.  Here we discuss three types of zero forcing: standard, introduced in \cite{AIM08}; positive semidefinite (PSD), introduced in \cite{smallparam}; and skew, introduced in \cite{IMA10}.   In each case, $B$ denotes the set of  (currently) blue vertices and the set of (currently) white vertices  $W = V(G)\setminus B$. 
\bit 
\item {\bf Standard color change rule:} Any blue vertex $u\in B$ can change the color of a white vertex $w$ to blue if $w$ is the only white neighbor of $u$, i.e., $N(u)\cap W=\{w\}$. 
\item {\bf PSD color change rule:}  Let $W_1,\dots, W_k$ be the sets of vertices of the $k\ge 1$ components of $G[W]$.  If $u\in B$, $w\in W_i$, and $w$ is the only white neighbor of  $u$ in $G[W_i\cup B]$, then $u$ can change the color  $w$ to blue. 
\item {\bf Skew color change rule:} Any vertex $u\in V(G)$ can change the color of a white vertex $w$ to blue if $w$ is the only white neighbor of $u$, i.e., $N(u)\cap W=\{w\}$. 
\eit
Given an initial set $S$ of blue vertices of    $G$, the  \emph{standard, PSD, or skew final coloring}    is the subset of vertices that  are blue after applying the relevant color change rule until no more changes are possible (the final coloring is unique for each of these color change rules).  A \emph{standard, PSD, or skew zero forcing set}   for  $G$  is a subset $S$ of vertices $S$ such that the final coloring is $V(G)$ when starting with exactly the vertices in $S$ blue. The \emph{standard, PSD, or skew  zero forcing number} of a graph $G$ is  the minimum cardinality of a standard, PSD, or skew zero forcing set; this parameter is denoted by $\Z(G)$, $\Zp(G)$, or $\Zs(G)$, respectively. All these parameters (and power domination) are also discussed in \cite[Chapter 9]{HLS22book}.  TAR  reconfiguration graphs of standard zero forcing were studied in \cite{XTARiso}, where it was noted that the main results also apply to positive semidefinite TAR reconfiguration graphs; these parameters are discussed here in Sections \ref{ss:ZF} and \ref{s:PSD}, respectively. TAR  reconfiguration graphs of skew zero forcing has not been studied previously and is discussed here in Section \ref{s:skew}.

  A set $S\subseteq V(G)$ is a \emph{power dominating set} of $G$ if $N[S]$ is a standard zero forcing set of $G$.  Power domination on graphs was introduced to model the placement of Phase Monitoring Unities (PMUs), which are used to monitor electric networks to avoid catastrophic failures (see, for example, \cite{HHHH02} and \cite{BH}). 
TAR  reconfiguration graphs of power domination were studied in \cite{PDrecon} and are discussed Section \ref{ss:PD}.

The study of irredundance was introduced in 1978 by Cockayne, Hedetniemi and Miller \cite{CHM78} as part of a study of minimal dominating sets. 
 Let $G$ be a graph and $T\subseteq V(G)$. 
A \emph{private neighbor} of $x\in T$ (relative to $T$)  is a vertex $w$ such that $x$ is the unique vertex in $T$ that dominates $w$.  The set $T$ is \emph{irredundant} or an \emph{Ir-set} if every $x\in T$ has a private neighbor relative to $T$.  
The  \emph{upper Ir number} is $\IR(G) = \max\{|T|: T \text{ is a maximal Ir-set}\}$. 

Forts play a fundamental role in blocking standard zero forcing and were introduced in \cite{FH18}.  Let $G$ be a graph. A nonempty set $F\subseteq V(G)$ is a \emph{fort} if $|F\cap N(v)|\ne 1$ for all $v\in V(G)\setminus F$.  Forts were used in \cite{ZIR} to define standard zero forcing irredundance. For $T\subseteq V(G)$ and $x\in T$, a fort $F$ of $G$ is a \emph{private fort} of $x$ relative to $T$ provided that $T\cap F = \{x\}$. As defined in \cite{ZIR}, the set $T\subseteq V(G)$ is a \emph{$\Z$-irredundant set} or \emph{ZIr-set} provided every element of $T$ has a private  fort, the  \emph{upper ZIr number} is $\ZIR(G) = \max\{|T|: T \text{ is a maximal ZIr-set}\}$. 

 A set $T$ of vertices in a graph $G$ is \emph{independent} (or is an \emph{independent set}) of $G$ if no two distinct vertices in $T$ are adjacent.  The \emph{independence number} of $G$, denoted by $\alpha(G)$, is the maximum cardinality of an independent set of vertices of $G$. 
A set $S$ of vertices in a graph $G$ is a \emph{vertex cover} of $G$ if every edge of $G$ has at least one of its endpoints in $S$.  The \emph{vertex cover number} of $G$, denoted by $\vc(G)$, is the minimum cardinality of a vertex cover of vertices of $G$. The independence number and vertex cover number are widely studied parameters. It is
well-known (and easy to see) that a set $T\subseteq V(G)$ is independent if and only if $V(G)\setminus T$ is a vertex cover of $G$.  TAR  reconfiguration graphs of vertex cover number and independence number are discussed Sections \ref{s:vtxcover} and \ref{ss:independent}.

 \section{$X$-set parameters (supersets and minimal sets): Axioms required for specific results}\label{sec:proof-requirements}

Numerous results for $X$-set parameters are established in \cite{PDrecon} and \cite{XTARiso}.  
 Many of these results still hold when some of the axioms in the definition of an $X$-set are removed.
In this section we  examine which  of these axioms  are necessary for each  of the aforementioned  results.  This is accomplished by introducing several universal variants  of  the original definition of $X$-sets, including the definition of robust $X$-set parameters  for which the main results of \cite{PDrecon, XTARiso} still hold. A 
benefit of this definition is that these results now apply to skew forcing, which had been missing from the family of parameters related to zero forcing.  By allowing consideration of properties that do not require every isolated vertex to be in an $X$-set, these results apply to vertex covers (without any assumption  about isolated vertices).
Just as \cite{PDrecon} extended some results on domination to original $X$-set parameters, later in this section we also extend additional results from domination to super  $X$-set parameters (sometimes requiring additional axioms),  including results on cut-vertices and Hamilton paths.

The next definition states the axioms that were assumed (including  implicit assumptions) for $X$-set parameters as defined in Definitions 2.1 in \cite{PDrecon} and 1.2 in \cite{XTARiso}. 

\begin{defn}\label{def:original+assumed}
 An {\em original $X$-set parameter} is  
   a graph parameter $X(G)$ defined to be the minimum cardinality of an $X$-set of $G$, where the $X$-sets of $G$  are defined by a given property and satisfy the following  axioms. 

\begin{enumerate}[(i)]
\item(Superset) 
If $S$ is an $X$-set and $S\subseteq S'$, then $S'$ is an $X$-set.

\item\label{axo:n-1}  ($(n-1)$-set) If $G$ has no isolated vertices, then every set of $|V(G)|-1$ vertices is an $X$-set. 
\item\label{ax:component}
(Component consistency) Suppose $G = G_1\du \dots \du G_k$ where $G_i,i=1,\dots,k$ are the connected components of $G$.  Then $S$ is an $X$-set of $G$ if and only if $S\cap V(G_i)$ is an $X$-set of $G_i$ for $i=1,\dots,k$.

\item\label{ax:whole}  (Isolated vertex) $\Xo(K_1)=1$. 

\item\label{ax:empty} (Empty set) The empty set is not an $X$-set.
\end{enumerate}
\end{defn}

Note that the Isolated vertex axiom \eqref{ax:whole} or similar was clearly intended in \cite{PDrecon, XTARiso} but not stated; without this or another assumption about $K_1$ having an $X$-set, a graph with no edges might have no $X$-sets. 
 The Connected component axiom \eqref{ax:component} was stated only for one direction: {\em If $S$ is an $X$-set of $G$, then $S\cap V(G_i)$ is an $X$-set of $G_i$ for $i=1,\dots,k$}.

 While the focus of this section is on parameters defined to be the minimum cardinality over sets satisfying a given property, we begin with the more general Definition \ref{def:vert set prop}. This definition will facilitate an extension of the $X$-set idea to parameters defined to be the maximum cardinality over sets satisfying a given property; these are called $Y$-sets and are discussed in Section \ref{s:maxYparam}.

 \begin{defn}\label{def:vert set prop}
A \emph{vertex set property} $W$ is a property that is defined on subsets of the vertex set of each graph such that if $\varphi: V(G) \to V(G')$ is a graph isomorphism, then for every $S\subseteq V(G)$ with property $W$, $\varphi(S)$ has property $W$. Let $W$ be a vertex set property and let $G$ be a graph. Then $S\subseteq V(G)$ is called an \emph{$W$-set} if it has property $W$.  

A vertex set property $W$ is called  \emph{cohesive}, or is a \emph{cohesive 
property}, if every graph has at least one $W$-set.  
A  graph parameter for which the value associated to each graph is solely determined by a cohesive property is called a \emph{cohesive parameter}.
\end{defn}

There are many natural graph parameters that can be associated to a cohesive  property. For example, the minimum or maximum cardinality of an $W$-set, or the total number of $W$-sets.  For a cohesive parameter, we can define the TAR graph  of a base graph.

\begin{defn}
  For a cohesive $W$-set property $W$,   the {\em token addition and removal  reconfiguration graph (TAR graph)}  of a  base graph  $G$ is the graph defined as follows:  The vertex set of the TAR graph is  the set of all  $W$-sets of $G$.  There is an edge between  two vertices  $S_1$ and $S_2$ of the TAR graph of $G$ if and only if  $|S_1 \ominus S_2|=1$.
\end{defn}

The Superset axiom is the most fundamental of the original $X$-set axioms (for parameters that take the minimum cardinality of an $X$-set) and is needed for almost every result established in \cite{PDrecon, XTARiso}. For example,  the Superset axiom ensures that the TAR reconfiguration graph is connected,  making the study of connectedness of certain subgraphs of the TAR graph meaningful.   
So we give a name to a cohesive  parameter that is the minimum size of an $X$-set and  satisfies the Superset axiom (and no other assumptions).   
Throughout what follows, we will always use an adjective  to modify \emph{$X$-set parameter}, such as \emph{original}, \emph{super}, or the name of one of the new variants (all of which are super $X$-set parameters).  However, we will use the same notation for the TAR graph and additional parameters derived from any type of $X$-set parameter.  Since  the only notation used for reconfiguration   graphs is for  TAR graphs, we have omitted the superscript \emph{TAR} used in \cite{PDrecon, XTARiso}.

\begin{defn}\label{d:superX}
 A \emph{super $X$-set 
parameter}  is a cohesive parameter $X$ such that $X(G)$ is defined to be the minimum cardinality of an $X$-set of $G$  where   the $X$-sets of $G$ satisfy the following condition:  
 \ben[(I)]
\item(Superset) 
If $S$ is an $X$-set and $S\subseteq S'$, then $S'$ is an $X$-set.
\een

 When $X$ is a super $X$-set parameter, the {\em $X$-TAR graph}  of a  base graph  $G$ is denoted by  $\xtar(G)$, and the \emph{upper $X$ number}, denoted by  $\xbar(G)$, is the maximum cardinality of a minimal $X$-set.
  \end{defn}


The next set of axioms allows recovery of all main $X$-set parameter results in \cite{PDrecon, XTARiso}, and all these axioms are needed for  many of the main results in full generality. Removal of the Empty set axiom does not  require changes to most of the proofs, and allows application to skew forcing.  Removal of the Isolated vertex axiom (which states that  $X(K_1) =1$) requires more changes but also allows the  main results to be  generalized to a larger class of parameters, including vertex cover number. Although the inclusion of the Superset axiom in the next definition is redundant, we think it is preferable to explicitly list all required axioms.

\begin{defn}\label{Robust-parameters-2}
A \emph{robust  $\Xr$-set parameter} is a super $X$-set parameter $X$  such that $X(G)$ and  the $X$-sets of $G$ satisfy the following conditions:  
\ben[(I)]
\item\label{axr:super}  (Superset)
If $S$ is an $\Xr$-set of  $G$ and  $S\subseteq S'$, then $S'$ is an $\Xr$-set of $G$.
\item\label{axr:n-1} ($(n-1)$-set) If $G$ is a connected graph of order  $n\ge 2$, then every set of $n-1$ vertices is an $\Xr$-set. 
\item\label{axr:component}
(Component consistency) Suppose $G = G_1\du \dots \du G_k$ where $G_i,i=1,\dots,k$ are the connected components of $G$.  Then $S$ is an $X$-set of $G$ if and only if $S\cap V(G_i)$ is an $X$-set of $G_i$ for $i=1,\dots,k$.\een  
\end{defn}

\begin{obs}\label{r:n-1-isol}
    Let $X$ be a robust $X$-set parameter. Note that if $G$ has no isolated vertices, then every set of $|V(G)|-1$ vertices is an $\Xr$-set.
\end{obs}

  From the previous observation, we see that the definition of robust $X$-set  is equivalent to that obtained by replacing axiom \eqref{axr:n-1} of Definition \ref{Robust-parameters-2} by \eqref{axo:n-1} from Definition \ref{def:original+assumed} (and keeping the other two axioms of Definition \ref{Robust-parameters-2} unchanged). The change in the $(n-1)$-set axiom was made because the behavior of $X$-sets  is different when $X(K_1)=0$ compared to that of  an original $X$-set parameter (which implies $X(K_1)=1$). 
 Both cases   $X(K_1)=1$ and $X(K_1)=0$ are discussed in the next remark.


\begin{rem}\label{no-isol-vtx-2}    Let $X$ be a robust $X$-set parameter.
Then one of two things happens regarding isolated vertices for every graph $G$:
\ben[$(1)$]
\item\label{XK1=1}    \cite[Remark 2.3]{PDrecon} $X(K_1)=1$: Every isolated vertex of $G$ is in every $X$-set of $G$.  In this case, if $G'=G\sqcup rK_1$ (where $G$ has no isolated vertices), then $X(G')=X(G)+r$ and 
$\xtar(G')\cong\xtar(G)$.  
\item\label{XK1=0} $X(K_1)=0$: No isolated vertex of $G$ is in any minimal $X$-set of $G$.  Every set of $|V(G)|-1$ vertices is an $X$-set (even if $G$ has isolated vertices). Suppose $G'=G\sqcup rK_1$ (where $G$ has no isolated vertices).  Then $X(G')=X(G)$.  However, $\xtar(K_1)\cong K_2$ when $X(K_1)=0$.  We see in  the next proposition 
 that 
$\xtar(G')\cong\xtar(G)\Box K_2 \Box \cdots\Box K_2$ ($r$ copies of $K_2$). 
\een
 Therefore, it is sufficient to study  TAR reconfiguration graphs of base graphs with no isolated vertices  for robust $X$-set parameters. 
\end{rem}

Case \eqref{XK1=1} of the previous remark applies  to all parameters studied in \cite{PDrecon} and \cite{XTARiso}.  The vertex cover number  is an example of  a robust $X$-set parameter where case \eqref{XK1=0} applies.

The next result was established for original $X$-set parameters in Proposition 1.4 of \cite{XTARiso}. It reduces consideration of a TAR reconfiguration problem to connected base graphs. 

\begin{prop}\label{p:disjoint- cart-u}
Let $X$ be a super $X$-set parameter that satisfies the Component consistency  axiom and let $G = G_1\du G_2$. Then $\xtar(G)  \,\cong\, \xtar(G_1)\,\Box\,\xtar(G_2)$.
\end{prop} 
{\blu \begin{proof}
Let $S$ and $T$ be $X$-sets of $G$. Then $S = S_1\sqcup S_2$ and $T = T_1 \sqcup T_2$, where $S_1$ and $T_1$ are $X$-sets of $G_1$, and $S_2$ and $T_2$ are $X$-sets of $G_2$. Observe that $S$ and $T$ are adjacent in $\xtar(G)$ if and only if there exists a vertex $v\in V(G_1)\sqcup V(G_2)$ such that
\begin{itemize}
    \item $S_1 = T_1$ and ($ T_2=S_2\setminus \{v\}$ or $S_2 = T_2\setminus \{v\}$), or
    \item $S_2 = T_2$ and ($ T_1=S_1\setminus \{v\} $ or $S_1 = T_1\setminus \{v\}$).
\end{itemize}
Further, $ T_2=S_2\setminus \{v\}$ or $S_2 = T_2\setminus \{v\}$ if and only if $S_2$ is adjacent to $T_2$ in $\xtar(G_2)$, and $ T_1=S_1\setminus \{v\} $ or $S_1 = T_1\setminus \{v\}$ if and only $S_1$ is adjacent to $T_1$ in $\xtar(G_1)$.
\end{proof}}

 As noted earlier, every original $X$-set parameter is a robust $X$-set parameter. This includes (standard) zero forcing, PSD forcing, domination, and power domination, as noted in \cite{XTARiso}. Skew forcing is a robust $X$-set parameter that is not an original $X$-set parameter. The vertex cover number is another (very different) robust $X$-set parameter.

Of course, one can also define examples of robust $X$-set parameters that have no known use, as in the next example.

\begin{ex}
    Given a connected graph $G$ of order $n\ge 2$, a set $S\subseteq V(G)$ is a  \emph{star-set} of $G$ 
    if and only if $|S|\ge |V(G)|-1$.  The only star-set of $K_1$ is its vertex. Given a disconnected graph $G$, a set $S\subseteq V(G)$ is a star-set of $G$ if and only if $S\cap V(C)$ is a star-set of $C$ for every connected component $C$ of $G$. Define $X^*(G)$ to be the minimum cardinality of a star-set of $G$.  It is clear that for any graph $G$ of order $n$, $X^*(G)=n-c$ where $c$ is the number of connected components of $G$ of order at least two.  Furthermore, for every connected graph of order $n$, the star reconfiguration graph is $\X^*(G)\cong K_{1,n}$.
\end{ex}

 It is worth noting that 
the Component consistency axiom is irrelevant for connected graphs, which leads to the next definition and observation.

 \begin{defn}\label{Connected-parameters}
  A \emph{connected $X$-set parameter} is a super $X$-set parameter $X$  such that $X(G)$ and  the $X$-sets of $G$ satisfy the following conditions:  
\ben[(I)]
\item (Superset)
If $S$ is an $\Xc$-set of  $G$ and  $S\subseteq S'$, then $S'$ is an $\Xr$-set of $G$.

\item  (($n-1$)-sets) If $G$ is a connected graph of order  $n\ge 2$, then every set of $ n-1$ vertices is an $\Xr$-set. 
\een  
\end{defn}

\begin{obs}\label{r:component}  
If $\Xc$ is a connected graph parameter, 
then every result that is true for a robust $X$-set parameter and graphs with no isolated vertices holds for $\Xc$ for connected graphs  of order at least two. 
\end{obs}

While most of the main results of this section are stated for robust parameters and graphs with no isolated vertices, they remain true for connected graphs and connected parameters as noted in Observation \ref{r:component}.

\subsection{Initial Results}

Most of the results in this section  were  established in \cite{PDrecon} for original $X$-set parameters.   
The first statement in the next remark is  Observation 2.4 of \cite{PDrecon}.

\begin{rem}\label{r:sameXsets}     Let $\Xb$ be a super $X$-set parameter.  Knowledge of all the minimal $X$-sets is sufficient to determine $\xtar(G)$.

Suppose $G$ is a graph and $X_1$ and $X_2$ are super $X$-set parameters   such that $M\subseteq V(G)$ is a minimal $X_1$-set of $G$ if and only if minimal $X_2$-set of $G$.  Since the TAR graph is determined by the minimal $X$ sets, $\xtar_1(G)=\xtar_2(G)$.
\end{rem}

 Since several of the parameters we study in the following sections have the same minimal $X$-sets for certain graph families, this allows the transfer of results about TAR graphs from one parameter to another for such families of graphs.


  \begin{rem}\label{r:minXsets}
         Suppose $\Xb$ is  a super $X$-set parameter that also satisfies the $(n-1)$-set axiom and let $G$ be a graph of order at least two.  If $X(K_1)=0$ or if $G$ has no isolated vertices, then no vertex is in every minimal $\Xb$-set. If in addition the empty set is not an 
    $\Xb$-set of $G$, then $G$ has more than one minimal $\Xb$-set. 
 \end{rem}

The inequality $\Delta(\xtar (G))\le n$ is true  for all super $X$-set parameters and graphs, without additional assumptions, as seen in the proof of Proposition 2.5(1) in \cite{PDrecon}.  The next result also
shows how closely related the $(n-1)$-set axiom is to $\Delta(G)=n$.

 \begin{prop}\label{need-n-1}
Let $X$ be a super $X$-set parameter and $G$ be a graph of order $n$. Then $\Delta(\xtar(G))=n$ if and only if every set of $n - 1$ vertices is an $X$-set. 
\end{prop}
\begin{proof}
 {\blu Let $S\subseteq V(G)$ be a vertex of $\xtar(G)$ and let $|S|=r$. Then there are at most $r$ vertices that can be deleted (one at a time) leaving an $X$-set, and at most $n-r$ vertices that can be added (one at a time).  Thus $\deg_{\xtar(G)}(S)\le n$.} It is  immediate that $ \Delta(\xtar(G))=n$  if every set of $n - 1$ vertices is an $X$-set. 

 Suppose that $ \Delta(\xtar(G))=n$. Then there exists an $X$-set $S$ such that $S\cup \{v\}$ and $S\setminus\{w\}$ are $X$-sets for every $v\notin S$ and $w\in S$. For $u\in V(G)$, define $R_u = V(G)\setminus \{u\}$. Then $R_u$ is a superset of $S$ if $u\notin S$ and is a superset of $S\setminus \{u\}$ if $u\in S$. In either case, $R_u$ is an $X$-set by the Superset axiom.
\end{proof}

 The  next result    is central to many of the results that follow, including the isomorphism results. If $X(K_1)=0$ or   $G$ has no isolated vertices,  then   every set of $n-1$ vertices is an $\Xr$-set by   Remark \ref{no-isol-vtx-2} or Observation \ref{r:n-1-isol}, so $\deg_{\xtar(G)}(V(G))= n$.

\begin{cor}\label{x-deg} 
 Let $X$ be a robust $X$-set parameter.  If $X(K_1)=0$ or   $G$ is a graph  of order $n$ with no isolated vertices, 
then $\Delta(\xtar (G))=n.$
\end{cor}

 Example \ref{ex:X=V} illustrates the necessity of the $(n-1)$-set axiom for the previous result   and for results that depend on it.

\begin{ex}\label{ex:X=V}  
    Define $X_V$ to be the property that for every graph $G$, $V(G)$ is the one and only $X_V$-set,  and define $X_V(G)$ to be the minimum cardinality of an  $X_V$-set. Observe that  $X_V$ is a super $X$-set parameter that satisfies all of the axioms for an original $X$-set parameter except the $(n-1)$-set axiom. Since   $X_V(G)=|V(G)|$ and  $\xtar_V(G)=K_1$ for every graph $G$, for a graph $G$  with no isolated vertices,  neither the order of  $G$ nor $X_V(G)$ can be determined from $\xtar_V(G)$.  Many subsequent results also fail for  $X_V$.  
\end{ex}

 The connected domination number number, which is discussed in  Section \ref{s:conD}, provides an example of a parameter in the literature that is a super $X$-set parameter but not a robust $X$-set parameter, and for which Corollary \ref{x-deg} and TAR graph isomorphism results fail.

 The next result was established for original $X$-set parameters in Proposition 2.6 of \cite{PDrecon} (and used in the proof of what is Proposition \ref{minimal-u} here).   Since only the last statement was explicitly stated and proved, we included the brief revised proof here.

 \begin{prop}\label{p:umindegreeTAR}
Let $X$ be a super $X$-set parameter, let $G$ be a graph on $n$ vertices, and let  $S$ be an $X$ set. 
Then $\deg_{\xtar(G)}(S)\ge n-|S|$ and $S$ is minimal if and only if $\deg_{\xtar(G)}(S)=n-|S|$.  Furthermore, $\delta(\xtar(G))=n-\xbar(G)$.
\end{prop}

\bpf
Let
$\X=\xtar(G)$.   
Let $M$ be a minimal $X$-set of $G$, with $|M|=m$ for some $m\geq 0$.  Since $N_{\X}(M)=\{M \cup \{v\}: v \in V(G)\setminus M\}$,  $n-m=\deg_{\X}(M)$. 

 Suppose $S$ is not minimal and choose a minimal $X$-set  $M \subsetneq S$. Denote the elements of $S\setminus M$ by $U= \{u_1,\ldots,u_\ell\} $ and those of  $V(G)\setminus S$ by $W=\{w_1,\ldots, w_{n-m-\ell} \}$. Then the sets of the form $S\setminus \{u_i\}$ with $u_i \in U$ and $S \cup  \{w_j\}$ with $w_j \in W$ are neighbors of $S$ in ${\X}$, and $|N_{\X}(S)| \geq \ell > n-m$. 
The last statement is  now immediate.  
\epf

Additional minor results adapted  from \cite{PDrecon} can be found in Section \ref{ss:Xconn}.
\subsection{Hypercube representation and applications}

 For $
  d\ge 1$, the graph having   as vertices  all  $d$-tuples with entries in $\{0,1\}$ with two such sequences adjacent  if and only if  they differ in exactly one place is a characterization of $Q_d$, the hypercube of dimension $d$  (and $Q_0\cong K_1$). There is a well-known representation of any TAR reconfiguration graph as a subgraph of a hypercube.  Let $G$ be a graph with $ V(G)=\{v_1,\dots,v_n\}$. Any subset $S$ of $V(G)$ can be represented by a sequence $(s_1, \dots, s_n)$ where $s_i=1$ if $v_i\in S$ and $s_i=0$ if $v_i\not\in S$.     The first part of the next remark appeared in Remark 2.12 of \cite{PDrecon}

\begin{rem}\label{r:hypercube-bip} Let $X$ be a super $X$-set parameter. For a graph $G$ of order $n$, $\xtar(G)$ is isomorphic to a subgraph of $Q_n$, and thus $\xtar(G)$ is bipartite. If $G$ is a graph such that  any one vertex is an $X$-set and the emptyset is not an $X$-set, then    $\xtar(G)\cong Q_n-v$ for any $v\in V(Q_n)$.  
\end{rem}

\begin{rem}\label{r:noQn}    It is asserted in Remark 2.13 of \cite{PDrecon} that the $X$-TAR graph is not a hypercube for an original  $X$-set parameter and base graph of order at least two. This is not technically true: If $X$ is  an original $X$-set parameter, then $\xtar(nK_1)\cong Q_0$.  However, the intended statement,  the $X$-TAR graph is not a hypercube for an original  $X$-set parameter and base graph that has an edge, is true and the argument in \cite[Remark 2.13]{PDrecon} justifies that statement.

 However, a robust $X$-set parameter such as skew forcing, may allow the $X$-TAR graph of a base graph with an edge to be a hypercube: Any graph $G$ such that  $\emptyset$ is an $X$-set has $\xtar(G)\cong Q_n$ for a robust $X$-set parameter.   The $(n-1)$-set  axiom also plays a role, since the super $X$-set parameter $X_V$ in Example \ref{ex:X=V} satisfies the  (original) Isolated vertex, Connected component, and Empty set axioms but not the $(n-1)$-set  axiom, and $\xtar_V(G)\cong Q_0$ for all graphs $G$.
 \end{rem}

  The next lemma was established for original $X$-set parameters in Lemma 2.14 of \cite{PDrecon}; its corollary is immediate and appears for original $X$-set parameters in Corollary 2.15 of \cite{PDrecon}. 
 
\begin{lem}\label{L:uembedQt} Let $X$ be a super $X$-set parameter and let $G$ be a graph on $n$ vertices and let $t \leq n$. Then   $X(G)\leq n-t$ if and only if $\xtar(G)$ has an induced subgraph isomorphic to the hypercube $Q_t$. 
\end{lem}
{\blu
\bpf 
Let $S$ be a minimum $X$-set and let $W=V(G)\setminus S$.  The induced subgraph of $\xtar(G)$ having vertices consisting of sets of the form $S\sqcup W'$ over all subsets $W'\subseteq W$  is $Q_{n-X(G)}$.  Any hypercube $Q_{t}$ for $n-X(G)\ge t\ge 0$ is an induced subgraph of $Q_{n-X(G)}$.

 Suppose $H$ is an induced subgraph of $\xtar(G)$ isomorphic to $Q_t$ for some $1\leq t \leq n$. Choose $S \in V(H)$ such that $|S|$ is minimum over all $X$-sets in $H$.  Since no vertex in $H$ has fewer vertices than $S$, every one of the $t$ neighbors of $S$ in $H$ is obtained by adding a vertex of $G$ to $S$.  Thus $|V(G)\setminus S| \ge t$ and $X(G) \le |S|\le n-t$.
\epf 
}

\begin{cor}\label{c:X-Qt}Let $X$ be a super $X$-set parameter  and let $G$ be a graph on $n$ vertices.  
Then $d=n-X(G)$ is the  maximum dimension of a hypercube isomorphic to an induced subgraph of the reconfiguration graph 
 $\xtar(G)$. 
\end{cor}

 The next result follows from prior results and was established for original $X$-set parameters in Corollary 2.16 of \cite{PDrecon}. 
\begin{cor}\label{c:u-iso-props} Let $X$ be a robust $X$-set parameter  and
let $G$ and $G'$ be graphs such that $\xtar(G)\cong\xtar(G')$.  Suppose further that $X(K_1)=0$ or $G$ and $G'$ have no isolated vertices. Then 
\ben
\item $|V(G)|=|V(G')|$.
\item $X(G)=X(G')$.
\item $\xbar(G)=\xbar(G')$. 
\een
\end{cor}

{\blu \bpf  By Corollary \ref{x-deg}, $|V(G)|=\Delta(\xtar(G))=\Delta(\xtar(G'))=|V(G')|$. 
Corollary \ref{c:X-Qt} implies $X(G)=X(G')$. By Proposition \ref{p:umindegreeTAR}, $\xbar(G)=n-\delta(\xtar(G))=n-\delta(\xtar(G'))=\xbar(G')$. 
\epf}


\subsection{Isomorphism results}

 The original $X$-set parameter version of the  next result is Theorem 1.1 of \cite{XTARiso}, and we have verified  that the expanded version including $X(K_1)=0$ remains true for robust $\Xr$-set parameters.  The proof is established through a series of preliminary results, which are stated for robust $\Xr$-set parameters (or super $X$-set satisfying some additional axioms) after Theorem \ref{t:main}.  

\begin{thm}\label{t:main} Let $\Xr$ be a robust $X$-set parameter  and let $G$ and $G'$ be base graphs  such that $\xtar(G)\cong\xtar(G')$.  If $X(K_1)=0$ or $G$ and $G'$ have no isolated vertices, then there is a relabeling of the vertices of $G'$ such that $G$ and $G'$ have exactly the same $\Xr$-sets.  
\end{thm}

Observe that Proposition \ref{need-n-1}  shows that the $(n-1)$-set axiom (or at least something beyond the Superset  axiom)
is necessary for the previous result, as well as Theorem \ref{isomorph-isometry-u}  (which together with Theorem \ref{isomorph-same-zfs}  and Remark \ref{r:isomorph-same-zfs} implies Theorem \ref{t:main}).

 Next we reproduce some necessary definitions from \cite{XTARiso}. Let  $P=(Y,\leq)$ be a poset.   {\blu An element $u\in Y$ is \textit{maximal} (respectively, \textit{minimal}) if $u\not< y$ (respectively,  $u\not> y$) for each $y\in Y$.} Two elements $u$ and $v$ of $Y$ are \emph{comparable} provided $u\leq v$ or $v\leq u$. A poset such that each pair of elements is comparable is called a  \emph{chain}. An \emph{interval} in $P$ is a set $[u,v] = \{y : u\leq y \leq v\}$ where $u\leq v$.   The \emph{length of a chain} $C$ is $|C|-1$. The \emph{length of an interval} $[u,v]$ is the maximum length of a chain in $[u,v]$ and is denoted by $\ell(u,v)$. For a set $T$, the power set of $T$ is denoted by $\mathcal{P}(T)$.  If $X$ is a super $X$-set  property, then $\xtar(G)$ is a join-semilattice contained in $(\mathcal{P}(V(G)),\subseteq)$  and  $\ell(S,S')=\dist_{\xtar(G)}(S,S')$ for any interval $[S,S']$ in $\xtar(G)$ since $S$ and $S'$ are comparable.  
The next results was established for original $X$-set  parameters in Lemma 2.1 of \cite{XTARiso}.
 \begin{lem}\label{interval} Let $X$ be a super $X$-set  parameter, let $t\geq 0$ be an integer, let $G$ be a graph on $n$ vertices,  and let $H$ be an induced subgraph of $\xtar(G)$.  If $H\cong Q_t$, then $V(H)$ is an interval  of length $t$ in the poset $(\mathcal{P}(V(G)),\subseteq)$. \end{lem}  
 
 {\blu 
 
 \bpf Assume $H\cong Q_t$. The claim is obvious for $t=0,1$, so suppose that $t\geq 2$. We begin by showing that $(V(H),\subseteq)$ has exactly 1 maximal element. Assume, to the contrary, that $(V(H),\subseteq)$ has at least 2 maximal elements $u$ and $v$. Since $H$ is connected, there exists a path from $u$ to $v$ in $H$. Every path $P$ from $u$ to $v$ in $H$ can be written in the form $(u,y_1,\dots, y_i, v)$. For each such $P$, let $d_P$ to be the smallest index $k$ such that $y_{k-1}\subseteq y_k$, where $y_0 = u$ and $y_{i+1} = v$. Note that since $u$ and $v$ are maximal in $(V(H),\subseteq)$, $d_P\geq 2$ is always defined.
 
 Let $d$ be the minimum $d_P$ amongst all paths $P$ from $u$ to $v$ in $H$. Pick a path $P^*$ in $H$ of the form $(u,x_1,\dots, x_i, v)$ such that $d_{P^*} = d$, and let $x_0 = u$ and $x_{i+1} = v$. Then $x_{d-1}\subseteq x_{d-2}$  because $x_{d-2}\not\subseteq x_{d-1}$ and $x_{d-2}$ is adjacent to $x_{d-1}$. Since $x_{d-2}$ and $x_d$ have  $x_{d-1}$ as a common neighbor and each pair of vertices in a hypercube share exactly 0 or 2 common neighbors, there exists some vertex $w \not=x_{d-1}$ in $H$ that is adjacent to $x_{d-2}$ and $x_d$. By our choice of $d$, $w\subseteq x_{d-2}$ and $w\subseteq x_d$   or else $P'=(u,x_1,\dots, x_{d-2},w,x_{d},\dots, x_i,v)$ would be a path with $d_{P'}<d_{P^*}$ since  $w=x_{d-2}\cup x_d $. Hence $x_{d-2}\cap x_d = w$. This is absurd since $x_{d-2}\cap x_d = x_{d-1}$ and $w\not=x_{d-1}$.
 
 Thus, $H$ has exactly 1 maximal element $T$. A similar argument shows $H$ has exactly 1 minimal element $S$ in the subset partial ordering of $\xtar(G)$. Thus, $S\subseteq R\subseteq T$ for every $R\in V(H)$.  Since $H\cong Q_t$, $\dist(S,T)\le \diam(Q_t)=t$.  So there are at most $2^t$ elements in the interval $[S,T]$.  But $|V(H)|=|V(Q_t)|=2^t$. Therefore, $\dist(S,T)=t$ and $V(H)=[S,T]$.\epf}

The proof of Lemma 2.2 in \cite{XTARiso} also establishes  the next result.
 
  \begin{lem}\label{lem:ord_pres} Let $X$ be a super $X$-set  parameter and let  $G$ and $G'$ be graphs.   Suppose $\varphi:\xtar(G)\to \xtar(G')$ is a graph  isomorphism. Let $S^\prime = \varphi(V(G))$. If $M'$ is a minimal $X$-set of $G^\prime$, then $M'\subseteq S^\prime$. \end{lem} 

{\blu \bpf Let $M'$ be a minimal $X$-set of $G'$. Define $M = \varphi^{-1}(M')$. The interval $[M,V(G)]$ in $(\mathcal{P}(V(G),\subseteq)$ forms an induced $Q_t$ in $\xrtar(G)$ for some integer $t>0$. By Lemma \ref{interval}, $\varphi([M,V(G)])$ is an interval $[Z',W']$ in $(\mathcal{P}(V(G^\prime)),\subseteq)$. Since $M'$ is a minimal $X$-set of $G'$, $Z' = M'$. Thus, $M' \subseteq S'$. \epf}

Irrelevant vertices play a key role in the proofs of results about   TAR-graph isomorphsims of original $X$-set parameters in \cite{XTARiso}.  Here we extend the definition of irrelevant set to super $X$-set  parameters  and reproduce  results that remain true (some results require the parameter to be robust).  Let $G$ be a graph and let $X$ be a super $X$-set  parameter.   A vertex  $v\in  V(G)$ is \emph{$X$-irrelevant}  if $v\not\in S $ for every minimal $X$-set $S$ of $G$.  Observe that if $X(K_1)=0$  and $X$ satisfies the Component consistency axiom, 
 then every isolated vertex is irrelevant. A set $R\subseteq V(G)$ is an \emph{$X$-irrelevant set}  if every vertex of $R$ is $X$-irrelevant.
For  a graph $G$ and $R\subseteq V(G)$, define the map 
$\nu_R: V(\xtar(G)) \to \mathcal{P}(V(G))$ by $\nu_R(S) = S\ominus R$. The proof of Theorem 2.7  in \cite{XTARiso} for original $X$-set parameters also establishes the next result.

\begin{thm}\label{t:irrel-isom}
Let $X$ be a super $X$-set  parameter, let $G$ be a graph, and let  $R\subseteq V(G)$. Then $\nu_R$ is a graph automorphism of $\xtar(G)$ if and only if $R$ is $X$-irrelevant.
\end{thm}

\bpf\blu 
Suppose that $\nu_R$ is an automorphism of $\xtar(G)$. By Lemma \ref{lem:ord_pres}, every minimal $X$-set of $G$ is  a subset of $\nu_R(V(G)) = V(G)\setminus R$. Thus, $R$ is $X$-irrelevant.

Suppose that $R$ is $X$-irrelevant. Let $S$ be an $X$-set of $G$. Then there exists some minimal $X$-set $T\subseteq S$. Since $R$ is $X$-irrelevant, $T\subseteq S\setminus R \subseteq S$. Thus, $\nu_R(S)\supseteq S\setminus R \supseteq T$ is an $X$-set of $G$. Adjacency in $\xtar(G)$ is preserved by $\nu_R$: Consider the adjacent $X$-sets $S$ and $S'=S\cup\{u\}$.  If $u\ne R$, then  $S'\setminus R = (S\setminus R)\cup\{u\}$, so $\nu_R(S')$ is adjacent to $\nu_R(S)$.  If $u\in R$, then  $S\setminus R = (S'\setminus R)\cup\{u\}$, so $\nu_R(S)$ is adjacent to $\nu_R(S')$.  
Therefore, $\nu_R$ is an automorphism.
\epf

 Theorems \ref{isomorph-isometry-u} and \ref{isomorph-same-zfs}  (and Remark \ref{r:isomorph-same-zfs}) imply Theorem \ref{t:main}. The proof of Theorem 2.8 
 \cite{XTARiso} also establishes Theorem \ref{isomorph-isometry-u} (even with the addition of $X(K_1)=0$ case).
  We supply  a proof for Theorem  \ref{isomorph-same-zfs} that is modified from the proof of Theorem 2.9 in 
 \cite{XTARiso} to weaken the hypotheses so that $X$ need only be a super $X$-set parameter rather than a robust $X$-set parameter provided the two graphs have the same order.

\begin{thm}\label{isomorph-isometry-u}  Let $\Xr$ be a robust  $\Xr$-set parameter, let $G$ and $G'$ be graphs, and let  $\tilde\varphi:\xrtar(G)\to\xrtar(G^\prime)$ be an isomorphism.  Suppose further that  $X(K_1)=0$ or  $G$ and $G'$ have no isolated vertices. Then $R^\prime = V(G^\prime)\setminus \tilde\varphi(V(G))$ is $\Xr$-irrelevant and  $\varphi = \nu_{R^\prime}\circ \tilde\varphi$ is an isomorphism  such that $|\varphi(S)|=|S|$ for every $S\in V(\xrtar(G))$.
\end{thm}
\bpf \blu Let $n=|V(G)|$.
If $\tilde\varphi(V(G))=V(G')$, then $R'=\emptyset$ is $X$-irrelevant. 
So suppose that $\tilde\varphi(V(G)) = S^\prime$, where $S^\prime\not=V(G^\prime)$. 
Since $\tilde\varphi$ is an isomorphism, Lemma \ref{lem:ord_pres} implies every minimal $X$-set of $G^\prime$ is a subset of $S^\prime$ and hence $R^\prime$ is $X$-irrelevant. By Theorem \ref{t:irrel-isom}, $\nu_{R^\prime}$ is an automorphism of $\xtar(G^\prime)$. Thus, $\varphi = \nu_{R^\prime}\circ \tilde\varphi$ is an isomorphism such that $\varphi(V(G)) = V(G^\prime)$.

Note that $|V(G')|=n$ and $\Xr(G')=\Xr(G)$ by Corollary \ref{c:u-iso-props}. Let $S\in V(\xrtar(G))$.  The interval $H=[S,V(G)] \in \xrtar(G)$ is an induced hypercube  in $\xrtar(G)$, so $\varphi(H)$ is an induced hypercube in $\xrtar(G')$.  By Lemma \ref{interval} and since $\varphi(V(G))=V(G')$, $\varphi(H)= [S',V(G')]$ in $\xrtar(G')$ and $\dist(S',V(G'))=\dist(S,V(G))$. We show by induction on $|S|$ that $|\varphi(S)|=|S|$.
We say $S$ is a $k$-$\Xr$-set if $S$ is an $\Xr$-set and $|S|=k$.

For the base case, assume $|S|=\Xr(G)$, so $\dist(S',  V(G')) =\dist(S,  V(G)) =n-X(G)=|V(G')|-\Xr(G')$.  This implies $|S'|=\Xr(G')=\Xr(G)=|S|$. The same reasoning applies using $\varphi^{-1}$, since $\varphi^{-1}(V(G'))=V(G)$. 
Thus $\varphi$ defines a bijection between minimum  $\Xr$-sets of $G$  and  minimum  $\Xr$-sets  of $G'$.

Now assume $\varphi$ defines a bijection between  $i$-$\Xr$-sets of $G$ and  $i$-$\Xr$-sets of $G'$ for $\Xr(G)\le i\le k$ and let $S$ be a $(k+1)$-$\Xr$-set of $G$.  This implies  $|\varphi(W)|\ge k+1$ for $W\in [S,V(G)]$. By Lemma \ref{interval}, $\varphi([S,V(G)])= [S',V(G')]$ in $\xrtar(G')$ and $\dist(S',V(G'))=n-k-1$.  Thus $S'$ is a $(k+1)$-$\Xr$-set of $G$.
\epf

For any map $\psi:A \to A'$ and subset $B\subseteq A$ we write $\psi(B)$ to mean the image of $B$, i.e.,\ $\psi(B) = \{\psi(b) : b\in B\}$.  
This is particularly useful when working with a map $\psi:V(G)\to V(G')$ that maps $X$-sets of $G$ to $X$-sets of $G'$, since this convention naturally induces a map $\psi:V(\xtar(G))\to V(\xtar(G'))$.

 \begin{thm}\label{isomorph-same-zfs} 
Let $X$ be a super $X$-set parameter, let $G$ and $G'$ be graphs of order $n$, and suppose $\varphi:\xtar(G)\to \xtar(G')$ is a graph isomorphism. 

Then  $|\varphi(S)|=|S|$ for every $\Xr$-set $S$ of $G$ if and only if there exists a bijection $\psi:V(G)\to V(G')$ such that $\psi(S)=\varphi(S)$ for every $X$-set $S$ of $G$.

If $\varphi:\xtar(G)\to \xtar(G')$ is a graph isomorphism such that $|\varphi(S)|=|S|$ for every $\Xr$-set $S$, 
then there is a relabeling of the vertices of  $G'$ such that the relabeled graph has the same $\Xr$-sets as $G$ and the same $\Xr$-TAR graph as $G$.
\end{thm}
\bpf
If there exists a bijection $\psi:V(G)\to V(G')$ such that $\psi(S)=\varphi(S)$ for every $X$-set $S$ of $G$, then it is immediate that $|\varphi(S)|=|S|$ for every $\Xr$-set $S$ of $G$.  Assume $|\varphi(S)|=|S|$ for every $X$-set $S$ of $G$. For any graph $H$ and $v\in V(H)$, define $R_v=V(H)\setminus\{v\}$. If $R_v$ is an $X$-set of $G$, let $v'$ be the unique vertex such that $\varphi(R_v)=R_{v'}$. Let $W$ be the set of all vertices $w\in V(G)$ such that $R_{w}$ is not an $X$-set of $G$ and let $W'$ be the set of all vertices $w'\in V(G')$ such that $R_{w'}$ is not an $X$-set of $G'$. Since $|\varphi(S)|=|S|$ for every $X$-set $S$ of $G$, $|W| = |W'|$ and so there exists a bijection $\theta: W\to W'$. Let $\psi:V(G)\to V(G')$ be given by
\[
\psi(v)=
\begin{cases}
v' & \text{if $S_v$ is an $X$-set of $G$,}\\
\theta(v) & \text{otherwise.}
\end{cases}
\]
Observe that $\psi:V(G)\to V(G')$ is a bijection. 
  The proof that $\psi(S)=\varphi(S)$ for every $\Xr$-set $S$ of $G$ is as  the proof of Theorem 2.9 in \cite{XTARiso}, proceeding iteratively from $|S|=n$ to $|S|=X(G)$. 
  
  {\blu   By the choice of $\varphi$ and the definition of $\psi$, we have $\varphi(S)=\psi(S)$ for $|S|=n,n-1$.  Assume $\varphi(S)=\psi(S)$ for each $\Xr$-set $S$ of order $k$ for some $k$ with $n-1\ge k>X(G)$.  Let $S$ be an $\Xr$-set of order $k-1$.   
 Since $k-1\le n-2$, there exist distinct vertices $a,b\in V(G)\setminus S$ such that $S\cup \{a\}$ and $S\cup \{b\}$ are $\Xr$-sets of order $k$. Since $S$ is adjacent to $S\cup\{a\}$ and $S\cup\{b\}$ in $\xtar(G)$, and $|\varphi(S\cup\{a\})|=|\varphi(S\cup\{b\})| > |\varphi(S)|$, there exist distinct $a',b'\in V(G')\setminus \varphi(S)$ such that $\varphi(S\cup\{a\}) = \varphi(S)\cup\{a'\}$ and $\varphi(S\cup\{b\}) = \varphi(S)\cup\{b'\}$. Thus,
 \[
\varphi(S)=\varphi(S\cup \{a\})\cap \varphi(S\cup \{b\})=\psi(S\cup \{a\})\cap \psi(S\cup \{b\})=\psi(S),
 \]
 where the last equality follows since $\psi$ is a bijection.}

Similarly, the proof that 
 if $\varphi:\xtar(G)\to \xtar(G')$ is a graph isomorphism such that $|\varphi(S)|=|S|$ for every $\Xr$-set $S$, then there is a relabeling of the vertices of  $G'$ such that the relabeled graph has the same $\Xr$-sets as $G$ and the same $\Xr$-TAR graph as $G$ follows the proof of Theorem 2.9 in \cite{XTARiso}.
 {\blu Finally, suppose  $|\varphi(S)|=|S|$ for every $\Xr$-set $S$  and $\psi:V(G)\to V(G')$ is a bijection such that $\psi(S)=\varphi(S)$  for every $\Xr$-set $S$ of $G$. Define $G''$ from $G'$ by relabeling vertices of $G'$ so that $v'\in V(G')$ is labeled by $\psi^{-1}(v')$.  Then $G''\cong G'$ and $G''$ and $G$ have the same $\Xr$-sets. }
 \epf\vspace{-8pt}

  \begin{rem}\label{r:isomorph-same-zfs}
    Let $X$ be a robust $X$-set parameter, and let $G$ and $G'$ be graphs such that $\varphi:\xtar(G)\to\xtar(G')$ is a graph isomorphism.  If $X(K_1)=0$ or $G$ and $G'$ have no isolated vertices, then  Corollary \ref{c:u-iso-props} implies that the hypotheses of Theorem \ref{isomorph-same-zfs} are satisfied. 
\end{rem}


\subsection{$\Xr$-irrelevant vertices and automorphisms of $\Xr$-TAR graphs.}\label{ss:irrelevant}

In this section we point out that the proof of the characterization of automorphisms of $X$-TAR graphs in Theorem 2.13 in \cite{XTARiso}  remains true for robust parameters. 
 Let $M_{\Xr}(G)$ denote the set of bijections 
 $\psi:V(G)\to V(G)$ that send minimal $\Xr$-sets of $G$ to minimal $\Xr$-sets of $G$ of the same size.

\begin{thm}\label{t:auto-univ}
Let $\Xr$ be a robust $X$-set parameter and let $G$ be a graph.  If $X(K_1)=0$ or $G$ has no isolated vertices, then the automorphism group of $\xrtar(G)$ is generated by
\[
\{\nu_R : R\text{ is $\Xr$-irrelevant}\}\cup M_{\Xr}(G).
\]
\end{thm}
 The proof of the previous result needs two preliminary results, 
 The next proposition  follows immediately from  Proposition \ref{p:umindegreeTAR}, which implies that whether or not a set is minimal can be determined from its cardinality and its degree; it is also  stated and proved for $X$-sets as Proposition 2.10 in \cite{XTARiso}.

\begin{prop}\label{minimal-u} Let $X$ be a super $X$-set parameter, and let $G$ and $G'$ be graphs. 
 Suppose $\varphi:\xrtar(G)\to \xrtar(G')$ is a graph isomorphism  such that $|\varphi(S)|=|S|$. Then $\varphi$ maps minimal $X$-sets to minimal $X$-sets (of the same size).
\end{prop}

The next result is established by the proof of Proposition 2.11 in \cite{XTARiso} .
\begin{prop}\label{isomorph-minimal-u} Let $X$ be a super $X$-set parameter, and let $G$ and $G'$ be graphs.  
Suppose $\psi:V(G)\to V(G')$ is a bijection.
\ben[$(1)$]
\item\label{c3}  Suppose $\psi$ maps  $X$-sets of $G$ to  $X$-sets of $G'$. Then the induced mapping $\psi:V(\xtar(G))\to V(\xtar(G'))$ is an isomorphism of $\xtar(G)$ and $\psi(\xtar(G))$.   If every  $X$-set of $G'$ is the image of an $X$-set of $G$, then $\psi:V(\xtar(G))\to V(\xtar(G'))$ is an isomorphism of $\xtar(G)$ and $\xtar(G')$.
\item\label{c2} Suppose $\psi$  maps minimal $X$-sets of $G$ to minimal $X$-sets of $G'$.  Then  $\psi$ maps $X$-sets of $G$ to $X$-sets of $G'$. 
If every minimal $X$-set of $G'$ is the image of a minimal $X$-set of $G$, then  $\psi$ is a bijection from $X$-sets of $G$ to $X$-sets of $G'$. 
\item\label{c4} Suppose $\psi$  maps minimal $X$-sets of $G$ to minimal $X$-sets of $G'$.   Then the induced mapping\break $\psi:V(\xtar(G))\to V(\xtar(G'))$ is an isomorphism of $\xtar(G)$ and $\psi(\xtar(G))$.   If every minimal  $X$-set of $G'$ is the image of a minimal $X$-set of $G$, then $\psi:V(\xtar(G))\to V(\xtar(G'))$ is an isomorphism of $\xtar(G)$ and $\xtar(G')$.\een
\end{prop}

\bpf \blu  
\eqref{c3}: Since $\psi$ maps $X$-sets of $G$ to $X$-sets of $G'$,   $\psi$ induces a bijection between the vertices of $\xtar(G)$ and a subset of the vertices of  $\xtar(G')$ ($X$-sets of $G'$ of the form $\psi(S)$ where $S$ is an $X$-set of $G$). Assume that $S_1,S_2\in V(\xtar(G))$ are adjacent in $\xtar(G)$. Without loss of generality, $|S_1\setminus S_2|=1$. Since $\psi$ is a bijection, $|\psi(S_1)\setminus \psi(S_2)| = 1$. Thus, $\psi(S_1)$ and $\psi(S_2)$ are adjacent in $\xtar(G')$. Hence $\psi$ is an isomorphism from $\xtar(G)$ to $\psi(\xtar(G))$.

\eqref{c2}: 
Let $S\in V(\xtar(G))$ be an $X$-set. There is a minimal $X$-set  $T\subseteq S$ of $G$ and $\psi(T)\subseteq \psi(S)$. Since $\psi(T)$ is a minimal $X$-set  of $G'$, $\psi(S)$ is an $X$-set   of $G'$.

Statement \eqref{c4} is immediate from statements \eqref{c2} and \eqref{c3}
\epf

{\blu
\noi{\bf Proof of Theorem \ref{t:auto-univ}:}
\bpf 
By Theorem \ref{t:irrel-isom} $\nu_R\in\text{aut}(\xrtar(G))$ for every $X$-irrelevant set $R$. 
 By Proposition \ref{isomorph-minimal-u}\eqref{c4}, $\psi\in\text{aut}(\xrtar(G))$ for every  $\psi\in M_{\Xr}(G)$ (since $G'=G$ here, having $\psi$ map minimal $\Xr$-sets to minimal $\Xr$-sets is sufficient). 
 
We now show that $\{\nu_R : R\text{ is $\Xr$-irrelevant}\}\cup M_{\Xr}(G)$ generates aut$(\xrtar(G))$. Let $\varphi$ by an automorphism of $\xrtar(G)$. Suppose first that $V(G)$ is fixed by $\varphi$. Then   $|S|=n-\dist(V(G),S)=n-\dist(V(G),\varphi(S))=|\varphi(S)|$  for each $S\in V(\xrtar(G))$.   By Theorem \ref{isomorph-same-zfs} there exists a bijection $\psi:V(G)\to V(G)$ such that $\psi(S) = \varphi(S)$ for every $\Xr$-set $S$ of $G$.  By  Proposition \ref{minimal-u}, $\psi$ maps minimal $\Xr$-sets to minimal $\Xr$-sets of the same size.  Thus $\psi\in M_{\Xr}$.

Suppose that $V(G)$ is not fixed by $\varphi$. By Theorem \ref{isomorph-isometry-u} and the preceding argument, there exists a bijection $\psi\in M_{\Xr}(G)$ such that $\psi = \nu_{R}\circ\varphi$, where $R=V(G)\setminus \varphi(G)$. Thus, $\varphi = \nu_R^{-1}\circ\psi$.
\epf
}

 Recall that if $X(K_1)=0$, then every isolated vertex is irrelevant (see Remark \ref{no-isol-vtx-2})
\subsection{Connectedness}\label{ss:Xconn}

 A main question in reconfiguration is: For which $k$ is the subgraph of the TAR reconfiguration graph induced by sets of no more than $k$ vertices  connected?  For a super $X$-set parameter $X$, $\xtar (G)$ is always connected (every $X$-set can be augmented one vertex at a time to get $V(G)$).  
 Parameters relating to connectedness were defined for domination and bounds on these parameters were established in \cite{HS14}. The definitions and bounds were extended to original $X$-sets in \cite{PDrecon}. 
 In this section we  further extend  these definitions and results   to  super  $X$-set parameters, showing they remain valid for a larger class of parameters that the isomorphism results. For example, these results apply to connected domination,  which is not robust (see  Section \ref{s:conD}).  
 A method involving twins was used in \cite{XTARiso} to construct examples  with strict inequality in one of the bounds, and we extend this from zero forcing to super $X$-set parameters that satisfy additional conditions on twins.

\begin{defn}  \label{def: k token X}
Suppose $X$ is a  super $X$-set parameter.      The \emph{$k$-token addition and removal (TAR) reconfiguration graph for $X$}, denoted by $\xtar_k(G)$, is the subgraph of $\xtar(G)$ induced by the set of all  $X$-sets of cardinality at most $k$ as its vertex set.

The least $k_0$ such that $\xtar_k(G)$ is connected for all $k\ge k_0$ is denoted by $\xxo(G)$, and the least $k$ such that $\xtar_k(G)$ is connected is denoted by $\ulxo(G)$.
\end{defn}

 Knowing $\xxo(G)$ or $\ulxo(G)$ allows us to work within a smaller TAR reconfiguration graph when modifying one solution to another.
Since  $\xtar_k(G)$ is an induced subgraph of $\xtar(G)$,  many elementary results concerning $\xtar(G)$ also apply to $\xtar_k(G)$ . 


 The next result was established for  original $X$-set parameters in  Proposition 2.8 in \cite{PDrecon}.  Although that result does not include Proposition \ref{basic1}\eqref{h2}, that statement is immediate since adding one vertex at a time to the one minimal $X$-set does not disconnect the graph.  The proof of \cite[Proposition 2.8]{PDrecon} remains valid for any super $X$-set parameter (and the hypothesis that $G$ has no isolated vertices is not needed).

\begin{prop}\label{basic1}
    Let  $X$ be a  super $X$-set parameter and let $G$ be a graph of order $n$.
    \ben[$(1)$]
    \item\label{h1} Then $X(G)\le \ulxo(G)\le\xxo(G)$.
    \item\label{h2} If $G$ has only one minimal $X$-set, then $X(G)=\ulxo(G)=\xxo(G)$.
    \item\label{h3} If $G$ has more than one minimal $X$-set, then 
    $\xbar(G)+1\le \xxo(G)\le\min\{\xbar(G)+X(G),n\}$.
    \item\label{h4} If $G$ has more than one minimum $X$-set, then 
    $X(G)+1\le \ulxo(G)$.
    \een
\end{prop}

{\blu
\bpf
\eqref{h1}: This relationship follows immediately from the definitions.

\eqref{h2}: Building onto one minimal $X$-set does not disconnect the graph, therefore $\xtar_k(G)$ is always connected for $k \geq X(G)$ when there is only one minimal $X$-set.

\eqref{h3}: Let $\hat S\subset V(G)$ be minimal $X$-set with $|\hat S|=\xbar(G)$.  Then $\hat S$ is an isolated vertex of $\xtar_{\xbar(G)}(G)$ (because we can't add a vertex, and removal results in a set that is not an $X$-set).  Thus $\xbar(G)+1\le \xxo(G)$. 
It is immediate from the definition of $\xxo(G)$ that $\xxo(G)\le n$. Suppose that  $\xbar(G)+X(G)< n$ and let $k_0=\xbar(G)+X(G)$. Let $S\subset V(G)$ be a minimal $X$-set of $G$ and $S'\subset V(G)$ be a minimum $X$-set of $G$.    To ensure $\xtar_k(G)$ is connected for all $k\ge k_0$, it is sufficient to show that  every such pair of vertices $S$ and $S'$ is connected in $\xtar_{k_0}(G)$.
Define $S''=S\cup S'$ and observe that $|S''|\le k_o$.  Then each of $S$ and $S'$ is connected by a path to $S''$ by adding one vertex at a time.  Thus $\xxo(G)\le\xbar(G)+X(G)$. 

\eqref{h4}: Each minimum $X$-set is an isolated vertex in $\xtar_{X(G)}(G)$.
\epf
}

 For each of the parameters discussed in Sections \ref{ss:Dom}--\ref{s:vtxcover}, it is easy to find examples of graphs for which $\xbar(G)+1= \xxo(G)$ and $\ulxo(G)=\xxo(G)$.  However, strict inequalities are also possible, and examples are presented for each of the parameters discussed. The next result provides some such examples.

\begin{cor}\label{c:X1}
Let  $X$ be a  super $X$-set parameter. 
 If $G$ has more than one minimal $X$-set and $X(G)=1$, then $\xxo(G)=\xbar(G)+1$ .
\end{cor}

The   next result expands and extends Lemma 4 in \cite{HS14} to super $X$-set parameters; this lemma  established Proposition \ref{basic2}\eqref{conn3} for domination.
 
\begin{prop}\label{basic2}
    Let  $X$ be a  super $X$-set parameter and let $G$ be a graph of order $n$. 
     \ben[(1)]
\item\label{conn1} If for every pair of minimal $X$-sets $M_1$ and $M_2$, there is a path between $M_1$ and $M_2$ in $\xtar_k(G)$, then $\xtar_k(G)$ is connected. 

  \item\label{conn3}  If $k\ge \xbar(G)$ and $\xtar_k(G)$ is connected, then $\xxo(G)\le k$.
\item If  $|M_1\cup M_2|\le k$ for every pair of minimal $X$-sets $M_1$ and $M_2$, then $\xxo(G)\le k$.  \een
\end{prop}
\bpf  Suppose first that for every pair of minimal $X$-sets $M_1$ and $M_2$, there is a path between $M_1$ and $M_2$ in $\xtar_k(G)$.   Given two $X$-sets $S_1,S_2\in V(\xtar_k(G))$, each $S_i$  contains a minimal $X$-set $M_i$. There are paths in $\xtar_k(G)$ from $S_1$ to $M_1$, $M_1$ to $M_2$, and $M_2$ to $S_2$, so $\xtar_k(G)$ is connected.

Now assume    $k\ge \xbar(G)$, $\xtar_k(G)$ is connected, and $\ell>k$. Since $k\ge \xbar(G)$, $\xtar_k(G)$ contains every minimal $X$-set and there is a path between every pair of minimal $X$-sets in $\xtar_k(G)$, which is a subgraph of $\xtar_\ell(G)$.  Thus $\xtar_\ell(G)$ is connected by \eqref{conn1}. 

If   $|M_1\cup M_2|\le k$ for every pair of minimal $X$-sets $M_1$ and $M_2$, then there  is a path through $M_1\cup M_2$ in $\xtar_k(G)$ for every pair of minimal $X$-sets $M_1$ and $M_2$ and $\xxo(G)\le k$.
\epf

 It is sometimes useful to exhibit examples where inequalities are not equalities. The next remark (based on Remark 2.10 in \cite{PDrecon}) presents conditions under which a graph $G$  
satisfies $ \xxo(G)> \xbar(G)+1$.
 
\begin{rem}Let $X$ be a super $X$-set parameter  and let $G$ be a graph  that has at least two minimal $X$-sets. Suppose $G$ has a minimal $X$-set $S$ such that $|S|=\xbar(G)$ and for every $u$ and $v$ such that $u\in S$ and  $v\in V(G) \setminus S$, the set $B=(S\setminus \{u\} )\cup \{v\}$ is not an $X$-set of $G$.  Since there is another minimal $X$-set that is not $S$ and it must be in  another component of $\xtar_{\xbar(G)+1}(G)$, $ \xxo(G)> \xbar(G)+1$.\end{rem}

 The next result is a generalization of Proposition  2.11 in \cite{PDrecon}.  Just as  Proposition  2.11 was used to bound $\pdo(H)$ for an example $H$ such that $\pdo(G)> \pdbar(H)+1$, Proposition \ref{p:part-min-sets} is used in the proof of Proposition \ref{p:upper-con-skew}  to show that $\zso(H(r))>\zsbar(H(r))+1$.  



\begin{prop}\label{p:part-min-sets}
Let $X$ be a  super $X$-set parameter. 
Let  $G$ be a graph such that there is a  partition of the minimal $X$-sets  into two sets $\{S_1,\dots,S_k\}$ and $\{T_1,\dots,T_\ell\}$.   
Let $s=\max_{i=1}^k|S_i|$, $t=\max_{j=1}^\ell|T_j|$, and $p=\max_{i=1,j=1}^{k,\ell}|S_i\cap T_j|$.  Then 
\ben
\item $x_0(G)\ge s+t-p$.
\item If $|M_1\cup M_2|\le s+t-p$ for any minimal $X$-sets $M_1$ and $M_2$, then $x_0(G)= s+t-p$. 
\een
\end{prop}

\bpf 
Let $S \in \{S_1,\dots,S_k\}$ with $|S|=s$, and $T \in \{T_1,\dots,T_\ell\}$ with $|T|=t$.   Since $|S\cap T|\le p$,  $\dist_{\xtar(G)}(S,T)\ge s+t-p$. 
Thus there is no path between $S$ and $T$ in $\xtar_k(G)$ for $k$ with $\max\{s,t\}\le k\le s+t-p-1$.

Now assume $|M_1\cup M_2|\le s+t-p$ for any minimal $X$-sets $M_1$ and $M_2$.  Then $\xxo(G)\le s+t-p$ by Proposition \ref{basic2}.
\epf

 Twin vertices were used in \cite{XTARiso}  to  expand one graph $G$ with $\ulzo(G)< \zzo(G)$ to a family of graphs with that property. This twin method can be extended to super $X$-set parameters that satisfy additional conditions  related to twins. Vertices $u$ and $w$ in a graph $G$ are called \emph{independent twins} if  $N_G(u) = N_G(w)$. A \emph{set of independent twins} is a set  $\{u_1,\dots,u_r\}\subseteq V(G)$ such that $u_i$ and $u_j$ are independent twins  for all $1\le i<j\le r$. For a graph $G$ with  a set of independent twins $T$,  define $G_u = G-u$ for $u\in T$.
The next remark and lemma 
 require a  property   that is parameter specific (this property is established below for skew forcing in Proposition \ref{p:skew-twin}). 

\begin{defn} 
 A super $X$-set parameter $\Xb$ has the \emph{twins property} if it satisfies the following conditions  for every graph $G$ and  set of independent twins $T$ with $|T|\ge 3$: 
 \ben
 \item Any $\Xb$-set must contain at least $|T|-1$ elements of $T$.
 \item If $S_u$ is an $\Xb$-set of $G_u$, then $S=S_u\cup\{u\}$ is an $\Xb$-set of $G$.  
 \item If $S$ is an $\Xb$-set of $G$ and $ u\in S$, then  $S_u=S\setminus\{u\}$ is an $\Xb$-set of $G_u$.
 \een
\end{defn}

 \begin{rem}\label{r:twin} Suppose $\Xb$ is a super $X$-set parameter that satisfies the twins property and $G$ is a graph with a set $T$ of at least three independent twins.  Then there is a bijection between $\Xb$-sets of $G_u$ and $\Xb$-sets of $G$ that contain $u$, and an $\Xb$-set 
 $S_u$ of $G_u$ is minimal if and only if $S=S_u\cup\{u\}$ is a minimal $\Xb$-set of $G$.  Furthermore, $\Xb(G-u)=\Xb(G)-1$ and  $\xbbar(G-u)=\xbbar(G)-1$  for $u\in T$. 
 \end{rem}
 
 The proof of the next result is very similar to that of Lemma 3.8 in \cite{XTARiso} and is omitted.
\begin{lem}\label{l:twin-conn}
Let $\Xb$ be a super $X$-set parameter such that has the twins property. Let $G$ be a graph
that has a set of  twins $T$ with $|T|\ge 3$.  If ${\xbtar}_{k}(G_u)$ is connected for some $u$, then ${\xbtar}_{k+1}(G)$ is connected.
\end{lem}

{\blu  \bpf For $u\in T$, the graphs  $G_u$ are isomorphic and hence the ${\xbtar}_{k}(G_u)$ are isomorphic as well.  Thus if ${\xbtar}_{k}(G_u)$ is connected for some $u\in T$,  then ${\xbtar}_{k}(G_w)$ is connected for all $w\in T$.  For  $S\subseteq V(G)$ with $u\in S$, let $S_u=S\setminus \{u\}$.  

Since $r\ge 3$ and $\Xb$ has the twins property,  any $\Xb$-set $S$ of $G$ must contain at least two vertices in $T$, say $u$ and $w$.  A set  $S\subseteq V(G) $  that contains $u,w$  is an $\Xb$-set of $G$  if and only if  $S \setminus \{u\}$ is an-$X$-set of $G_u$ and  $S \setminus \{w\}$ is an-$X$-set of $G_w$.  

Let $S, S'$ be  two $X$-sets  of $G$ of size $k+1$ or less.  Since $r\ge 3$  and each can omit at most one  vertex in $T$, their intersection must contain at least one  $u \in T$. Then $S\setminus \{u\}, S'\setminus \{u\} $ are $\Xb$-sets for $G_u$. Since by hypothesis ${\xbtar}_{k}(G_u)$ is connected, that means that there is a path between $S\setminus \{u\}$ and  $S'\setminus \{u\} $ in ${\xbtar}_{k}(G_u)$ and hence a path between $S, S'$ in ${\xbtar}_{k+1}(G)$.
\epf}

The proof of the next result uses ideas from Proposition 3.9 
in \cite{XTARiso}, but that result is specific to one graph family (in addition to being stated for zero forcing), so we provide the brief proof here.  

\begin{prop}\label{p:xso-lowxso} Suppose $\Xb$ is a super $X$-set parameter that has the twins property and $G$ is a graph such that  $G$ has a set $T$ of $t\ge 2$ twin vertices,      $\Xb(G)<\xbbar(G)$, and  $\xbtar(G)_{\Xb(G)+1}$ is connected. 
 Define  $G(t)=G$ and for $r\ge t$, construct $G(r+1)$ from $G(r)$ by adding one additional independent twin of a vertex in $T$.  Then $\Xb(G(r))=\Xb(G)+(r-t)$, $\xbbar(G(r))=\xbbar(G)+(r-t)$, $\ulxbo(G(r))=\Xb(G(r))+1$, and $\xbo(G(r))>\ulxbo(G(r))$.
\end{prop} 

\bpf  The proof is by induction.  The base case for induction is $r=t$ and is assumed.  For $r>t$, apply Remark \ref{r:twin} to show that $G(r)$ has minimal $X$-sets of sizes $\Xb(G(r))=\Xb(G)+(r-t)$ and $\xbbar(G(r))=\xbbar(G)+(r-t)$; let $k=\Xb(G)+(r-t)$ and $\ell=\xbar(G)+(r-t)$.  Apply Lemma \ref{l:twin-conn} to show ${\xbtar}_{k+1}(G(r))$ is connected. Since $\xbbar(G(r))>\Xb(G(r))$, $G(r)$ has at least two minimal $\Xb$-sets, implying ${\xbtar}_{k}(G(r))$ and ${\xbtar}_{\ell}(G(r))$ are disconnected by Proposition \ref{basic1}.  Thus
$\ulxbo(G(r))=k+1$ and $\xbo(G(r))>\ulxbo(G(r))$.  
\epf

 Token jumping (TJ) reconfiguration (also called token exchange), which  involves exchanging a vertex between sets of vertices of the same size, has been studied for some cohesive vertex set parameters, particularly for minimum $X$-sets  (see, e.g., \cite{GHH} for standard zero forcing and \cite{PDrecon} for power domination).
Kami\'nski, Medvedev, and Milani\v c study token jumping reconfiguration of independent sets and establish the equivalence of the connectedness of the token jumping reconfiguration graph of independent sets of $G$ of size $k$ and the connectedness of the independence $(k-1)$-TAR graph of $G$ in \cite{KMM12}.
Here we show this equivalence extends to   super $X$-set parameters (and later present the more obvious extension to sub $Y$-set parameters in Proposition \ref{p:YvsetTJ}).  

\begin{defn}
Let $W$ be a cohesive vertex set parameter and let $G$ be a graph.  The  \emph{$k$-token jumping} or  \emph{$k$-TJ reconfiguration graph of $G$ for $W$} takes as vertices  the $W$-sets of size $k$, with an edge between two $W$-sets $R_1$ and $R_2$ if and only if $R_2$ can be obtained from $R_1$ by 
exchanging exactly one vertex.  
\end{defn}

\begin{prop} \label{p:vsetTJ}
Let $X$ be a super $X$-set parameter, let $G$ be a graph, and let  $S_1$ and $S_2$ be $X$-sets of $G$  with $|S_i|=k, i=1,2$. Then there is a path between  $S_1$ and $S_2$ in the $k$-TJ reconfiguration graph of $G$  if and only if there is a path between  $S_1$ and $S_2$ in $\xtar_{k+1}(G)$.     
\end{prop}
\bpf If $(S_1=R_1,R_2,\dots,R_{r-1},R_r=S_2)$ is a path in  $k$-TJ reconfiguration graph of $G$, then $(S_1=R_1,R_1\cup R_2,R_2, R_2\cup R_3, R_3,\dots,R_{r-1},R_{r-1}\cup R_r,R_r=S_2)$ a path between  $S_1$ and $S_2$ in $\xtar_{k+1}(G)$.

We show that any path in between $S_1$ and $S_2$ can be replaced by a path that uses only $X$-sets of sizes $k$ and $k+1$.  Such a path starts with a $k$-set and must alternate additions and deletions.  Each addition-deletion pair  can then be replaced by an edge in the $k$-TJ graph of $G$. 

So suppose $(S_1=R_1,R_2,\dots,R_{r-1},R_r=S_2)$ is a path  between  $S_1$ and $S_2$ in $\xtar_{k+1}(G)$. If $|R_i|\ge k$ for $i=1,\dots,r$ then there is nothing to show.  So assume that is not the case and let $j$ be an   index such that $|R_j|$ is minimized. Necessarily $|R_j|=|R_{j-1}|-1$,  $R_j$ is obtained from $R_{j-1}$ by deleting some vertex $w_{j-1}$, $|R_{j+1}|=|R_{j}|+1$, and $R_{j+1}$ is obtained from $R_{j}$ by adding some vertex $u_{j}$.  Then replace $R_j$ by $R'_j=R_{j-1}\cup\{u_j\}$ so $R'_j$ is obtained from $R_{j-1}$ by adding $u_j$ and $R_{j+1}$ is obtained from $R_j$ by deleting $w_{j-1}$.  Since $X$ is a super $X$-set parameter and $R_{j-1}$ is an $X$-set, necessarily $R'_j$ is an $X$-set. This process can be repeated as needed to obtain a path in which all $X$-sets have size $k$ or $k+1$. 
\epf

\subsection{$X$-TAR graph uniqueness}\label{ss:uni}
Let $X$ be a robust $X$-set parameter.  If $X(K_1)=1$, for a graph $H$ with no isolated vertices, we say its  TAR  graph  is \emph{unique} if  $\xtar(G)\cong\xtar(H)$  implies $G\cong H$  for any graph $G$ with no isolated vertices. If $X(K_1)=0$, for a graph $H$, we say its  TAR  graph  is \emph{unique} if  $\xtar(G)\cong\xtar(H)$  implies $G\cong H$  for any graph $G$.   The concept of  unique TAR graphs was introduced in \cite{PDrecon} for power domination. Certainly whether a graph has a unique TAR graph is parameter specific, and is discussed for various parameters in Sections \ref{s:Dom-PD-Z}-- \ref{s:vtxcover}.  For power domination and standard zero forcing, examples of graphs with unique TAR graphs were presented in \cite{PDrecon} and \cite{XTARiso}.  The fact that the order of $G$,  $X(G)$, and $\xbar(G)$ can all be determined from $\xtar(G)$ was used in \cite{PDrecon}, and the fact that $\xtar(G)\cong\xtar(H)$ implies we can relabel so $G$ and $H$ have the same $X$-sets was essential to the uniqueness results in \cite{XTARiso}. This ability to assume that $G$ and $H$ have the same $X$-sets is true for robust $X$-set parameters, and we use this to  present examples of unique TAR graphs for domination, PSD forcing,  skew forcing, and vertex covers in Sections \ref{ss:Dom}, 
\ref{s:PSDunique},  \ref{ss:skewex}, and \ref{ss:vcex}.  We also present data for the frequency of unique TAR graphs among small graphs in those sections.
 We point out that vertex cover TAR graphs behave quite differently from the TAR graphs of the other parameters discussed here: it is shown in Proposition \ref{p:VCunique} that every vertex cover TAR graph is unique.
Which specific  graphs have unique TAR graphs naturally depends on the parameter, but there are some graph families that have unique TAR graphs for many of the parameters discussed, such as complete graphs  and complete bipartite graphs. 
Even when the same graph family has unique TAR graphs for multiple parameters, the proofs are usually parameter-specific (although some ideas recur).

\subsection{Hamilton paths and cycles and cut-vertices in TAR graphs}\label{ss:Ham}
 
The study of whether a TAR graph has a Hamilton path was initiated in \cite{WIGArecon} for domination.  As shown there, no domination TAR graph can have a Hamilton cycle (because any cycle in a bipartite graph is even and there are an odd number of dominating sets), but for some of the other parameters, such as skew forcing, it is possible to have a Hamilton cycle.

\begin{rem}\label{r:QnHam}  Let $X$ be a super $X$-set parameter.   It is well-known that a  hypercube of dimension at least two has a Hamilton cycle.  Thus $\xtar(G)$ has a Hamilton cycle whenever $\xtar(G)\cong Q_n$, which happens if the empty set is an $X$-set. The graph $Q_n-v$ for $v\in V(Q_n)$ is realized as a TAR graph whenever any one vertex of $G$ is an $X$-set. Since deleting a vertex from a cycle leaves a path and $Q_n-v$ has an odd number of vertices, $Q_n-v$ has an Hamilton path but not a Hamilton cycle.  
\end{rem}

 Lemma 6.7 in \cite{WIGArecon} shows that $V(G)$ 
is the only possible cut-vertex of  the domination TAR graph of $G$ 
  for any graph $G$. We extend this from domination to super $X$-set parameters. 
 
\begin{prop}\label{p:cut-vtx}
Let $X$ be a super $X$-set parameter and let $G$  be a  graph of order $n$. 
If $S$  is a cut-vertex of  $\xtar(G)$, then $S=V(G)$. 
\end{prop}

\bpf 
 We prove the contrapositive.  Assume that $S\ne V(G)$. If $S$ has exactly one neighbor, then $S$ is not a cut-vertex. So, suppose that $S$ has distinct neighbors $A,B\subseteq V(G)$ in $\xtar(G)$.  It suffices to show that there is a path in $\xtar(G)-S$ between $A$ and $B$.  
Since $S\ne V(G)$, there is some vertex $w$ of $G$ such that $w\not\in S$. Define $S'=S\cup\{w\}$, $A'=A\cup\{w\}$ and $B'=B\cup\{w\}$.  Note that $S'\ne S$.   If $w\not\in A$ then $ A\ne A',S'$ and  $A'\ne S'$, and similarly for $B$.  
Thus if  $w\not\in A$ and $w\not\in B$, then $(A,A',S', B',B)$ is a path from $A$ to $B$ in $\xtar(G)-S$.  If  $w\in A$ and  $w\not\in B$, then $A=A'=S'$ and $(A,B',B)$ is a path from $A$ to $B$ in $\xtar(G)-S$, and similarly when $w\not\in A$ and  $w\in B$.  It is not possible to have $w\in A$ and $w\in B$.
\epf

Observe that if $X$ is a connected $X$-set parameter, $G$ is a graph of order $n\ge 2$, and $\xbar(G)=n-1$, then $V(G)$ is a cut-vertex of $\xtar(G)$ because any minimal $X$-set with $n-1$ vertices is a connected component in $\xtar(G)-V(G)$ (and $G$ has more than one minimal $X$-set).

\section{TAR reconfiguration results for domination, power domination, and standard zero forcing }\label{s:Dom-PD-Z}

 TAR graphs have  been studied previously for domination, power domination, and zero forcing.  In this section we review these results and provide some additional material.
 
\subsection{Domination TAR graphs}\label{ss:Dom}
 As noted in \cite{XTARiso}, domination is an original $X$-set parameter, so Theorem \ref{t:main} applies. 
 We summarize known results and provide some additional material related to domination TAR graphs, specifically uniqueness, Hamiltonicity, and connectedness results.

 The \emph{domination TAR graph}  of a  base graph  $G$ is denoted by  $\dtar(G)$ and the upper domination number of $G$ is denoted here by $\dbar(G)$ (note that in the literature this parameter is usually dented by $\Gamma(G)$).

\begin{ex}\label{ex:KnDom}
    Since $K_n$ is the only graph of order $n$ for which any one vertex dominates, $K_n$ is the only graph such that $\dtar(K_n)\cong Q_n-v$ for $v\in V(Q_n)$ and $\dtar(K_n)$ is unique.    As was noted in \cite{WIGArecon}, $\dtar(K_n)$ has a Hamilton path but not a Hamilton  (see also Remark \ref{r:QnHam}).
\end{ex}

\begin{ex}\label{ex:KabDom}
 Let $2\le p\le q$ and let $A$ and $B$ be the partite sets of $K_{p,q}$. Then the minimal dominating sets of $K_{p,q}$ are  $A$, $B$ and $\{a_i,b_j\}$ where $a_i\in A$ and $b_j\in B$.   Thus $\gamma(K_{p,q})=2$ and $\dbar(K_{p,q})=q$. 
\end{ex}

  Complete biparite graphs are discussed further in Example \ref{ex:KabDom2} as well as in  the next result.
 
\begin{prop}
    For $1\le p\le q$, $\dtar(K_{p,q})$ is unique unless $q=p\ge 3$.  For $p=q\ge 3$, if $\dtar(G)\cong\dtar(K_{p,p})$, then $G\cong \dtar(K_{p,p})$ or $G\cong \dtar(K_p\Box K_2)$.
\end{prop}
\bpf   Suppose $G$  has no isolated vertices and $\dtar(G)\cong \dtar(K_{p,q})$.  This implies $G$ has the same dominating sets as $K_{p,q}$ (possibly after relabeling).  If $p=1$, then  $\gamma(G)=1$ and $\dbar(G)=q=|V(G)|-1$, which implies $G\cong K_{1,q}$.  Thus $\dtar(K_{1,q})$ is unique, and we assume $q\ge p\ge 2$. 

We begin by  showing that this implies $G[B]\cong qK_{1}$ or $G[B]\cong K_{q}$.  Suppose  that $G[B]\not\cong K_{q}$, i.e., there exist $k,\ell$ such that $b_k$ and $b_\ell$ are not adjacent.  Since $\{a_i,b_k\}$ is a dominating set for $i=1,\dots,p$, $A\subseteq N_G[b_\ell]$ (and similarly $A\subseteq N_G[b_k]$).  Suppose there is an edge in $G[B]$, i.e., there exist $i,j$ such that $b_i$ and $b_j$ are adjacent.  Necessarily $\{i,j\}\ne\{k,\ell\}$, so without loss of generality $i\ne k,\ell$.  Then $B\setminus \{b_i\}$ is a dominating set, contradicting  the minimal dominating sets of $G$.  
Thus $G[B]\cong qK_1$  or $G[B]\cong K_q$. 
Similarly, $G[A]\cong pK_1$ or $G[A]\cong K_p$.

Now suppose $G[B]\cong qK_1$; we show this implies $G\cong K_{p,q}$. By the argument above, $G[B]\cong qK_1$ implies  every vertex of $B$ is adjacent to every vertex of $A$. If $a_i$ and $a_j$ were adjacent, then $A\setminus \{a_i\}$ would be a dominating set, which is a contradiction.  Thus $G[A]\cong pK_1$ and $G\cong K_{p,q}$.  

The previous argument applies to $G[A]\cong pK_1$ also, so the only remaining case is $G[B]\cong K_q$ and $G[A]\cong K_p$, and we assume this.   We show that for every $k=1,\dots,q$, there must exist $\ell\in\{1,\dots,p\}$ such that $N_G[a_\ell]\cap B=\{b_k\}$.  Once this is established, necessarily $q=p$ and $G\cong K_p\Box K_2$. So suppose not, i.e., there is some $b_k$ such that $N_G[a_i]\cap B\ne \{b_k\}$  for every $i=1,\dots,p$.  Then $B\setminus \{b_k\}$ is a dominating set because $b_k$ is dominated by any other vertex of $B$ and $A$ is dominated by $B$. Since this is a contradiction, $G\cong K_p\Box K_2$.

Let $G\cong K_p\Box K_2$ with $p\ge 2$, and denote the vertices of one copy of $K_p$ by $A$ and the other by $B$,  Then minimal dominating sets of $G$ are  $A$, $B$ and $\{a_i,b_j\}$ where $a_i\in A$ and $b_j\in B$.  Thus $\dtar(G)\cong  \dtar(K_{p,q})$.  Note that for $p\ge 3$, $K_p\Box K_2\not\cong K_{p,p}$.  Since  $K_2\Box K_2\cong K_{2,2}$, $\dtar(K_{2,2})$ is unique. 
\epf

Table \ref{t:unique-tar-data-D} shows the number of (nonismorphic) graphs without isolated vertices of order at most  eight  that have  unique  domination TAR graphs (this data was computed in \cite{sage:recon}).\begin{table}[!h]
\begin{center}\begin{tabular}{ |c|c|c|c|c|c|c|c|c| } 
 \hline\#  vertices in $G$  & 2 & 3 & 4 & 5 & 6 & 7  & 8\\
\hline\# graphs with  unique $\dtar(G)$  & 1 & 2 & 5 &14 & 55 & 348 & 4275 \\
  \hline\# graphs with no isolated vertices  & 1 & 2 & 7 & 23 & 122 & 888 & 11302\\
  \hline ratio (\# unique/\# no isolated) &  1 & 1 & 0.7143 & 0.6087 &  0.4508 &  0.3919 & 0.3783 \\
 \hline\end{tabular}\vspace{-9pt}
\end{center}
    \caption{ Number of graphs having a unique domination TAR graph for graphs of  small order}
    \label{t:unique-tar-data-D}
\end{table}

 The question of whether $\dtar(G)$ has a Hamilton path is resolved for complete graphs and 
complete bipartite graphs in \cite{WIGArecon}, and for trees and cycles in \cite{WIGArecon2}.

\begin{thm}\label{t:KpqHdom}{\rm\cite{WIGArecon}}
    For $1\le p\le q$, $\dtar(K_{p,q})$ has a Hamilton path if and only if at least one of $p$ or $q$ is odd.
\end{thm}

 \begin{thm}\label{t:tree-dom}{\rm\cite{WIGArecon2}} For any tree $T$, $\dtar(T)$ has a Hamilton path.
\end{thm}

\begin{thm}\label{t:cycle-dom}{\rm\cite{WIGArecon2}} For  $n\ge 3$, $\dtar(C_n)$ has a Hamilton path if and only if n is not
a multiple of $4$.
\end{thm}

 Haas and Seyffarth (\cite{HS14}) were the first to investigate the connectedness of $\dtar_k(G)$.  The least $k_0$ such that $\dtar_k(G)$ is connected for all $k\ge k_0$ is denoted by $\ddo(G)$,  and the least $k$ such that $\dtar_k(G)$ is connected is denoted by $\uldo(G)$. 
 They established the bounds for $\ddo$ found in Proposition \ref{basic1}. It was also shown there that if $G$ is a graph of order $n$ with at least two disjoint edges, then $\dtar_{n-1}(G)$ is always connected. They  found graph families for which $\ddo$ is equal to the lower bound, such as the one listed next.  

\begin{thm}\label{dbarplus1}{\rm\cite{HS14}} For any nontrivial bipartite or nontrivial chordal graph, $\dtar_{\dbar(G)+1}(G)$ is connected. 
\end{thm}

Conditions that guarantee a graph will satisfy $\ddo(G)=\dbar(G)+1$  were presented in \cite{HS17}.
The first examples of a graph family for which $\dtar_{\dbar(G)+1}$ is not connected is due to Suzuki, Mouawad and Nishimura in \cite{SMN14}. They constructed a family of graphs $G_{(b,d)}$ from $d+1$ cliques of size $b$ with additional edges such that $\ddo(G_{(d,b)})=\dbar(G_{(d,b)})+2$. The next example shows a planar graph in this family.

\begin{figure}[!h]
\centering
\scalebox{1}{
\begin{tikzpicture}[scale=1.3,every node/.style={draw,shape=circle,outer sep=2pt,inner sep=1pt,minimum size=.2cm}]		
\node[fill=none]  (1) at (-1,1) {$1$};
\node[fill=none]  (2) at (0,0) {$2$};
\node[fill=none]  (3) at (0,2) {$3$};
\node[fill=none]  (4) at (0.5,1) {$4$};
\node[fill=none]  (5) at (1.5,0) {$5$};
\node[fill=none]  (6) at (1.5,2) {$6$};
\node[fill=none]  (7) at (2,1) {$7$};
\node[fill=none]  (8) at (3,0) {$8$};
\node[fill=none]  (9) at (3,2) {$9$};
		
\draw[thick] (1)--(2)--(3)--(1)--(4)--(5)--(6)--(4)--(7)--(8)--(9)--(7);
\draw[thick] (2)--(5)--(8);
\draw[thick] (3)--(6)--(9);
\end{tikzpicture}}
\caption{ The graph $G=K_3\Box P_3$, which has $\gamma(G)=\dbar(G)=3$ and  $\ddo(G)=5>\dbar(G)+1$ \label{f:K3boxP3}}\vspace{-10pt}
\end{figure}

\begin{ex}
The graph $G_{(2,3)}$ in \cite{SMN14} is isomorphic to $G=K_3 \Box P_3$.    With the vertices labelled as shown in Figure \ref{f:K3boxP3}, the minimal dominating sets are \{4, 5, 6\}, \{1, 2, 9\},
 \{2, 6, 9\},
 \{3, 5, 9\},
 \{2, 4, 8\},
 \{3, 4, 8\},
 \{2, 5, 9\},
 \{2, 6, 7\},
 \{3, 5, 8\},
 \{1, 6, 8\},
 \{1, 5, 7\},
 \{1, 8, 9\},
 \{1, 5, 9\},
 \{2, 5, 8\},
 \{2, 5, 7\},
 \{1, 4, 9\},
 \{1, 6, 9\},
 \{2, 7, 9\},
 \{3, 6, 8\},
 \{1, 5, 8\},
 \{3, 6, 9\},
 \{1, 4, 7\},
 \{1, 6, 7\},
 \{2, 3, 7\},
 \{2, 6, 8\},
 \{3, 6, 7\},
 \{1, 3, 8\},
 \{3, 4, 9\},
 \{3, 7, 8\},
 \{1, 4, 8\},
 \{2, 4, 9\},
 \{3, 4, 7\},
 \{2, 4, 7\},
 \{3, 5, 7\} (this is documented in \cite{sage:recon}).  Observe that a  minimal dominating set is a set of three vertices that is the vertices of middle $K_3$ (the set of vertices $\{4,5,6\}$) or  has exactly one element in each $K_3$. Thus  $\gamma(G)=\dbar(G)=3$ and $\ddo(G)=5>\dbar(G)+1$. 
\end{ex}

For any positive integer $k$, a \emph{$k$-tree} is a graph that can be constructed from $K_{k+1}$ by repeatedly (possibly zero times) adding a new vertex and joining it to an existing $k$-clique.  Thus, a $1$-tree is a tree, and a tree of order two or more is a 1-tree.  
For a graph $G$, the minimum $k$ such that $G$ is a subgraph of some $k$-tree is called the \emph{tree-width} of $G$.  A graph is \emph{$b$-partite} if 
the  vertex set can be partitioned into sets $V_1,\dots,V_b$  such that  for every edge the end points are in distinct partite sets.

\begin{thm}\label{dkupperbound}{\rm\cite{SMN14}}
    For every integer $b \ge 3$ there exists an infinite family of graphs of tree-width $2b-1$ such that for each $G$ in the family $\dtar_{\dbar(G)+1}(G)$ is not connected. There also exists an infinite family of $b$-partitie graphs such that $\dtar_{\dbar(G)+1}$ is not connected for any $G$ in the family.
\end{thm}

In \cite{MTR19}, Mynhardt, Teshima and Roux  constructed graphs with arbitrary domination number and   arbitrary upper domination number  at least three and greater than domination number that  realize the upper bound for $\ddo(G)$, as stated in the next theorem.  This construction was built on another construction in \cite{MTR19} that realizes $\ddo(G)=\gamma(G)+\dbar(G)-1$. 

\begin{thm}\label{ddomaximum}{\rm\cite{MTR19}}
For each integer $k \ge 3$ and each integer $r$ such that $1 \leq r \leq k-1$, there is a graph $Q_{k,r}$ such that $\gamma(Q_{k,r})=r$, $\dbar(Q_{k,r})=k$ and $\ddo(Q_{k,r})=k+r=\gamma(Q_{k,r})+\dbar(Q_{k,r})$.
\end{thm}

 The graph $K_n$ has $\ddo(K_n)=\uldo(K_n)$. The next graph is an example where $\uldo(G)< \ddo(G)$.

\begin{ex}\label{ex:KabDom2}
 Let $4\le p\le q$ and let $A$ and $B$ be the partite sets of $K_{p,q}$. As noted in Example \ref{ex:KabDom}, the minimal dominating sets of $K_{p,q}$ are  $A$, $B$ and $\{a_i,b_j\}$.  Thus $\uldo(K_{p,q})=3$ and $\ddo(K_{p,q})= q+1\ge 5$.
\end{ex}

For more information on domination TAR reconfiguration see the survey by Mynhardt and Nasserasr \cite{MS20}  and the references therein, as well as to \cite{HS14, HS17, MTR19, SMN14}.

\subsection{Power domination TAR graphs}\label{ss:PD}

Token addition and removal reconfiguration  for power domination and  token jumping reconfiguration of minimum power dominating sets were studied  in \cite{PDrecon}, where  the first definition  of an original $X$-set property was given and  results  were established for  TAR graphs of original $X$-set parameters.  As noted in \cite{XTARiso}, power domination is an original $X$-set parameter, so Theorem \ref{t:main} applies. 
Here we summarize parameter-specific results for power domination TAR graphs from \cite{PDrecon} and \cite{XTARiso}, and present some additional information.  The \emph{power domination TAR graph}  of a  base graph  $G$ is denoted by  $\pdtar(G)$ and the upper power domination number of $G$ is denoted by $\pdbar(G)$.

\begin{ex}
It is shown in \cite[Section 3.3]{PDrecon} that $\pdtar(K_{1,q})$ is unique for $q\ge 3$ but  $\pdtar(K_{2,q})$ is not for $q\ge 2$.   
\end{ex}
It is also proved in \cite{PDrecon} that $\pdtar(K_{3,3})$ is unique, and established computationally that $\pdtar(K_{3,4})$, $\pdtar(K_{3,5})$, $\pdtar(K_{4,4})$ and $\pdtar(K_{4,5})$ are unique.  It is conjectured there that $\pdtar(K_{p,q})$ is unique for $q\ge p\ge 3$. 
Table \ref{t:unique-tar-data-PD} shows the number of (nonisomorphic) graphs without isolated vertices of order at most  eight  that have  unique  power domination TAR graphs (this data was computed in \cite{sage:recon}).\begin{table}[!h]
\begin{center}
\begin{tabular}{ |c|c|c|c|c|c|c|c|c| } 
 \hline
\#  vertices in $G$  & 2 & 3 & 4 & 5 & 6 & 7 
 & 8\\
\hline
\# graphs with  unique $\pdtar(G)$  & 1 & 0 & 3 &4 & 13 & 25 &  
79\\
  \hline
\# graphs with no isolated vertices  & 1 & 2 & 7 & 23 & 122 & 888 & 11302\\
 \hline
 ratio (\# unique/\# no isolated) &  1 & 1 & 0.4286 & 0.1739 &  0.1066 &  0.0282 & 0.00699\\
 \hline\end{tabular}\vspace{-9pt}
\end{center}
    \caption{Number of graphs having a unique power domination TAR graph for graphs of  small order}
    \label{t:unique-tar-data-PD}
\end{table}

 Next we examine the Hamiltonicity or lack thereof for power domination TAR graphs,  starting with examples from \cite{PDrecon}.
\begin{ex}  If $G$ is a graph such that any one of its vertices is a power dominating set, then $\pdtar(G)$ is isomorphic to an $n$-dimensional hypercube with one vertex (the empty set) deleted.  Examples of such base graphs include $C_n, P_n, K_n,$ and the wheel $W_n=C_{n-1}\vee K_1$ \cite{PDrecon}.  Thus the power domination TAR graphs $\pdtar(G)$ for $G=C_n, P_n, K_n,W_n$ each have a Hamilton path but not a Hamilton cycle.
\end{ex}

\begin{ex}
The graph $G(r,1)$ is constructed by adding one leaf  vertex $\ell$ to $K_r$ 
and it is   shown in \cite[Proposition 2.16]{XTARiso} that the vertex $\ell$ is $\pd$-irrelevant. Thus $\pdtar(G(r,1))\cong (Q_r-v)\Box K_2$ where $v\in V(Q_r)$, which has a Hamilton cycle because $Q_r-v$ has a Hamilton path.
\end{ex}

Finally, we  present connectedness examples.  The least $k_0$ such that $\pdtar_k(G)$ is connected for all $k\ge k_0$ is denoted by $\pdo(G)$  and the least $k$ such that $\pdtar_k(G)$ is connected is denoted by $\ulpdo(G)$. There are many examples of graphs $G$ where $\pdo(G)=\pdbar(G)+1$ and $\pdo(G)=\ulpdo(G)$, including any graph $G$ of order $n$ that has $\pdtar(G)\cong Q_n-v$  for  $v\in V(Q_n)$. An example from \cite{PDrecon} that shows $\pdo(G)>\pdbar(G)+1$ is possible and a new example that shows $\pdo(G)>\ulpdo(G)$ is possible.  
For any integer $n\geq 3$, let $G_n=(V_n,E_n)$ be the graph defined as follows:
 $V_n=T_n\cup \big( \cup_{i=1}^n S_{n,i}\big)$ where $T_n=\{u_1,\ldots, u_{n-1}\}$ and $S_{n,i}=\{v^i_1,\ldots , v^i_n\}$, for every integer $i$,  $1\leq i\leq n$.
To define $E_n$: For integers $i=1,\ldots, n$ and $j=1,\ldots , n-1$, $N_{G_n}(v^i_j)=\{u_j\}\cup \{v^i_p: 1\leq p\leq n, p\not=j\}$; $N_{G_n}(v^i_n)=\{v^i_1,\ldots , v^i_{n-1}\}$; $N_{G_n}(u_j)=\{v^r_j: 1\leq r\leq n\}$ \cite[Definition 3.11]{PDrecon}.

\begin{thm}{\rm\cite[Theorem 3.12]{PDrecon}}
For $n\ge 3$, the graph $G_n$ has $\pdbar(G_n)=\pd(G_n)=n-1$ and $\ulpdo(G_n)=\pdo(G_n)=2n-2>\pdbar(G_n)+1=\pd(G_n)+1$.
\end{thm}


  Next we present an example of a graph $G$ having $\pdo(G)>\ulpdo(G)$.
For $H$ a   connected graph of order at least two and $q\ge 3$, define the graph $K^{2,q}(H)$ to be the graph obtained from $H$  deleting each edge of $H$  and adding $q$ vertices each adjacent exactly to the endpoints of the edge that was deleted.  Keep the vertex labels of $H$ in $K^{2,q}(H)$, so $V(H)\subseteq V(K^{2,q}(H))$ (this extends the definition of  a family of graphs introduced in \cite{PDrecon} in the study  of token exchange reconfiguration). 
The  next result was established for $q=3$ in Lemma 5.6 of \cite{PDrecon} and the proof remains valid for $q\ge 3$.

\begin{lem}\label{vc-K2qH} Let $H$ be a  connected graph of order at least two and let $q\ge 3$.   If $S$ is a minimum power dominating set  of $K^{2,q}(H)$, then  $S\subseteq V(H)$. A set $S\subseteq V(H)$ is a  minimum vertex cover of  $H$ if and only if $S$ is a minimum power dominating set of $K^{2,q}(H)$.
 \end{lem}
\bpf \blu
Denote the vertices of $H$ by $v_1,\dots,v_n$ and note that $\{v_1,\dots,v_n\}\subseteq V(K^{2,q}(H))$.   Suppose first that $S$ is a power dominating set of $K^{2,q}(H)$. For every  $i,j\in\{1,\dots,n\}$ with $v_iv_j\in E(H)$,  one of $v_i$ or $v_j$ must be in $S$ or at least $q-1$ of the $q$ vertices in $N_{K^{2,q}(H)}(v_i) \cap N_{K^{2,q}(H)}(v_j)$ must be in $S$.    Observe that $q-1\ge 2$.
Thus in order for $S$ to be a minimum power dominating set, $S$  must contain one of $v_i$ or $v_j$ for every $v_iv_j\in E(H)$ (and no vertices in $V(K^{2,q}(H))\setminus V(H)$).  Thus $S$ is a vertex cover  of $H$.

Now suppose $S$ is a minimum vertex cover of $H$ and let $v_iv_j\in E(H)$ with $v_i\in S$. Then the set of $q$ vertices in $N_{K^{2,q}(H)}(v_i) \cap N_{K^{2,q}(H)}(v_j)$ is dominated and any one of these vertices can force $v_j$ if $v_j\not\in S$.  Thus $S$ is a power dominating set of $K^{2,q}(H)$. Since every minimum  power dominating set of $K^{2,q}(H)$ is a subset of $V(H)$, $S$ is a minimum  power dominating set of $K^{2,q}(H)$.
 \epf

\begin{prop} For $r\ge 3$ and $q\ge 4$, $\pd(K^{2,q}(K_r))=r-1$, $\pdbar(K^{2,q}(K_r))\ge q\lp\frac{r(r-1)}2\rp-(r-1)$,
   $\ulpdo(K^{2,q}(K_r))=r$, and $\pdo(K^{2,q}(K_r))\ge q\lp\frac{r(r-1)}2\rp-(r-1)+1 > r=\ulpdo(K^{2,q}(K_r))$.
\end{prop}
\bpf
By Lemma \ref{vc-K2qH}, $\pd(K^{2,q}(K_r))=r-1$ and any set of $r-1$ vertices of $K_r$ is a minimum power dominating set.  Since each of these is adjacent to the vertex $V(K_r)$ in  $\pdtar_r(K^{2,q}(K_r))$ and the only minimal power dominating sets of at most $r$ vertices are the minimum power dominating sets, $\pdtar_r(K^{2,q}(K_r))$ is connected and $\ulpdo(K^{2,q}(K_r))=r$.

Next we exhibit a minimal power dominating set containing  of $q\lp\frac{r(r-1)}2\rp-(r-1)$ vertices, which shows that $\pdbar(K^{2,q}(K_r))\ge q\lp\frac{r(r-1)}2\rp-(r-1)$ and $\pdo(K^{2,q}(K_r))\ge q\lp\frac{r(r-1)}2\rp-(r-1)+1$.  It can be verified algebraically  that $q\lp\frac{r(r-1)}2\rp-(r-1)+1 > r$ for $r\ge 3$ and $q\ge 4$.  

Denote the vertices of $K_r$ by $v_1,\dots,v_r$ and denote the $q$ vertices of $K^{2,q}(K_r)$ obtained by replacing the edge $v_iv_j$ by $u^{i,j}_k$ for $k=1,\dots,q$.  Define $S$ to be the set of the first $q-1$ degree-two vertices associated with each of the $r-1$ edges incident with $v_1$ and all $q$ degree-two vertices associated with each of the  edges not incident with $v_1$, i.e.,  $S=\{u^{1,j}_k:k=1,\dots,q-1,j=2,\dots,r\}\cup\{u^{i,j}_k:k=1,\dots,q,j=i+1,\dots,r,i=2,\dots,r\}$.  Observe that $|S|=q\lp\frac{r(r-1)}2\rp-(r-1)$.  Starting with the vertices in $S$ blue, after the first (domination) step all vertex are blue except $v^{1,j}_q$ for $j=2,\dots, r$.  Then $v_j$ can observe $v^{1,j}_q$ in the next round.
\epf

\subsection{Standard zero forcing TAR graphs}\label{ss:ZF}
TAR reconfiguration was studied for standard zero forcing in \cite{XTARiso}.  The main isomorphism results were proved for original $X$-sets there.  Both those results and results for original $X$-set parameters were applied to standard zero forcing. 
Additional parameter-specific results for standard zero forcing TAR graphs are also presented in \cite{XTARiso}.   We briefly summarize the latter results here, and also examine the Hamiltonicity of some examples of standard zero forcing TAR graphs from \cite{XTARiso}.  The \emph{standard zero forcing TAR graph}  of a  base graph  $G$ is denoted by  $\ztar(G)$ and  the upper standard zero forcing number of $G$ is denoted here by $\zbar(G)$
. 


As noted in \cite[Section 3.1]{XTARiso},
     $K_n$ is the only graph $G$ of order $n$ that has no isolated vertices and for which $\Z(G)=n-1$, so  $\ztar(K_n)$ is unique.  Similarly, $P_n$ is the only graph $G$ of order $n$ that has no isolated vertices and for which $\Z(G)=1$, so  $\ztar(P_n)$ is unique.  It is also shown there that $\ztar(K_{1,q})$ is unique for $q\ge 2$ and that $\ztar(C_n)$ is not unique for $n\ge 4$.
 For the convenience of the reader in comparing across parameters, we reproduce the table of uniqueness data for standard zero forcing TAR graphs of small base graphs. 

\begin{table}[!h]
\begin{center}
\begin{tabular}{ |c|c|c|c|c|c|c|c|c| } 
 \hline
\#  vertices in $G$  & 2 & 3 & 4 & 5 & 6 & 7 
 & 8\\
\hline
\# graphs with  unique $\ztar(G)$  & 1 & 2 & 4 & 7 & 34 & 303 & 5318\\
  \hline
\# graphs with no isolated vertices  & 1 & 2 & 7 & 23 & 122 & 888 & 11302\\
 \hline
 ratio (\# unique/\# no isolated) &  1 & 1 & 0.5714 & 0.3043 & 0.2787 & 0.3412  
& 0.4705\\
 \hline
\end{tabular}\vspace{-9pt}
\end{center}
    \caption{\cite[Table 1]{XTARiso} Number of graphs with unique standard zero forcing TAR graph for small orders}
    \label{tab:unique-tar-data=Z}
\end{table}

\begin{ex} \label{ex:HamZsmall}
It is shown in \cite{XTARiso} that $\ztar(K_n)\cong K_{1,n}$ and $\ztar(K_{1,q})\cong K_{1,q}\Box K_2$.  Thus $\ztar(K_n)$ does not have a Hamilton path for $n\ge 3$.  Similarly, $\ztar(K_{1,q})$  has a Hamilton cycle if and only if $q\le 2$ if and only if $\ztar(K_{1,q})$ has a Hamilton path. 
\end{ex}

\begin{ex}\label{ex:bigHam}
    Observe that if $H$ is a graph of even order with a Hamilton cycle $C$, then $H\Box P_r$ also has a Hamilton cycle, by alternating edges of $C$ in the first and last copies of $H$ and moving between these two cycles on copies of $P_r$.  Note that $\ztar(K_{1,2})$ has a Hamilton cycle and is of even order and $\ztar(K_2)\cong P_3$ (Example \ref{ex:HamZsmall}).  So  $\ztar(K_{1,2}\du K_2\dots \du K_2)\cong \ztar(K_{1,2}\Box P_3\dots\Box P_3$) by Proposition \ref{p:disjoint- cart-u}.  Hence there are arbitrarily large order base graphs $K_{1,2}\du K_2\dots \du K_2$ that have standard zero forcing TAR graphs with Hamilton cycles. 
\end{ex}

  The least $k_0$ such that $\ztar_k(G)$ is connected for all $k\ge k_0$ is denoted by $\zzo(G)$  and the least $k$ such that $\ztar_k(G)$ is connected is denoted by $\ulzo(G)$.  It is common to have $\zzo(G)=\zbar(G)+1$ and $\ulzo(G)=\zzo(G)$, as illustrated by $K_n$ and in the next example.

  \begin{ex}\cite{XTARiso}
      Number  the vertices of the path $P_n$ with $n\ge 4$  in path order. A set $S\subseteq V(P_n)$ is a standard zero forcing set if and only if $S$ contains an endpoint or $S$ contains two consecutive vertices of the path. The set $\{2,3\}$ is a standard zero forcing set, but is not adjacent to any  zero forcing set in $\ztar_2(P_n)$. Thus $\zbar(G)=2$ and $\ztar_2(P_n)$ is not connected.  Adding an end vertex makes any set a standard zero forcing set, so $\ztar_k(G)$ is connected for $k\ge 3$.   Thus $\ulzo(P_n)=3=\zzo(P_n)=\zbar(P_n)+1$.
  \end{ex}
  
  Next we  present connectedness examples from \cite{XTARiso} showing $\zzo(G)>\zbar(G)+1$ and $\zzo(G)>\ulzo(G)$ are possible.  The graph $H(r)$ is constructed from two copies of $K_{r+2}$ by adding a matching between $r$ pairs of vertices. 
 
 \begin{prop}
   {\rm \cite[Proposition 3.4]{XTARiso}} For $r\ge 2$, $\Z(H(r))=\zbar(H(r))=r+2$ and $\ulzo(H(r))=\zzo(H(r)) = 2r+2= \zbar(H(r)) + r$.
 \end{prop}

 The next family of graphs provides examples of graphs $G$ such that $\zzo(G)>\ulzo(G)$.  For $r\ge 3$, construct $H_r$  from $H_2$ (shown in Figure \ref{f:z-1stconn}) by adding vertices $u_3,\dots,u_r$ with  $N_{H_r}(u_k)=N_{H_r}(u_1)$ for $k=3,\dots,r$ \cite{XTARiso}.  Note that    $|V(H_r)|=r+6$.  The next result is established by use of the {twins property} (although it is not called that in \cite{XTARiso}).
 \begin{figure}[h!]
\centering 

	\begin{tikzpicture}[scale=.8,every node/.style={draw,shape=circle,outer sep=2pt,inner sep=1pt,minimum
			size=.4cm}]
		
		\node[fill=none]  (0) at (-1.5,2) {};
		\node[fill=none]  (1) at (0,2.5) {};
		\node[fill=none]  (2) at (1.5,2) {};
		\node[fill=none]  (3) at (-1.5,0.5) {};
		\node[fill=none]  (4) at (0,0) {};
		\node[fill=none]  (5) at (1.5,0.5) {};
		\node[fill=none]  (6) at (-0.8,-1.8) {\scriptsize{$u_1$}};
		\node[fill=none]  (7) at (0.8,-1.8) {\scriptsize{$u_2$}};
		
		
		\draw[thick] (0)--(1)--(2)--(5)--(4)--(3)--(0)--(5)--(1)--(4)--(7)--(3)--(6)--(4)--(0);
		\draw[thick] (6)--(4)--(2)--(3)--(1);
	    \draw[thick](7)--(5)--(6);
		\end{tikzpicture}
		
    \caption{A graph $H_2$ satisfying $\zzo(H_2)>\ulzo(H_2)$ }  \label{f:z-1stconn}
\end{figure}

\begin{prop}  {\rm \cite[Proposition 3.9]{XTARiso}} For $r\ge 2$, $\Z(H_r)=r+2$, $\zbar(H_r)=r+4$, $\ulzo(H_r)=r+3$, and $\zzo(H_r)=r+5$.
\end{prop}

\section{TAR reconfiguration for positive semidefinite (PSD)  forcing}\label{s:PSD}
 It was noted in \cite{XTARiso} that PSD forcing is an original $X$-set parameter, so the original $X$-set parameter results in \cite{PDrecon, XTARiso} apply to PSD forcing. In particular, Theorem \ref{t:main} applies to PSD forcing.  It was noted there that there are no irrelevant vertices for PSD forcing \cite[Proposition 2.19]{XTARiso}.  Here we determine PSD forcing TAR graphs (also called $\Zp$-TAR graphs) for certain graph families, discuss uniqueness, Hamiltonicity,   and connectedness of $\Zp$-TAR graphs. Examples showing connectedness bounds need not be equalities are presented.   Denote the \emph{$\Zp$-TAR graph}  of a  base graph  $G$ by  $\zptar(G)$ and denote the upper PSD forcing number of $G$  by $\zpbar(G)$.

\subsection{Examples of PSD TAR graphs, uniqueness, and Hamiltonicity}\label{s:PSDunique}
 
 We  present examples of $\Zp$-TAR graphs that illustrate both uniqueness and nonuniqueness and the existence of graphs with Hamilton cycles, Hamilton paths but not cycles, and not having Hamilton paths.

\begin{ex}
     Since $\Zp(K_n)=n-1$ and $K_n$ is the only graph of order $n$ with this property, $\zptar(K_n)\cong K_{1,n}$ and $K_n$ has a unique  $\Zp$-TAR graph.  Thus $\zptar(K_n)$ does not have a Hamilton path for $n\ge 3$.   
 \end{ex}

 \begin{ex}
     For a tree $T$ of order $n$, any one vertex is a PSD zero forcing set (and the empty set is not a PSD forcing set of any graph).  Thus $\zptar(T)\cong Q_n-v$ for $v\in V(Q_n)$ (cf. Remark \ref{r:hypercube-bip}). Thus $\zptar(T)$ is not unique  for $n\ge 4$.   Furthermore, $\zptar(T)$ has a Hamilton path for every tree.
 \end{ex}

 \begin{ex}
     For a cycle $C_n$, any set of two vertices  is a PSD zero forcing set and $\Zp(C_n)=2$.  Thus $\zptar(C_n)$ is isomorphic to the graph obtained from $Q_n$ by removing the bottom two layers (corresponding to sets of cardinality zero and one).  That is, $\zptar(C_n)\cong Q_n-V(K_{1,n})$.  Since $Q_n-V(K_{1,n})$ has $2^{n-1}-1$ vertices of even cardinality and $2^{n-1}-n$ vertices of odd cardinality, $\zptar(C_n)$ does not have a Hamilton path for $n\ge 3$.
 \end{ex}

\begin{ex}\label{ex:Kab-PSD}
 Let $q\ge 1$. Since $K_{1,q}$ is a tree, $\Zp(K_{1,q})=1=\zpbar(K_{1,q})$ and $\zptar(K_{1,q})\cong Q_{q+1}-v$ for $v\in V(Q_{q+1})$.  Thus $\zptar(K_{1,q})$ is not unique for $q\ge 3$ and $\zptar(K_{1,q})$ does have a Hamilton path.
  
For $K_{2,2}$, every set of two vertices is a PSD forcing set and no one vertex can force.    For $q\ge 3$,  the minimal PSD forcing sets   of $K_{2,q}$ are  $A$, $B$ and $\{a_i,b_j\}$ where $a_i\in A$ and $b_j\in B$.   
 For $p\ge 3$, a set $S\subseteq V(K_{a,b})$ is a PSD forcing set if and only if $A\subseteq S$, $B\subseteq S$, ($|S\cap A|=a-1$ and $|S\cap B|\ge 1$) or ($|S\cap B|=b-1$ and $|S\cap A|\ge 1$).  Thus $\Zp(K_{p,q})=p$ and $\zpbar(K_{p,q})=q$ for $q\ge p\ge 2$. 
\end{ex}

 The next result follows from Theorem \ref{t:KpqHdom} together with Remark \ref{r:sameXsets}, since $K_{2,q}$ has the same minimal sets for PSD forcing and domination (cf. Examples  \ref{ex:Kab-PSD} and \ref{ex:KabDom}). 
\begin{cor}
  For $q\ge 3$, $\zptar(K_{2,q})$ has a Hamilton path if and only if $q$ is odd (and  $\zptar(K_{2,q})$ does not have a Hamilton cycle). 
\end{cor}

  \begin{prop}\label{p:Kab-unique-PSD}
For $q\ge 2$, $\zptar(K_{2,q})$ is unique.
\end{prop}
\bpf
 As usual,  $B$ is the partite set with $q$ vertices and  $A$ is the partite set with $p=2$ vertices. The cases $q=2, 3$ are straightforward to verify computationally; see \cite{sage:recon}. So assume $q\ge 4$. The minimal PSD forcing sets of $K_{2,q}$ are $A$, $B$ and $\{a_i,b_j\}$ where $a_i\in A$ and $b_j\in B$. 
Suppose $G$ has no isolated vertices and  $\zptar(G)\cong \zptar(K_{2,q})$.
By Theorem \ref{t:main}, $\Zp(G)=2$ and we can relabel the vertices of $G$ so that the PSD forcing sets are the same as for $K_{2,q}$. Recall that $\Zp(G)=2$ implies $G$ cannot contain a $K_4$ minor \cite{EGR11}.
We consider two cases, depending on whether or not the vertices in $A$ are adjacent.  

First  assume  $G[A]=2K_1$.
If $N_G(a_1)\ne B$, then $N_G(a_1)$ is a PSD forcing set, contradicting that $B$ is a minimal PSD forcing set.  So $N_G(a_1)=B=N_G(a_2)$.  Let $C$ be a connected component of $G[B]$.  Since $a_1$ and $a_2$ are each adjacent to every vertex of $B$ and $A$ is a PSD forcing set, $|C|=1$.  Thus $G\cong K_{2,q}$.

Now assume $G[A]=K_2$. 
 We show first that for any vertex $b_i\in B$,
$|N_G(b_i)\cap A|=1 \implies B\subseteq N_G[b_i]$. 
Suppose $|N_G(b_i)\cap A|=1$. The one neighbor of $b_i$ in $A$, say  $a_1$, can be  forced blue by $b_i$.  Then $a_1$ forces $a_2$ so $A$ is blue.  Thus $N_G[b_i]\cap B$ is a PSD forcing set and the miminality of $B$ as a PSD forcing set implies  $B\subseteq N_G[b_i]$.
Since $B$ is a PSD forcing set, there exists a vertex $b_k$ such that  $N_G(b_k)\cap A=\{a_1\}$   
 (possibly after relabeling $a_1$ and $a_2$). 
 Thus $N_G[b_k]=B\cup \{a_1\}$.  Since a leaf and its one neighbor cannot both be in a minimal PSD forcing set, and  for every pair of vertices there exists a minimal PSD forcing set that contains those vertices, 
 $\delta(G)\ge 2$.  
 Thus there exists a vertex $b_\ell\in B$ such that $b_\ell\in N_G(a_2)$  (necessarily $\ell\ne k$ since $N_G(b_k)\cap A=\{a_1\}$).  Let $b_i\in B$ and $i\ne k,\ell$. If  $A\subseteq N_G(b_i)$, then we see that $G[\{a_1,a_2,b_k,b_\ell,b_i\}]$ would have a $K_4$ minor   by contracting  the edge $b_kb_{\ell}$. 
 Since this is not allowed, $A\not\subseteq N_G(b_i)$.  If $A\cap N_G(b_i)=\emptyset$, then $B\setminus  \{b_i\} $  would be a  PSD forcing set.  So  $A\cap N_G(b_i)\ne \emptyset$.  This implies that $|N_G(b_i)\cap A|=1$,  which implies $B\subseteq  N_G[b_i]$. 
 This is a contradiction because then contracting edge $a_1a_2$ gives a $K_4$ minor.  Thus $G[A]\not\cong K_2$.
\epf

Table \ref{t:unique-tar-data-PSD} shows the number of (nonisomorphic) graphs without isolated vertices of order at most eight  that have  unique $\Zp$-TAR graphs (this data was computed in \cite{sage:recon}).\begin{table}[!h]
\begin{center}
\begin{tabular}{ |c|c|c|c|c|c|c|c|c| } 
 \hline
\#  vertices in $G$  & 2 & 3 & 4 & 5 & 6 & 7 
 & 8\\
\hline
\# graphs with  unique $\zptar(G)$  & 1 & 2 & 3 &10 & 48 & 398 & 6798 
\\
  \hline
\# graphs with no isolated vertices  & 1 & 2 & 7 & 23 & 122 & 888 & 11302\\
 \hline
 ratio (\# unique/\# no isolated) &  1 & 1 & 0.4286 & 0.4348 &  0.3934 &  0.4482 & 0.6015 \\
 \hline\end{tabular}\vspace{-9pt}
\end{center}
    \caption{Number of graphs having a unique PSD TAR graph for graphs of  small order}
    \label{t:unique-tar-data-PSD}
\end{table}

 We do not have an example of a base graph $G$ such that $\zptar(G)$ has a Hamilton cycle.
\subsection{Connectedness properties of the PSD  TAR graph}

 The least $k_0$ such that $\zptar_k(G)$ is connected for all $k\ge k_0$ is denoted by $\zpo(G)$  and the least $k$ such that $\zptar_k(G)$ is connected is denoted by $\ulzpo(G)$.
As usual, there are many example of graphs having $\zpbar(G)+1= \zpo(G)$  and $\ulzpo(G)= \zpo(G)$, such as complete graphs,  trees, and cycles. 
Here we provide examples of graphs $G$ with $\zpbar(G)+1< \zpo(G)$  and $\ulzpo(G)< \zpo(G)$

 \begin{prop}\label{p:upper-con-PSD}
For $q\ge p\ge 2$, $\Zp(K_{p,q})=p$ and  $\zpbar(K_{p,q})=q$.   If $p\ge 3$ then $\zpo(K_{p,q}) = p+q-2$. Thus for $q=p\ge 4$,
$\zpo(K_{p,p}) = 2p-2= \zpbar(K_{p,p}) + p-2> \zpbar(K_{p,p}) +1$.  For $q\ge 4$, $\ulzpo(K_{2,q})=3$ and $\zpo(K_{2,q}) = q+1$, so $\zpo(K_{2,q}) = q+1=\ulzpo(K_{2,q})+q-2>\ulzpo(K_{2,q})$.
\end{prop}
\bpf For $p\ge 3$, a set $S\subseteq V(K_{p,q})$ is a PSD forcing set if and only if $A\subseteq S$, $B\subseteq S$, ($|S\cap A|=p-1$ and $|S\cap B|\ge 1$) or ($|S\cap B|=q-1$ and $|S\cap A|\ge 1$).  Thus $\Zp(K_{p,p})=p=\zpbar(K_{p,p})$ and $\zptar(K_{p,p})$ is disconnected until a set $S$ can contain at least $p-1$ vertices from each part, i.e. $|S|\ge 2p-2$.  Furthermore, ${\zptar}_k(K_{p,p})$ is connected for $k\ge 2p-2$.   

For $q\ge 4$, a set $S\subseteq V(K_{2,q})$ is a PSD forcing set if and only if $A\subseteq S$, $B\subseteq S$, or ($|S\cap A|=1$ and $|S\cap B|\ge 1$). Thus the minimal PSD forcing sets are $A, B$, and $S$ such that $|S\cap A|=1$ and $|S\cap B|= 1$.  Thus ${\zptar}_2(K_{2,q})$ is disconnected, ${\zptar}_3(K_{2,q})$ is connected, ${\zptar}_q(K_{2,q})$ is disconnected, and ${\zptar}_k(K_{2,q})$ is connected for $k\ge q+1$.
\epf

\section{TAR reconfiguration for skew forcing}\label{s:skew}

 In this section  we  define the skew forcing TAR graph, apply results of Section \ref{sec:proof-requirements}, and obtain additional results for skew TAR graphs, including results about uniqueness, existence or nonexistence of Hamilton paths and cycles, and connectedness. 
   Denote the \emph{skew TAR graph}  of a  base graph  $G$ by  $\zstar(G)$ and deonte  the upper skew forcing number of $G$  by $\zsbar(G)$. 

Since skew forcing is a robust $X$-set parameter  (but not an original $X$-set parameter),  all the results that are true for robust $X$-set parameters and $\xrtar(G)$ apply to skew forcing and the skew TAR graph.  Here we restate only the main result for isomorphic skew TAR graphs.

\begin{thm}\label{t:main-skew}  Suppose base graphs $G$ and $G'$ have no isolated vertices and $\zstar(G)\cong\zstar(G')$.  Then $G$ and $G'$ have the same order and  there is a relabeling of the vertices of $G'$ such that $G$ and $G'$ have exactly the same skew forcing sets.  
\end{thm}

\subsection{Examples of skew TAR graphs, uniqueness, and Hamiltonicity}\label{ss:skewex}

  In this section, 
 we exhibit skew TAR graphs of some families of base graphs, including complete graphs, complete multipartite graphs, odd cycles, and every   graph $G$ with $\Zs(G)=0$ (which includes  even  paths).    We  show that the skew TAR graph of a complete multipartite graph is unique and examine Hamiltonicity.  
 
 \begin{ex}
  It is well known that $\Zs(K_n)=n-2$ and any set of $n-2$ vertices is a skew forcing set.  Thus $\zstar(K_n)$ is isomorphic to the top three `levels' of $Q_n$, i.e., all the subsets of $n-2$, $n-1$, and $n$ vertices.    It is shown in the next theorem that $\zstar(K_n)$ is unique. For $n\ge 4$, $\zstar(K_n)$ does not have a Hamilton path because there are $n$ sets of cardinality $n-1$ and $1+\frac{n(n-1)}2$ sets of cardinality $n$ or $n-2$ and $1+\frac{n(n-1)}2>n+2$.
 \end{ex}

For $t\ge 2$ and $n_i\ge 1, i=1,\dots,t$, a \emph{complete multiparite graph}  is a graph whose vertex set can be partitioned into sets of orders $n_1,\dots,n_t$ such that the edge set is every edge between two distinct partite sets; this is denoted by $K_{n_1,\dots,n_t}$.For a graph $G$ of order $n$, it is known that $\Zs(G)=n-2$ if and only if $G$ is a complete multipartite graph \cite[Theorem 9.76]{HLS22book}.  Note that the next results includes $K_n\cong K_{1,1,\dots,1}$ as well as complete bipartite graphs.
 
\begin{thm}\label{t:skew-comp-multi}
 The  complete multiparite graph $K_{n_1,\dots,n_t}$ has as its minimal skew forcing sets the sets of $n-2$ vertices where the two vertices omitted are in  different partite sets.  These are the vertices of the skew TAR graph that have cardinality $n-2$ (and degree two).  Every set of $n-1$ vertices of $K_{n_1,\dots,n_t}$ is a vertex of the skew TAR graph and has degree $n-n_i+1$ where the vertex omitted is in the $i$th partite set.
The skew TAR graph of a complete multiparite graph is unique.
\end{thm}
\bpf Suppose $\zstar(G)\cong\zstar(K_{n_1,\dots,n_t})$ with $t\ge 2$ and $n_i\ge 1$, and let $n=n_1+\dots+n_t$.  Then $|V(G)|=n$ and $\Zs(G)=n-2$.
Relabel the vertices of $G$ so that the skew forcing sets of $G$ are the same as the skew forcing sets of $K_{n_1,\dots,n_t}$ (by Theorem \ref{t:main-skew}); necessarily $V(G)=V(K_{n_1,\dots,n_t})$. Denote the partite sets of $K_{n_1,\dots,n_t}$ by $Y_i$ where $|Y_i|=n_i$. A set $S$ of $n-2$ vertices of $K_{n_1,\dots,n_t}$ is a skew forcing set if and only if the two vertices not in $S$ are not in the same partite set.  Consider a set $T=V(G)\setminus\{y_i\}$ of $n-1$ vertices  
where  $y_i\in Y_i$.  
Removing a vertex $v$ of $T$ results in a skew forcing set $K_{n_1,\dots,n_t}$ (and thus of $G$)  if and only if $v\not\in Y_i$.  Of course, $T\cup\{y_i\}$ is also a skew forcing set.  Thus  $\deg_{\zstar(G)}(T)=n-n_i+1$.  We can recover $n_1,\dots,n_t$ from the degrees of the  vertices of $\zstar(G)$ that contain $n-1$ elements, because omitting a vertex from $Y_i$ contributes $n_i$ vertices  of degree $n-n_i+1$.
\epf

At the other end of the range of values of $\Zs$, for every even order $n\ge 4$, there are examples known of nonisomporphic graphs $G$  and $H$ of order $n$ that have $\Zs(G) = \Zs(H)=0$. 
For instance,  every path of even order has $\Zs(P_{2k})=0$ and $\Zs(H\circ K_1)=0$ for any graph $H$.  The $s$th-\emph{half-graph}, denoted by Half$_s$, is the graph constructed from a copy of $K_s$ with vertices $\{x_1,\dots,x_s\}$ and a copy of $\overline{K_s}$ with vertices $\{y_1,\dots,y_s\}$ by adding exactly the edges $x_iy_j$ such that $i+j\leq s+1$; $\Zs($Half$_s)=0$. 

\begin{obs}
    The following are equivalent:
    \ben[$(1)$]
    \item $\Zs(G)=0$.
    \item $\emptyset$ is a skew forcing set.
    \item $\zstar(G)$ is a hypercube.
    \een
     As noted in Remark $\ref{r:QnHam}$, a hypercube has a Hamilton cycle.
\end{obs}

   The Leaf-Stripping Algorithm (reproduced below and implemented in \cite{sage:recon}) can be used to determine whether $\Zs(G)=0$, because $\Zs(G)=0$ if and only if the algorithm returns the empty set  \cite[Algorithm 9.79, Theorem 9.80]{HLS22book}.

\begin{alg}[Leaf-Stripping Algorithm]\label{alg-leaf-stripping-S}{\rm \cite{HPRY, King15}}\\
{\bf Input:} Graph $G$.\\
{\bf Output:} Graph $\hat G$ with $\delta(\hat G) \neq 1$, or $\hat{G} = \emptyset$.

\ms

\noi $\hat G := G$.\\
{\bf While} $\hat G$ has a leaf $u$ with neighbor $v$:
	
$\hat G := \hat G[V(\hat G) \setminus \{u,v\}]$.\\
{\bf Return} $\hat G$.
\end{alg}

 Next we present examples of graphs in which any one vertex is a skew forcing set, and again nonisomorphic base graphs have isomorphic skew TAR graphs.  It is well known that $\delta(G)-1\le \Zs(G)$ \cite[Remark 9.63]{HLS22book}.  For $r\ge 2, s\ge 3$, the \emph{flower} with $r$ petals of size $s$ (or the $(r,s)$-\emph{generalized friendship graph}), denoted by $F(r,s)$, is a union of $r$ copies of $C_s$ that share one common vertex.  Note that  $F(r,s)$ has $(s-1)r+1$ vertices and $sr$ edges. 
The \emph{center} vertex of a flower is the vertex that has degree at least four.  

\begin{prop}
    The following graphs $G$ have the property that any one vertex is a skew forcing set but the empty set is not a skew forcing set, so $\zstar(G)\cong Q_n-v$ for any $v\in V(Q_n)$:
    \ben[$(1)$]
    \item\label{odd-cyc}
    An odd cycle.
    \item\label{flower} A flower $F(r,2k+1)$ for $r\ge 2, k\ge 1$.
    \een
 As noted in Remark $\ref{r:QnHam}$, a  graph $Q_n-v$ for $v\in V(Q_n)$ has a Hamilton path (but not a Hamilton cycle).
\end{prop} 
 \bpf  In all cases, $\delta(G)=2$, so $\Zs(G)\ge 1$.

\eqref{odd-cyc}:    In an odd  cycle, coloring any one vertex blue leaves an even path, which can then force itself (it is well known that $\Zs(P_{2k})=0$).

\eqref{flower}:  If the center vertex of  is blue, the white vertices in each petal forms an even path, which can force itself.  So suppose a degree-two vertex $u$ is blue.  Since each petal has an odd number of vertices, there is a path with an even number of vertices from a white  neighbor $w$ of $u$  to the center vertex $c$ (including $w$ and $c$).  Using white vertex forcing along this path, $c$ can be forced.
\epf


 Returning to the question of uniqueness of skew TAR graphs, Table \ref{t:unique-tar-data} shows the number of (nonismorphic) graphs without isolated vertices of order at most eight  that have  unique $\Zs$-TAR graphs (this data was computed in \cite{sage:recon}).
\begin{table}[!h]
\begin{center}
\begin{tabular}{ |c|c|c|c|c|c|c|c|c| } 
 \hline
\#  vertices in $G$  & 2 & 3 & 4 & 5 & 6 & 7 
 & 8\\
\hline
\# graphs with  unique $\zstar(G)$  & 1 & 2 & 4 & 7 & 27 & 179 & 3026 
\\
  \hline
\# graphs with no isolated vertices  & 1 & 2 & 7 & 23 & 122 & 888 & 11302\\
 \hline
 ratio (\# unique/\# no isolated) &  1 & 1 & 0.5714 & 0.3043 &  0.2213 &  0.2016 & 0.2677 \\
 \hline
\end{tabular}\vspace{-9pt}
\end{center}
    \caption{Number of graphs having a unique skew TAR graph for graphs of  small order}
    \label{t:unique-tar-data}
\end{table}

For $n=2,3$ and $4$, the number of nonisomorphic graphs $G$ with unique $\zstar(G)$  is the same as the number of complete multipartite graphs of order $n$ (which is one less than the number of partitions of $n$).  But beginning with $n=5$ there exists at least one additional base graph $G$ of each order that has a unique skew TAR graph.  For $n=5$ there is only one, the Full House graph \cite{sage:recon}; see Example \ref{ex:FH}.

\begin{figure}[!h]
\centering
\scalebox{1}{
\begin{tikzpicture}[scale=1.3,every node/.style={draw,shape=circle,outer sep=2pt,inner sep=1pt,minimum size=.2cm}]		
\node[fill=none]  (1) at (1,-1) {$2$};
\node[fill=none]  (2) at (1,0) {$3$};
\node[fill=none]  (3) at (0.5,0.866) {$0$};
\node[fill=none]  (4) at (0,0) {$4$};
\node[fill=none]  (5) at (0,-1) {$1$};
		
\draw[thick] (1)--(2)--(3)--(4)--(2)--(5)--(1)--(4)--(5);
\end{tikzpicture}}
\caption{The Full House graph}
\label{f:FullHouse}
\end{figure}

\begin{ex}
    \label{ex:FH}
    The Full House graph,  shown in Figure \ref{f:FullHouse}, has minimal skew forcing sets $\{3\},\{4\},\{0,1,2\}$. 
\end{ex}

\subsection{Connectedness properties of the skew forcing TAR graph}\label{ss:skew-conn}

In this section we discuss connectedness properties of the skew forcing TAR graph. 
  The least $k_0$ such that $\zstar_k(G)$ is connected for all $k\ge k_0$ is denoted by $\zso(G)$  and the least $k$ such that $\zstar_k(G)$ is connected is denoted by $\ulzso(G)$. We begin with a simple observation.

\begin{obs}
Suppose the empty set is a skew forcing set of a graph $G$ of order $n$, which implies $\zstar(G)\cong Q_n$ and thus $\zstar_k(G)$ is connected for every $k\ge 0$. Therefore $\zso(G)=\ulzso(G)=\Zs(G)=0$ (cf. Proposition  \ref{basic1}).    
\end{obs}

 Many well known graphs have  $\zso(G)=\zsbar(G)+1$ and $\ulzso(G)=\zso(G)$.  Any graph for which $\Zs(G)=0$ satisfies $\ulzso(G)=\zso(G)$.
It follows from results  in Section \ref{ss:skewex} that  the graphs $K_{n_1,\dots,n_t}, C_n$, and $F_r$ all have both these properties.  We show that the family of graphs $H(r)$ presented in \cite{XTARiso} to show that $\zzo(G)$ can exceed the lower bound $\zbar(G)+1$  by an arbitrary amount also shows that $\zso(G)$ can exceed the lower bound $\zsbar(G)+1$  by an arbitrary amount, although the values are different. 
 We also apply the  twin vertex method (Proposition \ref{p:xso-lowxso})  to  construct a family of graphs that have $\ulzo(G)$  strictly less than $\zzo(G)$. 
  Recall the graph $H(r)$ (see Section \ref{ss:ZF}) is constructed from two copies of $K_{r+2}$ by adding a matching between $r$ pairs of vertices. 

\begin{prop}\label{p:upper-con-skew}
For $r\ge 2$, $\Zs(H(r))=r$, $\zsbar(H(r))=r+1$, and 
$\zso(H(r)) = 2r= \zsbar(H(r)) + r-1$.
\end{prop}

 \bpf Observe that $\delta(H(r))=r+1$.  Let the vertices of degree $r+1$ be denoted by $x,y,x',y'$ where $x$ and $y$ are adjacent, as are $x'$ and $y'$. Let $U=\{u_1,\dots,u_r\}$ be the remaining neighbors of $x$, $U'=\{u'_1,\dots,u'_r\}$ be the remaining neighbors of $x'$,  $V=U\cup\{x,y\}$, and $V'=U'\cup\{x',y'\}$. Then $U$ is a skew forcing set, so 
$r=\delta(H(r))-1\le\Zs(H(r))\le r.$ In fact, any set of $r$ vertices in $V$ that contains at most one of $x$ or $y$ is a skew forcing set (and similarly for $V'$). A set $S$ of $r$ vertices in $V$ that contains both $x$ and $y$ is not a skew forcing set, but adding a neighbor of $v\in (V\setminus S)$ produces a minimal skew forcing set of $r+1$ vertices. No set of $r+2$ or more vertices of $H(r)$ is minimal, so $\zsbar(H(r))=r+1$.  We partition the minimal skew forcing sets into those that have at least $r$ vertices in $V$ and  those that have at least $r$ vertices in $V'$.  Then by Proposition \ref{p:part-min-sets}, $\zso(H(r))\ge (r+1)+(r+1)-2=2r$. Since $|S\cup T|\le 2r$ for any minimal skew forcing sets $S$ and $T$, 
 $\zso(H(r))=2r$.
\epf

 A family of graphs $H_r$ having $\ulzo(H_r)<\zzo(H_r)$ was constructed in \cite{XTARiso} by using twins.  This family also has the property that $\ulzso(H_r)<\zso(H_r)$, but for skew forcing there are  graphs with this property of  smaller order, including  the family obtained by creating twins of the degree two vertex in the Full House graph.  

 The next result is the skew analog of \cite[Proposition 3.6]{XTARiso}, but follows from known results: Let $G$ be a graph 
with a set of $r$ independent twins.  It is immediate that 
 any skew forcing set of $G$ must contain at least $r-1$ of the vertices of $T$. The remaining two conditions for the twins property are established in the proof of Proposition 9.87 in \cite{HLS22book}. 

\begin{prop}\label{p:skew-twin} 
Skew forcing has the twins property.
\end{prop}

 Now we  construct the family $FH(r)$: Define $FH(1)=FH$ and construct $FH(r+1)$ from $FH(r)$ by adding another independent twin of a degree-two vertex; $FH(2)$ is shown in  Figure \ref{f:FH2}

\begin{figure}[!h]
\centering
\scalebox{1}{
\begin{tikzpicture}[scale=1.3,every node/.style={draw,shape=circle,outer sep=2pt,inner sep=1pt,minimum size=.2cm}]		
\node[fill=none]  (1) at (1,-1) {};
\node[fill=none]  (2) at (1,0) {};
\node[fill=none]  (3) at (1,1) {};
\node[fill=none]  (4) at (0,0) {};
\node[fill=none]  (5) at (0,-1) {};
\node[fill=none]  (6) at (0,1) {};
		
\draw[thick] (1)--(2)--(3)--(4)--(2)--(5)--(1)--(4)--(5);
\draw[thick] (2)--(6)--(4);
\end{tikzpicture}}
\caption{A graph $FH(2)$ satisfying $\zso(FH(2))>\ulzso(FH(2))$}
\label{f:FH2}
\end{figure}

\begin{prop}\label{p:zso-lowzso}  For $r\ge 1$, $\Zs(FH(r))=r$, $\zsbar(FH(r))=r+2$, $\ulzso(FH(r))=r+1$, and $\zso(FH(r))=r+3$.
\end{prop} \bpf The results for $r=1$ are in Example \ref{ex:FH} and are computed in \cite{sage:recon} for $r=2$.  
The first three equalities then follow from  Proposition \ref{p:xso-lowxso}. Since  $\zsbar(FH(r))+1\le \zso(FH(r))$,  ${\zstar}_{r+3}(FH(r))$ is connected, and $|V(FH(r))|=r+4$, we have $\zso(FH(r))=r+3$.
\epf

\subsection{Skew irrelevant vertices}\label{s:zs-irrel}

Recall that a vertex $v$ of a graph $G$  is \emph{skew irrelevant} if it is not in any minimal skew forcing set of $G$.   Irrelevant vertices play a key role in the automorphism group of an $X$-TAR graph. In this section we examine skew irrelvant vertices.

\begin{obs}\label{o:null-force-irrelevant}
   If a vertex $v$ of $G$ can be skew forced by the empty set, then $v$ is skew-irrelevant.
\end{obs}

We  see in the Example \ref{o:not-null-force-irrelevant} 
that the converse of Observation \ref{o:null-force-irrelevant} is not true.   However, it is true for trees. 

\begin{prop}\label{p:tree-skew}
    Let $T$ be a forest.  A vertex $v$ of $T$ is skew irrelevant if and only if it can be skew forced by the empty set.
\end{prop}
\bpf Observe that applying the Leaf-Stripping Algorithm \ref{alg-leaf-stripping-S} to a forest always returns a set (possibly empty) of isolated vertices.  During any application of the algorithm, $\hat T$ is a forest after each iteration. Let $v\in V(T)$. Since a forest that has an edge must have at least two leaves, we can always choose to  not have $v$ perform a force even if $v$ is a leaf of $\hat T$.  So at the end of the algorithm either $v$ was forced by the empty set or $v\in V(\hat T)$ and $V(\hat T)$ is a minimum skew forcing set.
\epf

\begin{ex}\label{o:not-null-force-irrelevant}
    Let $F(r)$ denote the graph constructed from $r$ copies of $K_3$ and one $K_{1,r}$, where each leaf of the $K_{i,r}$ is a vertex of one $K_3$. The graph $F(3)$ is shown in Figure \ref{fig: xxxxxxxxx}. Then every minimal skew forcing set of $F(r)$ has $r-1$ vertices taken from distinct copies of $K_3$.  Notice that the center vertex of the $K_{1,r}$ is irrelevant but no vertex can be forced by the empty set. The irrelevant vertex means that the skew TAR graph is a Cartesian product: Let $W$ be the set of  skew forcing sets of $F(r)$ that do not contain the irrelevant vertex.  Then
    $\zstar(F(r))\cong \lp \zstar(F(r))[W]\rp\Box P_2$.  
\end{ex}

\begin{figure}[!h]
\centering
\scalebox{.7}{
\begin{tikzpicture}[scale=1.3,every node/.style={draw,shape=circle,outer sep=2pt,inner sep=1pt,minimum size=.285cm}]		
\node[fill=none]  (1) at (0,0) {};
\node[fill=none]  (2) at (0,1) {};
\node[fill=none]  (3) at (-0.5,1.866) {};
\node[fill=none]  (4) at (0.5,1.866) {};
\node[fill=none]  (5) at (0.866,-0.5) {};
\node[fill=none]  (6) at (1.866,-0.5) {};
\node[fill=none]  (7) at (1.366,-1.366) {};
\node[fill=none]  (8) at (-0.866,-0.5) {};
\node[fill=none]  (9) at (-1.366,-1.366) {};
\node[fill=none]  (10) at (-1.866,-0.5) {};
		
\draw[thick] (5)--(7)--(6)--(5)--(1)--(2)--(3)--(4)--(2);
\draw[thick] (1)--(8)--(9)--(10)--(8);
\end{tikzpicture}}
\caption{The graph $F(3)$}
\label{fig: xxxxxxxxx}
\end{figure}

\section{TAR reconfiguration for vertex covers}\label{s:vtxcover}
    In this section  we  define the vertex cover TAR graph, apply results of Section \ref{sec:proof-requirements}, and obtain additional results for such TAR graphs, including uniqueness, connectedness, and irrelevant vertex results.  Recall that the vertex cover  number of $G$ is denoted by $\tau(G)$, so the upper vertex cover  number of $G$ is denoted  by $\vcbar(G)$.  
   Denote \emph{vertex cover TAR graph}  or \emph{VC TAR graph} of a  base graph  $G$ by  $\vctar(G)$. 
 As with the other parameters, our focus here is on isomorphisms of TAR graphs, connectedness, existence of Hamilton cycles or paths, etc.  There is a close relationship between vertex cover reconfiguration  and independent set reconfiguration (see Section \ref{ss:independent}), and there is
is extensive prior work on reconfiguration of both independent sets and vertex covers, including TAR reconfiguration; often the work on independent sets was done first.
Most of this work is focused on complexity, approximation, and algorithms (see, for example, \cite{vtxcov-recon, IDHPSUU11, INZ15, KMM12}).  However, we point out that some of these papers provide useful tools for the type of results we emphasize.  For instance, results from \cite{INZ15} and \cite{KMM12} are used to establish equality in  connectedness bounds for certain families of graphs in Propositions \ref{vctau0ehfg} and \ref{p:ind-evenholefree}.

Since the vertex cover number is a robust $X$-set parameter  (but not an original $X$-set parameter because $\vc(K_1)=0$),  all the results that are true for robust $X$-set parameters and $\xrtar(G)$ apply to vertex covers and the VC TAR graph.  
Here we restate only the main result for isomorphic VC TAR graphs.

\begin{thm}\label{t:main-vc}  Let $G$ and $G'$ be graphs. If $\vctar(G)\cong\vctar(G')$, then  $G$ and $G'$ have the same order and  there is a relabeling of the vertices of $G'$ such that $G$ and $G'$ have exactly the same vertex covers.  
\end{thm}


The robust $X$-set parameter vertex cover number is the base graph complement parameter to the robust $Y$-set parameter independence number (in the sense that $S$ is independent if and only if $V(G)\setminus S$ is a vertex cover);  more information about the base graph complement parameter relationship can be found in Section \ref{ss:comp-param} and more information about independence number  can be found in Section \ref{ss:independent}.

\subsection{Examples of VC TAR graphs, uniqueness, and  Hamiltonicity}\label{ss:vcex}

 In this section, 
 we exhibit vertex cover TAR graphs of some families of base graphs, including  complete graphs, complete bipartite graphs,  and empty graphs, including examples with no Hamilton path, Hamilton path but no Hamilton cycle, and Hamilton cycle.    We also show every vertex cover TAR graph is unique.
 
\begin{ex}
    Every  vertex cover for $K_n$ needs $n-1$ vertices.  Thus $\vc(K_{n})=n-1=\vcbar(K_{n})$, and  $\vctar(K_n)\cong K_{1,n}$.  Furthermore, $K_n$ does not have a Hamilton path for $n\ge 3$.
\end{ex}

\begin{ex}\label{ex:vc-Kpq}
The   two partite sets $A$ and $B$  of $K_{p,q}$ are minimal vertex covers and these are the only minimal vertex covers. Thus $\vctar(K_{p,q})$ is the vertex sum of hypercubes  $Q_p$ and $Q_q$,  $V(K_{p,q})=A\du B$ is a cut-vertex of $\vctar(K_{p,q})$, $\vc(K_{p,q})=p$,  $\vcbar(K_{p,q})=q$, and  $\vco(K_{p,q})=p+q$.  Note that $\vctar_q(K_{p,q})$ is disconnected, but if $p<q$, then $\ulvco(K_{p,q})=p$. Furthermore, 
$\vctar(K_{p,q})$ has a Hamilton path but not a Hamilton cycle.
\end{ex}

\begin{ex}
A path on $n$ vertices has $\vc(P_n)=\lf\frac{n}{2}\rf$. Label the vertices of $P_n$ in  path order.  Using the hypercube representation, a subset $S$ of $V(P_n)$ is represented by  the sequence $(s_1, \dots, s_n)$ where $s_i=1$ if $i\in S$ and $s_i=0$ if $i\not\in S$. A set $S$ is a vertex cover  for $P_n$ if and only if the sequence for $S$ has no consecutive zeros.  As is standard in the hypercube representation, two vertices of $\vctar(P_n)$ are adjacent if and only if they differ in exactly one digit.
In particular, $P_4$ has vertex cover sets $1111$, $1110$, $1101$, $1011$, $0111$, $1010$, $0101$ and $0110$.  A Hamilton cycle for $\vctar(P_4)$ is  described by the  sequence 
\[(1111,1011,1010,1110,0110,0111,0101,1101,1111) \]
\end{ex}

\begin{ex}
Since the vertex cover number is a robust $X$-set parameter  and $\vc(K_1)=0$,    $\vctar(\ol{K_n})\cong K_2\Box\dots\Box K_2\cong Q_n$ for $n\ge 1$ by Proposition \ref{p:disjoint- cart-u}.  Thus  $\vctar(\ol{K_n})$ has a Hamilton cycle.
\end{ex}

\begin{prop}\label{p:VCunique} 
Every vertex cover TAR graph is unique.     \end{prop}
\bpf  Let $G$ be a graph.  For $u,w\in V(G)$, the set $V(G)\setminus\{u,w\}$ is a vertex cover if and only if $uw\not\in E(G)$.  Thus the 
vertex cover TAR graph $\vctar(G)$ determines the base graph.
Suppose $\vctar(G)\cong\vctar(G')$.  By Theorem \ref{t:main-vc}, we may relabel the vertices of  $G'$ to obtain $G''$ such that $\vctar(G)$ and 
 $\vctar(G'')$   have the same vertex covers.  Thus  $G$ and $G''$  are the same graph, and $G= G''\cong G'$. So the VC TAR graph is unique.
\epf\vspace{-12pt}

\subsection{Connectedness properties of the VC TAR graph}\label{ss:vc-conn}

   The least $k_0$ such that $\vctar_k(G)$ is connected for all $k\ge k_0$ is denoted by $\vco(G)$  and the least $k$ such that $\vctar_k(G)$ is connected is denoted by $\ulvco(G)$.  Although complexity and algorithms are main interests of Ito, Nooka, and Zhou in \cite{INZ15},  they also establish equality in a lower bound for $\vco(G)$  for certain families of graphs $G$. 

\begin{rem}
Let $G$ be a graph and let $C_0$ and $C_t$ be two vertex covers of $G$.  As defined \cite{INZ15}, the \emph{minmax vertex cover reconfigutation problem} is to determine the least index $k$ such that there is a path between $C_0$ and $C_t$ in $\vctar_k(G)$. The minmax vertex cover reconfigutation problem is one of the main problems studied in \cite{INZ15}.  When  maximized over all pairs $C_0, C_t$, 
this is equivalent to determining $\vco(G)$.
\end{rem}

   There are many examples of graphs that having $\vco(G)>\vcbar(G)+1$, including complete bipartite graphs; see, e.g.,  Example \ref{ex:vc-Kpq} with $p\ge 2$ or $p=1,q\ge 3$.  As shown in that example,  $ \vco(K_{p,q})=\vcbar(K_{p,q})+\vc(K_{p,q})$, so $\ulvco(K_{p,q})=p<p+q=\vco(K_{p,q})$. 
The next result provides examples where $\vco(G)=\vcbar(G)+1$.
An {\em even--hole--free graph} is a graph which contains no induced even cycles. The class of even--hole--free graphs includes trees, chordal graphs and interval graphs  (see \cite[Section 5.5]{Diestel17} for definitions of chordal graphs and interval graphs). 

\begin{prop} \label{vctau0ehfg}
Let $G$ be an even--hole--free graph. Then $\tau_0(G)=\overline{\tau}+1$.
\end{prop}

\bpf
     Lemma 2 in \cite{INZ15} uses an analogous result for independent sets in \cite{KMM12} (see Theorem  \ref{t:ind-evenholefree}) to show that if $S_1$ and $S_2$ are vertex covers of a graph $G$, then there exists a path  in $\vctar_k(G)$ for $k= \max(|S_1|,|S_2|)+1$.
  Let $\ell=\vcbar(G)$ and let $S_1,S_2 \in V(\vctar_{\ell+1}(G))$. Then there exists minimal vertex covers $M_1$ and $M_2$ such that $M_i \subseteq S_i, i=1,2$. Since $|M_i| \leq \ell,i=1,2$, there exists a path in $\vctar_{\ell+1}(G)$ from $M_1$ to $M_2$. Then starting at $S_1$ and removing the vertices in $S_1\setminus M_1$ one at a time, traversing the path from $M_1$ to $M_2$, and then adding the vertices in $S_2\setminus M_2$ one at a time gives a path from $S_1$ to $S_2$. Since $\vctar_{\ell+1}(G)$ is connected, $\vco(G)=\vcbar(G)+1$ by Propositions \ref{basic2}\eqref{conn3} and \ref{basic1}\eqref{h3}.
    \epf
\subsection{VC-irrelevant vertices}\label{s:vc-irrel}

 A vertex $v$ is \emph{VC-irrelevant} if it is not in any minimal vertex cover.  

\begin{prop} \label{p:vc-irrel}
  Let $G$ be a graph 
  and $v\in V(G)$.  Then $V$ is $VC$-irrelevant if and only if $v$ is an isolated vertex.
\end{prop}
\bpf 
An isolated vertex is not in any minimal vertex cover.  Now suppose $u$ is incident with an edge.  A minimal vertex cover containing $u$ can be  constructed by starting with $S=\{u\}$ and repeatedly adding one vertex at a time, choosing a vertex  $w$ that is not an endpoint of any edge that has already been covered, until all edges are covered. 
\epf


\section{Connected domination}\label{s:conD}

Connected domination is a parameter that has been studied in the literature \cite{conn-Dom, conn-Dom-orig} and  is an example of a super $X$-set parameter that is not robust, and in fact, not a connected $X$-set parameter, because it does not satisfy the $(n-1)$-set axiom. Furthermore, graphs of different orders can have the same connected domination TAR graphs.  For a connected graph $G$, a \emph{connected dominating set} is a (standard) dominating set $S$ of a graph $G$ such that $G[S]$ is connected, and the \emph{connected domination number}, $\Dc(G)$, is the minimum cardinality of a connected dominating set.  
As noted in \cite{conn-Dom},  every superset of a connected dominating set is a connected dominating set, because if $S$ is a connected dominating set then  every vertex not in $S$  is a neighbor of a vertex in $S$. Thus connected domination is a super $X$-set parameter, and the connected domination TAR reconfiguration graph of a base graph $G$ is denoted by $\dctar(G)$.  However, as seen in the next example, connected domination does not satisfy the $(n-1)$-set axiom, so it is not  a connected $X$-set parameter (and is not robust). The next example shows it is possible to have $\dctar(G)\cong Q_r$ for some $r<|V(G)|$. Example \ref{ex:double-broom} shows that base graphs of different orders can have the same connected domination TAR graph.

\begin{ex}\label{ex:conndom1}  Consider the star $K_{1,n-1}$.  For $n\ge 3$,  every connected dominating set contains the center vertex. 
Thus the set of $n-1$ leaf vertices is not a connected dominating set.  Furthermore,  $\dctar(K_{1,n-1})\cong Q_{n-1}$.  
\end{ex}

For $r,s,t\ge 2$, a \emph{double-broom} $DB_r(s, t)$ is a    tree obtained by attaching $s$ leaves to one end vertex of a
path $P_r$ and $t$ leaves to the other end vertex of the path.

\begin{ex}\label{ex:double-broom}    Observe that  every connected dominating set of $DB(r,s,t)$ must contain the $r$ vertices of the path. 
Furthermore,  $\dctar(DB_r(s,t))\cong Q_{s+t}$.  Thus base graphs of different orders can have isomorphic connected domination TAR graphs.
\end{ex}

\section{$Y$-set parameters (subsets and maximal sets)}\label{s:maxYparam}

The original $X$-set parameters were strongly motivated by parameters related to standard zero forcing  and domination. As such, the abstraction of these parameters naturally led to the Superset axiom. Another natural family of cohesive parameters is obtained by replacing the Superset axiom with a ``Subset axiom.'' 
Many of the results for super $X$-set parameters, including connectedness results, can be naturally adapted to sub $Y$-set parameters. The main isomorphism theorem, Theorem \ref{t:main},  is extended to certain sub $Y$-set parameters through  a complementation technique, described in Section \ref{ss:comp-param}.  This allows us to apply results to additional parameters such as independence number, irredundance, and zero forcing irredundance. 

\begin{defn}
 A \emph{sub $Y$-set 
parameter}  is a cohesive parameter $Y$ such that $Y(G)$ is defined to be the maximum cardinality of a $Y$-set of $G$  where   the $Y$-sets of $G$ satisfy the following condition:  
\ben[(I)]
\item (Subset) If $T$ is a $Y$-set of $G$ and $T'\subseteq T$, then $T'$ is a $Y$-set of $G$.
\een

\begin{rem} 
Let $Y$ be a sub $Y$-set parameter and let $G$ be a  graph. Recall that an \emph{abstract simplicial complex} $\mathcal{C}$ is a collection of finite sets such that if $A\in \mathcal{C}$ and $B\subseteq A$, then $B\in\mathcal{C}$ (see, for example, \cite[Definition 4.1.2]{Hal22}).  Furthermore, the set of \emph{vertices of $\mathcal{C}$} is $V(\mathcal{C})=\cup_{A\in \mathcal{C}}A$. In the setting of matroids an abstract simplicial complex is referred to as an independence system. The $Y$-sets of $G$ form an abstract simplicial complex.    Note that $V(\mathcal{C})\subseteq V(G)$ but these sets need not be equal (See Example \ref{ex:FZ-path}).
\end{rem}

When $Y$ is a sub $Y$-set parameter, the {\em $Y$-TAR graph}  of a base graph  $G$ is denoted by  $\ytar(G)$  and the \emph{lower $Y$ number}, denoted by  $\ybar(G)$, is the minimum cardinality of a maximal $Y$-set.
\end{defn}

\subsection{Connectedness}\label{ss:Yconnect}

 The study of $k$-TAR reconfiguration graphs in the setting of sub $Y$-set parameters requires minor modifications to Definition \ref{def: k token X}.
 
\begin{defn}  Suppose $Y$ is a  sub $Y$-set parameter.      The \emph{$k$-token addition and removal (TAR) reconfiguration graph for $Y$}, denoted by $\ytar_k(G)$, is the subgraph of $\ytar(G)$ induced by the set of all  $Y$-sets of cardinality at least $k$. 

The greatest $k_0$ such that $\ytar_k(G)$ is connected for all $k\le k_0$ is denoted by $\yyo(G)$, and the greatest $k$ such that $\ytar_k(G)$ is connected is denoted by $\ulyo(G)$.
\end{defn}

 The next result is the $Y$-set parameter version of Proposition \ref{basic1}. The proofs of \eqref{h1y}, \eqref{h2y}, and \eqref{h4y} are analogous; we prove \eqref{h3y} here since it needs more significant adjustment.

\begin{prop}\label{Ybasic1}
    Let  $Y$ be a  sub $Y$-set parameter and let $G$ be a graph of order $n$.
    \ben[$(1)$]
    \item\label{h1y} Then $Y(G)\ge \ulyo(G)\ge\yyo(G)$.
    \item\label{h2y} If $G$ has only one maximal $Y$-set, then $Y(G)=\ulyo(G)=\yyo(G)$.
    \item\label{h3y} If $G$ has more than one  maximal $Y$-set, then 
    $\ybar(G)-1 \ge \yyo(G)\ge\max\{\ybar(G)+Y(G)-n,0\}$.
    \item\label{h4y} If $G$ has more than one  maximum $Y$-set, then 
    $Y(G)-1\ge \ulyo(G)$.
    \een
\end{prop}

\bpf
{\blu\eqref{h1y}: This relationship follows immediately from the definitions.

\eqref{h2y}: Taking away one vertex at a time from one maximal $Y$-set does not disconnect the graph, therefore $\ytar_k(G)$ is always connected for $k \leq Y(G)$ when there is only one maximal $Y$-set.}

\eqref{h3y}: Let $\hat T\subseteq V(G)$ be a maximal $Y$-set with  $|\hat T|=\ybar(G)$.  Then $\hat T$ is an isolated vertex of $\ytar_{\ybar(G)}(G)$ (because we can't remove a vertex, and adding a vertex results in a set that is not a $Y$-set).  Thus $\ybar(G)-1\le \yyo(G)$.
  It is immediate from the definition of $\yyo(G)$ that $\yyo(G)\ge 0$. Suppose that  $\ybar(G)+Y(G)- n\ge 0$ and let $k_0=\ybar(G)+Y(G)-n$. Let $ T\subseteq V(G)$ be a maximal $Y$-set of $G$ and let $ T'\subseteq V(G)$ be a maximum $Y$-set of $G$.    To ensure $\ytar_k(G)$ is connected for all $k\le k_0$, it is sufficient to show that every such pair of vertices $ T$ and $T'$ is connected in $\ytar_{k_0}(G)$.
Define $T''=T\cap T'$ and observe that $|T''|\ge k_o$.  Then each of $T$ and $T'$ is connected by a path to $T''$ by removing one vertex at a time.  Thus $\yyo(G)\ge\ybar(G)+Y(G)-n$. 

{\blu \eqref{h4y}: Each minimum $Y$-set is an isolated vertex in $\ytar_{Y(G)}(G)$.}
\epf

  It is often easy to find examples of graphs for which $\ybar(G)-1= \yyo(G)$ and $\ulyo(G)=\yyo(G)$.    The next result, which is the $Y$-set parameter version of Corollary \ref{c:X1}, provides some such examples.

\begin{cor}\label{c:Y1}
Let  $Y$ be a  sub  $Y$-set parameter.    
If $G$  has more than one maximal $Y$-set and $Y(G)=n-1$, then $\yyo(G)=\ybar(G)-1$.
\end{cor}

The next result is the $Y$-set parameter version of Proposition \ref{basic2}; the proof is analogous.
 
\begin{prop}\label{Ybasic2}
    Let  $Y$ be a  sub $X$-set parameter and let $G$ be a graph of order $n$. 
     \ben[$(1)$]
\item\label{Yconn1} If for every pair of maximal $X$-sets $M_1$ and $M_2$, there is a path between $M_1$ and $M_2$ in $\ytar_k(G)$, then $\ytar_k(G)$ is connected. 

  \item\label{Yconn3}  If $k\le \ybar(G)$ and $\ytar_k(G)$ is connected, then $\yyo(G)\ge k$.
\item If  $|M_1\cap M_2|\ge k$ for every pair of maximal $Y$-sets $M_1$ and $M_2$, then $\yyo(G)\ge k$.  \een
\end{prop}
{\blu \bpf  Suppose first that  for every pair of maximal $X$-sets $M_1$ and $M_2$, there is a path between $M_1$ and $M_2$ in $\ytar_k(G)$.   Given two $Y$-sets $T_1,T_2\in V(\ytar_k(G))$, each $T_i$  is contained in a maximal $Y$-set $M_i$. There are paths in $\ytar_k(G)$ from $T_1$ to $M_1$, $M_1$ to $M_2$, and $M_2$ to $T_2$, so $\ytar_k(G)$ is connected.

Now assume    $k\le \ybar(G)$, $\ytar_k(G)$ is connected, and $\ell<k$. Since $k\ge \ybar(G)$, $\ytar_k(G)$ contains every maximal $Y$-set and there is a path between every pair of maximal $Y$-sets in $\ytar_k(G)$, which is a subgraph of $\ytar_\ell(G)$.  Thus $\ytar_\ell(G)$ is connected by \eqref{Yconn1}. 

If   $|M_1\cap M_2|\ge k$ for every pair of minimal $X$-sets $M_1$ and $M_2$, then there  is a path through $M_1\cap M_2$ in $\ytar_k(G)$ for every pair of maximal $Y$-sets $M_1$ and $M_2$ and $\yyo(G)\ge k$.
\epf

}

The next result extends Theorem 7 in \cite{KMM12} from independent sets to $Y$-sets, and the same proof remains valid (the proof of Proposition \ref{p:vsetTJ} is essentially a complemententary version  of that proof in \cite{KMM12}).

\begin{prop} \label{p:YvsetTJ}
Let $Y$ be a sub $Y$-set parameter, let $G$ be a graph, and let  $T_1$ and $T_2$ be $Y$-sets of $G$  with $|T_i|=k, i=1,2$. Then there is a path between  $T_1$ and $T_2$ in the $k$-TJ reconfiguration graph of $G$  if and only if there is a path between  $T_1$ and $T_2$ in $\ytar_{k-1}(G)$.     
\end{prop}

\subsection{Base graph complement parameters}\label{ss:comp-param}
A key to studying sub $Y$-set parameters is observing the following correspondence between sub $Y$-sets and super $X$-sets obtained by taking set complements in the base graph.

\begin{defn} (Base graph complement parameters)
Let $Y$ be a sub $Y$-set parameter and let $G$ be a graph.  We say that a subset $S\subseteq V(G)$ is an $X_Y$-set if and only if $V(G)\setminus S$ is a $Y$-set of $G$, and we define the parameter $X_Y(G)$ to be the minimum cardinality of an $X_Y$-set of $G$.

Let $X$ be a super $X$-set parameter.  We say that a subset $T\subseteq V(G)$ is a $Y_X$-set if and only if $V(G)\setminus T$ is an $X$-set of $G$, and we define the parameter $Y_X(G)$ to be the maximum cardinality of a $Y_X$-set of $G$.

The parameters $Y$ and $X_Y$ (or $X$ and $Y_X$) are called \emph{base graph complement  parameters}. The  $X_Y$-TAR  and $Y_X$-TAR graphs of a base graph $G$ are denoted by
$\xytar(G)$ and  $\yxtar(G)$ 
\end{defn}

The following definition introduces a notion of robustness for sub $Y$-set parameters that is analogous to super $X$-set parameters  via base graph complement parameters. A  set of one vertex is called a \emph{singleton}. 

\begin{defn}\label{Robust-parameters-Y}
A \emph{robust  $Y$-set parameter} is a sub $Y$-set parameter $Y$  such that $Y(G)$ and  the $Y$-sets of $G$ satisfy the following conditions:  
\ben[(I)]
\item\label{y:sub} (Subset)
If $T$ is a $Y$-set of $G$ and $T'\subseteq T$, then $T'$ is a $Y$-set of $G$.
\item\label{y:n-1} 
(Singletons)
If $G$ is a connected graph of order  at least two, 
then every singleton is a $Y$-set. 
\item \label{y:component}
(Component consistency) 
Let $G$ be a graph with connected components $G_1,\ldots,G_k$.  Then $T$ is a $Y$-set of $G$ if and only if $T\cap V(G_i)$ is a $Y$-set of $G_i$ for $i=1,\dots,k$.
\een  
\end{defn}

\begin{rem}\label{p:comp-params} The following statements are  immediate from the definitions:   If  $Y$ is a sub (respectively, robust) $Y$-set parameter, then $X_Y$ is a super (respectively, robust) $X$-set parameter.
If $X$ is a super (respectively, robust) $X$-set parameter, then $Y_X$ is a sub (respectively, robust) $Y$-set parameter.
The $Y_{X_Y}$-sets of $G$ are the $Y$-sets of $G$ and $Y_{X_Y}(G)=Y(G)$ and analogous statements hold when $Y$ and $X$ are interchanged. 
\end{rem}

 The definition of base graph complement  parameters  is illustrated in the next example.

\begin{ex}\label{ex:star-Ztar-Y_Ztar}
  Label the vertices of the star $K_{1,3}$ with  $\{1,2,3,4\}$ where $1$ is the vertex of degree three. The minimal standard zero forcing sets of $K_{1,3}$ are  $\{2,3\}, \{2,4\},$ and $\{3,4\}$ and $\ztar(K_{1,3})$ is shown in Figure \ref{f:star-comp-Z}(a) with the part of the hypercube of all subsets that are not $\Z$-sets shown in gray.  The maximal $Y_{\Z}$-sets of $K_{1,3}$ are $\{1,4\}, \{1,3\},$ and $\{1,2\}$ and $\ytar_{\Z}(K_{1,3})$ is shown in Figure \ref{f:star-comp-Z}(b) with the part of the hypercube of all subsets that are not $Y_{\Z}$-sets shown in gray. 
\end{ex}

\begin{figure}[!h]
\centering
\scalebox{0.7}{
\begin{tikzpicture}[scale=1.75,every node/.style={draw,shape=circle,outer sep=2pt,inner sep=1pt,minimum size=.2cm}]		
\node[fill=black, label={[yshift=-20pt]$\{1,2,3,4\}$}]  (1) at (0,2) {};
\node[fill=black, label={[yshift=-15pt, xshift=-10pt]$\{1,2,3\}$}]  (2) at (-1.5,1) {};
\node[fill=black, label={[yshift=-23pt, xshift=-20pt]$\{1,2,4\}$}]  (3) at (-0.5,1) {};
\node[fill=black, label={[yshift=-23pt, xshift=20pt]$\{1,3,4\}$}]  (4) at (0.5,1) {};
\node[fill=black, label={[yshift=-15pt, xshift=10pt]$\{2,3,4\}$}]  (5) at (1.5,1) {};
\node[fill=none, opacity = 0.3, label={[yshift=-10pt, xshift=-6pt, opacity = 0.3]$\{1,2\}$}]  (6) at (-2.5,0) {};
\node[fill=none, opacity = 0.3, label={[yshift=-15pt, xshift=-15pt, opacity = 0.3]$\{1,3\}$}]  (7) at (-1.5,0) {};
\node[fill=black, label={[yshift=-20pt, xshift=16pt]$\{2,3\}$}]  (9) at (0.5,0) {};
\node[fill=none, opacity = 0.3, label={[yshift=-20pt, xshift=-16pt, opacity = 0.3]$\{1,4\}$}]  (8) at (-0.5,0) {};
\node[fill=black, label={[yshift=-15pt, xshift=15pt]$\{2,4\}$}]  (10) at (1.5,0) {};
\node[fill=black, label={[yshift=-10pt, xshift=6pt]$\{3,4\}$}]  (11) at (2.5,0) {};
\node[fill=none, opacity = 0.3, label={[yshift=-28pt, xshift=-5pt, opacity = 0.3]$\{1\}$}]  (12) at (-1.5,-1) {};
\node[fill=none, opacity = 0.3, label={[yshift=-28pt, xshift=-5pt, opacity = 0.3]$\{2\}$}]  (13) at (-0.5,-1) {};
\node[fill=none, opacity = 0.3, label={[yshift=-28pt, xshift=5pt, opacity = 0.3]$\{3\}$}]  (14) at (0.5,-1) {};
\node[fill=none, opacity = 0.3, label={[yshift=-28pt, xshift=5pt, opacity = 0.3]$\{4\}$}]  (15) at (1.5,-1) {};
\node[fill=none, opacity = 0.3, label={[yshift=-25pt, opacity = 0.3]$\emptyset$}, label={[yshift=-60pt]\Large(a)}]  (16) at (0,-2) {};
		
\draw[thick] 
(1)--(2)--(9)--(5)--(1)--(3)--(10)--(5)--(11)--(4)--(1);
\draw[thick, opacity=0.3]
(1)--(2)--(6)--(12)--(16)--(15)--(11)--(5)--(1)--(3)--(6)--(13)--(16)--(14)--(11)--(4)--(8)--(15)--(10)--(5)--(9)--(2)--(7)--(12)--(8)--(3)--(10)--(13)--(9)--(14)--(7)--(4)--(1);

\node[fill=black, label={[yshift=-2pt]$\emptyset$}]  (1) at (6,2) {};
\node[fill=black, label={[yshift=-7pt, xshift=-5pt]$\{4\}$}]  (2) at (4.5,1) {};
\node[fill=black, label={[yshift=-7pt, xshift=-5pt]$\{3\}$}]  (3) at (5.5,1) {};
\node[fill=black, label={[yshift=-7pt, xshift=5pt]$\{2\}$}]  (4) at (6.5,1) {};
\node[fill=black, label={[yshift=-7pt, xshift=5pt]$\{1\}$}]  (5) at (7.5,1) {};
\node[fill=none, opacity = 0.3, label={[yshift=-10pt, xshift=-6pt, opacity = 0.3]$\{3,4\}$}]  (6) at (3.5,0) {};
\node[fill=none, opacity = 0.3, label={[yshift=-15pt, xshift=-15pt, opacity = 0.3]$\{2,4\}$}]  (7) at (4.5,0) {};
\node[fill=none, opacity = 0.3, label={[yshift=-20pt, xshift=-16pt, opacity = 0.3]$\{2,3\}$}]  (9) at (5.5,0) {};
\node[fill=black, label={[yshift=-20pt, xshift=16pt]$\{1,4\}$}]  (8) at (6.5,0) {};
\node[fill=black, label={[yshift=-15pt, xshift=15pt]$\{1,3\}$}]  (10) at (7.5,0) {};
\node[fill=black, label={[yshift=-10pt, xshift=6pt]$\{1,2\}$}]  (11) at (8.5,0) {};
\node[fill=none, opacity = 0.3, label={[yshift=-37pt, xshift=-10pt, opacity = 0.3]$\{2,3,4\}$}]  (12) at (4.5,-1) {};
\node[fill=none, opacity = 0.3, label={[yshift=-30pt, xshift=-20pt, opacity = 0.3]$\{1,3,4\}$}]  (13) at (5.5,-1) {};
\node[fill=none, opacity = 0.3, label={[yshift=-30pt, xshift=20pt, opacity = 0.3]$\{1,2,4\}$}]  (14) at (6.5,-1) {};
\node[fill=none, opacity = 0.3, label={[yshift=-37pt, xshift=10pt, opacity = 0.3]$\{1,2,3\}$}]  (15) at (7.5,-1) {};
\node[fill=none, opacity = 0.3, label={[yshift=-40pt, opacity = 0.3]$\{1,2,3,4\}$}, label={[yshift=-60pt]\Large(b)}]  (16) at (6,-2) {};
		
\draw[thick] 
(1)--(2)--(8)--(5)--(1)--(3)--(10)--(5)--(11)--(4)--(1);

\draw[thick, opacity = 0.3] (1)--(2)--(6)--(12)--(16)--(15)--(11)--(5)--(1)--(3)--(6)--(13)--(16)--(14)--(11)--(4)--(9)--(15)--(10)--(5)--(8)--(2)--(7)--(12)--(9)--(3)--(10)--(13)--(8)--(14)--(7)--(4)--(1);
\end{tikzpicture}}
\caption{The graphs $\ztar(K_{1,3})$ and $\ytar_{\Z}(K_{1,3})$, illustrating the isomorphism}
\label{f:star-comp-Z}
\end{figure}

\begin{rem}\label{rem: isom compl param}
    Let $Y$ be a sub $Y$-set parameter.  For every graph $G$,   $\ytar(G)\cong \xytar(G)$ by the isomorphism $T\to V(G)\setminus T$.  This isomorphism is illustrated in Figure \ref{f:star-comp-Z} by inverting the hypercube of all subsets of the base graph.
\end{rem}

Next we 
 translate  the main  isomorphism result to robust $Y$-set parameters.  
 
\begin{thm}\label{t:main-comp} 
Let $Y$ be a robust $Y$-set parameter and let $G$ and $G'$ be base graphs such $\ytar(G)\cong\ytar(G')$. If $Y(K_1)=1$ or $G$ and $G'$ have no isolated vertices, then $G$ and $G'$ have the same order and there is a relabeling of the vertices of $G'$ such that $G$ and $G'$ have exactly the same $Y$-sets.  
\end{thm}
\bpf Note that $\ytar(G)\cong \xytar(G))$ and $\ytar(G')\cong \xytar(G'))$ by the isomorphism $T\to V(G)\setminus T$, which implies $\xytar(G)\cong \xytar(G')$.
Since
$Y(K_1)=1$ or $G,G'$ have no isolated vertices, we see that
$X_Y(K_1)=0$ or $G,G'$ have no isolated vertices.
Since $Y$ is a robust $Y$-set parameter, $X_Y$ is a robust $X$-set parameter.  Thus $G$ and $G'$ have the same order and there is a relabeling of the vertices of $G'$ such that $G$ and $G'$ have exactly the same $X_Y$-sets by  Theorem \ref{t:main}.  Therefore, $G$ and $G'$ have exactly the same $Y$-sets.
\epf

Recall that for a super $X$-set parameter, a vertex $v$ of $G$ is $X$-irrelevant if $v\not\in S $ for every minimal $X$-set $S$ of $G$. The next definition translates irrelevance  from super $X$-set parameters to sub $Y$-set parameters using complementary parameters.
 
 \begin{defn}
     Let $G$ be a graph and let $Y$ be a sub $Y$-set parameter.  A vertex  $v\in  V(G)$ is \emph{$Y$-irrelevant}  if $v\in T $ for every maximal $Y$-set $T$ of $G$.
 \end{defn}

\begin{obs}
Let $G$ be a graph and let $Y$ be a sub $Y$-set paramter. Then $v \in V(G)$ is $Y$-irrelevant if and only if $v$ is $X_Y$-irrelevant. Similarly for a super $X$-set parameter $X$, $v \in V(G)$ is $X$-irrelevant if and only if $v$ is $Y_X$-irrelevant.
\end{obs}

 \begin{rem}
A connected $Y$-set parameter definition analogous to Definition \ref{Connected-parameters} could also be given, allowing the application of these results to connected graphs when the $Y$-set parameter lacks component consistency.    
\end{rem}

\subsection{Independent Sets}\label{ss:independent}

In this section  we  define the independence TAR graph, list some known results for such graphs, and apply results of about robust $Y$-set parameters.  
As with the other parameters, our focus here is on isomorphisms of TAR graphs, connectedness, existence of Hamilton cycles or paths, etc.  There 
is extensive prior work on reconfiguration of  independent sets, including TAR reconfiguration (see, for example \cite{IDHPSUU11,KMM12});  most of this work is focused on complexity, approximation, and algorithms.  However, we point out that some of these papers provide useful tools for the type of results we emphasize. In particular,  \cite{KMM12} contains structural results, including one that  implies  equality in a connectedness bound   for certain families of graphs   (see Proposition \ref{p:ind-evenholefree}).
  
   Observe that the independence number  $\alpha$  
   is a robust $Y$-set parameter.   Furthermore, the independence number is the base graph complement  parameter of the vertex cover number.
  Denote the \emph{independence TAR graph}   of a  base graph  $G$ by  $\itar(G)$. 
 All the properties of a robust $Y$-set parameter apply to independence number.   Here we state only the main theorem (noting that $\alpha(K_1)=1$).

 \begin{thm}\label{t:main-indep}  Let $G$ and $G'$ be graphs. If $\itar(G)\cong\itar(G')$, then  $G$ and $G'$ have the same order and  there is a relabeling of the vertices of $G'$ such that $G$ and $G'$ have exactly the same independent sets.  \end{thm}
 
Note that a vertex is $\alpha$-irrelevant if and only if it is an isolated vertex. Since independence number is complementary to vertex cover number, we have the next result.

 \begin{prop}\label{p:Iunique} 
Every independence TAR graph is unique.     \end{prop}

Hamiltonicity and connectedness results for independent set reconfiguration also parallel those for vertex cover reconfiguration. 
 Let $\ibar$ denote the minimum cardinality of a maximal independent set.  
  The greatest $k_0$ such that $\itar_k(G)$ is connected for all $k\le k_0$ is denoted by $\io(G)$, and the greatest $k$ such that $\itar_k(G)$ is connected is denoted by $\ulio(G)$.  
The next example shows that strict inequality is possible in the bound $\ibar(G)-1\ge \alpha_0$ (cf. Proposition \ref{Ybasic1}). It could be derived from Example \ref{ex:vc-Kpq} for vertex covering, but is also easy to see directly.

\begin{ex}
    The bipartite graph $K_{p,q}$ for $2 \leq p \leq q$ with partite sets $A$ and $B$ has  exactly $A$ and $B$ as its maximal independent sets. Then $\alpha(K_{p,q})=q$, $\ibar(K_{p,q})=p$, and the vertex of $\itar(K_{p,q})$ corresponding to the empty set is a cut-vertex. Thus $\alpha_0(K_{p,q})=0 <p-1$; in fact, $\alpha_0(K_{p,q})=\alpha(K_{p,q})+\ibar(K_{p,q})-n$.  Furthermore, $K_{p,q}$ has a  Hamilton path 
    but not a Hamilton cycle. 
\end{ex}

\begin{rem} Let $G$ be  a graph. Since the independence and vertex cover numbers are complementary parameters, $\itar(G)$ has a Hamiltonian path (cycle) if and only if $\vctar(G)$ has a Hamiltonian path (cycle) by Remark \ref{rem: isom compl param}. Examples illustrating equality or strict inequality in a connectedness bound for vertex covering illustrate the same for independence. See Section \ref{ss:vcex} for more examples.
\end{rem}

 Recall that an even-hole-free graph is one with no induced even cycles. 
 Proposition \ref{p:ind-evenholefree}, which shows many graphs $G$ satisfy $\alpha_0(G)=\ibar(G)-1$,  follows immediately from Proposition \ref{vctau0ehfg}, but to obtain it that way would distort the historical record and deny the authors of the original work on independent sets their due. Instead, we present the original independence results that are used in all the reconfiguration results on even-hole-free graphs discussed in this paper.  

  Theorem \ref{t:ind-evenholefree} is based on Theorem 7 in \cite{KMM12}, which is stated   for token jumping (TJ) reconfiguration. The translation between token jumping and TAR reconfiguration for indepemdent sets is established in Theorem 1 of \cite{KMM12} (which is the basis for Proposition \ref{p:vsetTJ}).

\begin{thm}\label{t:ind-evenholefree} {\rm \cite{KMM12}}
    If  $T_1$ and $T_2$  are independent sets of size $k$ in a graph $G$ and $G[(T_1\setminus T_2)\cup (T_2\setminus T_1)]$ is  even-hole-free, then there is a path between independent sets $T_1$ and $T_2$ in $\itar_{k-1}(G)$.
\end{thm}

\begin{prop}\label{p:ind-evenholefree}
    Let $G$ be an even-hole-free graph. Then $\alpha_0(G)=\ibar(G)-1$.
\end{prop}

\bpf Let $\ell=\ibar(G)$ and let $T_1, T_2 \in V(\itar_{\ell-1}(G))$. Then there exist maximal independent sets $M_1$ and $M_2$ such that $T_i \subseteq M_i, i=1,2$. Without loss of generality, assume $|M_1|\ge |M_2|$ and let $k=|M_2|$; note that $k\ge \ell$.  If $|T_1|> k$, choose $M'_1\subseteq T_1$ such that $|M'_1|=k$; otherwise, 
Choose  $M'_1$ such that $ T_1\subseteq M'_1$ and $|M'_1|=k$. Then by Theorem \ref{t:ind-evenholefree}, there exists a path between $M'_1$ and $M_2$ in $\itar_{k-1}(G)$. 
As in  the proof of Proposition \ref{vctau0ehfg}, there is there is a path from $T_1$ to $T_2$ in $\itar_{\ell-1}(G)$. Thus $\itar_{\ell-1}(G)$ is connected, and $\alpha_0(G)=\ibar(G)-1$ by Propositions \ref{Ybasic2}\eqref{Yconn3} and \ref{Ybasic1}\eqref{h3y}.
\epf

\subsection{Irredundance and zero forcing irredundance}\label{ss:irred}
In this section we apply  base graph complement parameter results to  irredundance number and  (standard) zero forcing irredundance number, which are defined in Section \ref{ss:param}.  The irredundance TAR graph of a base graph $G$ is denoted by $\irtar(G)$ and  the zero forcing irredundance TAR graph of a base graph $G$ is denoted by $\zirtar(G)$. The lower irredundance number of $G$ is denoted by $\ir(G)$ and  the  lower zero forcing irredundance number of $G$ is denoted by $\zir(G)$ (these follow the literature rather than our underline convention).

\begin{prop}\label{p:irred-robust}
The upper irredundance number $\IR$ and the  upper zero forcing irredundance number $\ZIR$  are robust $Y$-set parameters.
\end{prop}
\bpf 
Let $G$ be a graph.

\eqref{y:sub}: Any subset of an Ir-set or ZIr-set of $G$ is an  Ir-set or ZIr-set of $G$, respectively. Thus $\IR$ and $\ZIR$ are sub $Y$-set parameters.

\eqref{y:n-1}:  For $v\in V(G)$,  $N[v]$  is a private neighborhood of $v$ relative to $\{v\}$, so $\{v\}$ is an Ir-set.  Since $V(G)$ is a fort, it is a private fort of any one vertex $v$, and $\{v\}$ is a ZIr-set.

\eqref{y:component}:     Now assume $G_1, \dots, G_k$ are the connected components of $G$. For any vertex $v\in V(G)$, there is some $i$ such that $v\in V(G_i)$ and thus $N_G[v]\subseteq V(G_i)$.  Let $F$ be a fort. Then $F\cap V(G_i)$ is a fort of $G_i$ if and only if $F\cap V(G_i)$ is nonempty.  Thus  $T$ is a Ir-set of $G$ (respectively, a ZIr-set of $G$) if and only if $T\cap V(G_i)$ is an Ir-set of $G$ (respectively, a ZIr-set of $G_i$) for $i=1,\ldots,k$.
\epf

 By Proposition \ref{p:irred-robust}, the main isomorphism results hold for $\irtar(G)$ and $\zirtar(G)$ (note that $\IR(K_1)=1$ and $\ZIR(K_1)=1$).
 \begin{thm}\label{t:main-irred}  Let $G$ and $G'$ be graphs. If $\irtar(G)\cong\irtar(G')$, then  $G$ and $G'$ have the same order and  there is a relabeling of the vertices of $G'$ such that $G$ and $G'$ have exactly the same Ir-sets.  \end{thm}

  \begin{thm}\label{t:main-zir}  Let $G$ and $G'$ be graphs. If $\zirtar(G)\cong\zirtar(G')$, then  $G$ and $G'$ have the same order and  there is a relabeling of the vertices of $G'$ such that $G$ and $G'$ have exactly the same ZIr-sets.  \end{thm}

    Note that $\IR$ and $\ZIR$ are \emph{not} the base graph complement parameters to domination and standard zero forcing, even though we have used base graph complementation to establish the previous two results.
Next we discuss some examples of irredundance  TAR graphs and zero forcing irredundance  TAR graphs, including graphs with unique  TAR graphs and nonunique  TAR graph.

\begin{rem} 
Let $G$ be a graph on $n$ vertices. It is immediate that $\ir(G) = \IR(G) = n$ if and only if $G = \ol{K_n}$.  
By \cite[Remark 5.1]{ZIR}, $\zir(G) = \ZIR(G) = n$ if and only if $G = \ol{K_n}$. 
Thus, $\irtar(\ol{K_n})$ and   $\zirtar(\ol{K_n})$  are unique by  Theorems \ref{t:main-irred} and \ref{t:main-zir}. Moreover, $\irtar(\ol{K_n})\cong\zirtar(\ol{K_n)} \cong Q_n$  and therefore  $\irtar(\ol{K_n})$ and   $\zirtar(\ol{K_n})$ have Hamiltonian cycles.    

Similarly, let $G$ be  a graph of order $n$ with no isolated vertices.  If  $u$ and $v$ are not adjacent, then $\{u,v\}$ is an Ir-set, so $\IR(G)=1$ if and only if $G\cong K_n$.  Thus, $\irtar({K_n})$ is unique by   Theorem \ref{t:main-irred}. Moreover, $\irtar({K_n)} \cong K_{1,n}$  and hence does not have a Hamiltonian path for $n\geq 3$.
By \cite[Remark 5.2]{ZIR}, $\zir(G) = \ZIR(G) = n-1$ if and only if $G = {K_n}$. Thus, $\zirtar({K_n})$ is unique by   Theorem \ref{t:main-zir}. Moreover, $\zirtar({K_n)} \cong Q_n-v$  and therefore has a Hamiltonian path but not a Hamiltonian cycle.
\end{rem}

\begin{prop}
Let $G$ be a graph and $n\geq 1$. If $\zirtar(G) \cong \zirtar(P_n)$, then $G\cong P_n$. If $\zirtar(G) \cong \zirtar(K_{1,n-1})$, then $G\cong K_{1,n-1}$.
\end{prop}
\begin{proof} Assume $\zirtar(G) \cong \zirtar(P_n)$ or $\zirtar(G) \cong \zirtar(K_{1,n-1})$.  In \cite{ZIR} it is shown that $\zir(G) = 1$ if and only if $G\cong P_n$ or $G\cong K_{1,n-1}$.
 By  Theorem \ref{t:main-irred}, $G\cong P_n$ or $G\cong K_{1,n-1}$. Note that if $n\leq 3$, then $P_n \cong K_{1,n-1}$, so it suffices to show that $\zirtar(P_n) \not\cong \zirtar(K_{1,n-1})$  for $n\ge 4$. 

Let $v_1$ and $v_n$ be the vertices of degree 1 in $P_n$. Then every fort of $P_n$ contains both $v_1$ and $v_n$. Thus, $\{v_1\}$ and $\{v_n\}$ are maximal ZIr-sets and hence have degree 1 in $\zirtar(P_n)$.

Let $v$ be the vertex of degree $n-1$ in $K_{1,n-1}$. The only fort that contains $v$ is $V(K_{1,n-1})$. Thus, $\{v\}$ is a maximal ZIr-set and has degree 1 in $\zirtar(K_{1,n-1})$. Note that every pair of vertices in $V(G)\setminus \{v\}$ is a fort of $K_{1,n-1}$. Thus  for $n\ge 4$, $\{v\}$ is the only maximal ZIr-set of size 1. It follows that  $\zirtar(P_n) \not\cong \zirtar(K_{1,n-1})$ for $n\ge 4$.
\end{proof}

The next example is motivated by \cite[Example 3.3]{XTARiso}.
\begin{ex}
Let $n\geq 4$. Then $\zirtar(C_n) = \zirtar(C_n + vw)$ for any two vertices $v$ and $w$ that are not adjacent in $C_n$. A set $S$ is standard zero forcing set if and only it it contains two vertices that are adjacent on the cycle; a set that does not contain two such vertices cannot perform any force. 
Thus, $C_n$ and $C_n + vw$ have the same forts and so $\zirtar(C_n) = \zirtar(C_n + vw)$.
\end{ex}

\section{Concluding remarks}\label{s:sum}

 In Section \ref{sec:proof-requirements} we have established that with  a few  exceptions,  the  results for original $X$-set parameters in \cite{PDrecon} and \cite{XTARiso} remain valid for robust $\Xr$-set parameters (with the same proofs). 
Connectedness results are central to the study of reconfiguration (transforming one solution to another with every intermediate intermediate step being a  solution).  Although connectedness examples tend to be very parameter specific, in  Section \ref{ss:Xconn} most connectedness results  established in \cite{PDrecon} for original $X$-set parameters are extended to  super $X$-set parameters. 

Results for specific parameters are surveyed and extended in Sections \ref{s:Dom-PD-Z}-\ref{s:vtxcover}.  For many of these parameters, there are open questions remaining to investigate.  For example, we do not have an example of a PSD TAR graph that has a Hamilton cycle.  

We introduce the universal perspective on TAR reconfiguration for maximizing parameters (sub $Y$-set parameters)  in Section \ref{s:maxYparam}, and many avenues remain unexplored.  What does simplical complex theory tell us about robust $Y$-set parameters?  We introduced a way of complementing a sub $Y$-set parameter to obtain a super $X$-set parameter and vice versa, called base graph complements.  There may be additional useful ways to connect sub $Y$-set parameters and super $X$-set parameters. There is more work to be done on TAR reconfiguration of specific sub $Y$-set parameters.

Complexity of reconfiguration problems has been studied extensively  (see, for example, \cite{IDHPSUU11, KMM12, MS20, Nish18})  but is not surveyed here. Most of this work focuses on the complexity of reconfiguration of specific parameters.  It would be interesting to develop a universal theory of reconfiguration complexity for $X$-set and $Y$-set parameters; we hope the work here will be useful for this.


\section*{Acknowledgements}

The research of all authors was partially supported by NSF grant  2331634.  The research of B.~Curtis  is also partially supported by NSF grant  1839918. This research began at the American Institute of Mathematics (AIM) and the authors thank AIM and NSF for their support.

\end{document}